\documentclass[a4paper, 11pt]{article}
\usepackage{amssymb, amsmath,latexsym,amsfonts,amsbsy, amsthm,mathtools,graphicx,CJKutf8,CJKnumb,CJKulem,color,geometry}
\usepackage{verbatim}
\usepackage{hyperref}
\usepackage{newtxtext, newtxmath}
\usepackage[misc]{ifsym}
\usepackage[marginal]{footmisc}
\usepackage[titletoc,title]{appendix}
\def\l{\left}
\def\r{\right}
\def\f{\frac}

\def\az{\alpha}
\def\lz{\lambda}
\def\Lz{\Lambda}
\def\az{\alpha}
\def\rz{\rho}

\def\ez{\epsilon}
\def\bz{\beta}
\def\dz{\delta}
\def\gz{\gamma}
\def\tz{\theta}
\def\sz{\sigma}
\def\hat{\widehat}

\def\beq{\begin{equation}}
\def\eeq{\end{equation}}
\def\be{\begin{equation*}}
\def\ee{\end{equation*}}
\def\beqn{\begin{eqnarray}}
\def\eeqn{\end{eqnarray}}
\def\ben{\begin{eqnarray*}}
\def\een{\end{eqnarray*}}

\theoremstyle{plain}
\theoremstyle{plain}\newtheorem{thm}{Theorem}[section]
\theoremstyle{plain}\newtheorem{prop}{Proposition}[section]
\theoremstyle{plain}
\theoremstyle{plain}\newtheorem{lem}{Lemma}[section]
\theoremstyle{plain}\newtheorem{rem}{Remark}[section]

\numberwithin{equation}{section}

\makeatletter
\def\thanks#1{\protected@xdef\@thanks{\@thanks
        \protect\footnotetext{#1}}}
\makeatother

\begin{document}

\title{Global stability of the Dirac-Klein-Gordon system in two and three space dimensions}
\author{Qian Zhang}

\date{}

\maketitle

\noindent {\bf{Abstract}}\ \ In this paper we study global nonlinear stability for the Dirac-Klein-Gordon system in two and three space dimensions for small and regular initial data. In the case of two space dimensions, we consider the Dirac-Klein-Gordon system with a massless Dirac field and a massive scalar field, and prove global existence, sharp time decay estimates and linear scattering for the solutions. In the case of three space dimensions, we consider the Dirac-Klein-Gordon system with a mass parameter in the Dirac equation, and prove uniform (in the mass parameter) global existence, unified time decay estimates and linear scattering in the top order energy space.

\bigskip


\bigskip


\section{Introduction}\label{s1}

Let $n=2$ or $3$. We consider the Dirac-Klein-Gordon system in $n$ space dimensions:
\beq\label{s2: nDKG}
\l\{\begin{split}
i\gz^\mu\partial_\mu\psi+M\psi&=F_\psi:=v\psi, \ \  \quad\quad (t,x)\in[2,\infty)\times\mathbb{R}^n,\\
-\Box v+v&=F_v:=\psi^*\gz^0\psi, \quad (t,x)\in[2,\infty)\times\mathbb{R}^n
\end{split}\r.
\eeq
with prescribed initial data at $t=2$
\beq\label{s2: ninitial}
(\psi,v,\partial_tv)|_{t=2}=(\psi_0,v_0,v_1).
\eeq
In the above, $i\gz^\mu\partial_\mu=i\gz^0\partial_t+i\gz^1\partial_1+\cdots+i\gz^n\partial_n$ is the Dirac operator, $\partial_a=\partial_{x_a}$ for $a\in\{1,\cdots,n\}$, $\psi=\psi(t,x):\mathbb{R}^{1+n}\to\mathbb{C}^{N_0}$ is a spinor field with mass $M\ge 0$, $v=v(t,x): \mathbb{R}^{1+n}\to\mathbb{R}$ is a scalar field, $\psi^*$ denotes the complex conjugate transpose of the vector $\psi$, and $\gz^\mu$ are the Dirac matrices, where $N_0:=2^{[(n+1)/2]}$ ($[s]$ is the largest integer which is smaller than or equal to $s$, for $s\in\mathbb{R}$). Dirac matrices are defined by the identities 
\beq\label{sn: 1.3}
\gz^\mu\gz^\nu+\gz^\nu\gz^\mu=-2g^{\mu\nu}{\boldsymbol{I}},\quad\quad (\gz^\mu)^*=-g_{\mu\nu}\gz^\nu,
\eeq
where $\mu,\nu\in\{0,1,\cdots,n\}$, $g=(g_{\mu\nu})=\mathrm{diag}(-1,1,\cdots,1)$ denotes the Minkowski metric in $\mathbb{R}^{1+n}$, $(g^{\mu\nu})$ is the inverse matrix of $(g_{\mu\nu})$, ${\boldsymbol{I}}:=I_{N_0}$ is the $N_0\times N_0$ identity matrix and $A^*$ is the  conjugate transpose of the matrix $A$. Throughout, Einstein summation convention is adopted. As usual, $\Box=g^{\mu\nu}\partial_\mu\partial_\nu=-\partial_t^2+\Delta$ denotes the wave operator. The Dirac-Klein-Gordon system is a basic model of proton-proton interaction or neutron-neutron interaction in physics. This system and the Maxwell-Dirac equations form the foundations of Relativistic Electrodynamics (see Bjorken and Drell \cite{Bjo}).

For $n=3$ and fixed Dirac mass parameter $M$, the global well-posedness and global behavior of the solutions to \eqref{s2: nDKG}-\eqref{s2: ninitial} have been widely studied both in the low regularity setting and in the high regularity setting; see, for example, Bournaveas \cite{Bo99}, D'Ancona, Foschi and   Selberg \cite{AFS}, Bejenaru and Herr \cite{BH17}, Bachelot \cite{Ba} and Choquet-Bruhat \cite{CB}. For $n=2$, Tsutsumi \cite{Ts03} and Delort, Fang and Xue \cite{DFX} proved small data global existence with optimal pointwise decay and scattering for the massive Dirac-Klein-Gordon system. On the other hand, for $(n,M)=(2,0)$, Bournaveas \cite{Bo} proved local existence for \eqref{s2: nDKG} with low regular initial data, and the global existence result was first proved by Gr\"unrock and Pecher \cite{GP}. Dong and Wyatt \cite{DW} first proved sharp asymptotic decay results for \eqref{s2: nDKG} with $(n,M)=(2,0)$ and small, regular initial data with compact support, using the method of hyperboloidal foliation of space-time. Recently Dong, Li, Ma and Yuan \cite{DLMY} removed the compactness condition on the support of the initial data in \cite{DW}. Precisely, by applying the ghost weight method introduced by Alinhac \cite{A}, they proved global existence, sharp time decay and linear scattering for \eqref{s2: nDKG} with $(n,M)=(2,0)$ and small, regular initial data. The hyperboloidal method used in \cite{DW} was introduced by Klainerman in \cite{K2} to prove global existence results for nonlinear Klein-Gordon equations by using commuting vector fields. Later this method was developed by Klainerman, Wang and Yang \cite{KWY, W} and by LeFloch and Ma \cite{LM14, LM18}; see also \cite{Dk, DL, DM}.

In this paper, we study the Cauchy problem \eqref{s2: nDKG}-\eqref{s2: ninitial} in two or three space dimensions. We consider the following two cases depending on the spatial dimension $n$ and the Dirac mass parameter $M$: $i)$ $(n,M)=(2,0)$; $ii)$ $n=3$ and $M\in[0,1]$. In both cases we consider small and regular initial data which decay at suitable rates at infinity. In each case we prove global existence with sharp time decay estimates and linear scattering in the high order energy space. In particular, in the case $i)$, we give an alternative proof of the results in \cite{DLMY}. In the case $ii)$, our global existence and pointwise decay results are uniform in the Dirac mass parameter $M\in[0,1]$ in the sense that, the smallness condition on the initial data is independent of $M$, and we give explicit dependence of the pointwise decay estimates on the parameter $M$. We use unified method to prove the global stability results in the cases $i)$-$ii)$. This method is inspired by a recent work by Huneau and Stingo \cite{HS}, where the authors studied global existence of solutions to a certain class of quasilinear systems of wave equations on a product space with null nonlinearities.

\vspace{0.5em}

\noindent${\mathbf{Major\ difficulties\ and\ key\ ideas.}}$ We set the initial data at time $t=2$ to simplify the computations and this is not essential by the invariance of \eqref{s2: nDKG} under time translations. We first consider \eqref{s2: nDKG}-\eqref{s2: ninitial} with $(n,M)=(2,0)$. The main difficulties in proving global existence of solutions in this case include: $a)$ the free Dirac field and  Klein-Gordon field decay at rates $t^{-\f{1}{2}}$ and $t^{-1}$ respectively, which means that neither of the nonlinearities in \eqref{s2: nDKG} is integrable in time; $b)$ there is no compactness assumption on the support of the initial data, and we cannot use the hyperboloidal method as in \cite{DW} directly. For this, we adopt an idea by Huneau and Stingo \cite{HS}, i.e., we divide the space-time into the exterior region $\mathscr{D}^{\rm{ex}}=\{(t,x): t\ge 2, |x|\ge t-1\}$ and the interior region $\mathscr{D}^{\rm{in}}=\{(t,x): t\ge 2, |x|<t-1\}$, and use weighted energy and Sobolev estimates on flat time slices in the exterior region, while using estimates on truncated hyperboloids in the interior region. However, there are several difficulties to overcome in using this method. Firstly, to close the bootstrap in the exterior region, we need to make use of the weighted space-time $L^2_{tx}$ estimates of the solution $(\psi,v)$. Considering the growth of the top order energy, it is essential to treat the energy of both $\psi$ and $v$ carefully. Secondly, to bound the lower order energy of $v$ on truncated exterior hyperboloids, we need the improved pointwise estimate of $[\psi]_-:=\psi-(x_a/r)\gz^0\gz^a\psi$. Usually this is obtained from the conformal energy estimate of $\Psi$ which is the solution to $\Box\Psi=F_\psi$. We can bypass the conformal energy and obtain pointwise estimates of $|L_0\Psi|$ and $|\Gamma\Psi|$ from weighted energy estimates of $\Psi$, where $L_0$ is the scaling vector field and $\Gamma\in\{\partial, L,\Omega\}$ (see Section \ref{sno}). Thirdly, due to the growth of the top order energy of the solution $(\psi,v)$ in the interior region, we need to bound the lower order energy of the solution by performing nonlinear transformations and give control of boundary integrals appearing in the energy estimates of the functions in these transformations. Furthermore, some new ingredients are required to prove the scattering results, since the linear scattering relies on flat $L^2$-type estimates. However, in the interior region, we only have $L^2$-type estimates on truncated hyperboloids. The main ideas are described below.

For any fixed $t\ge 2$, let $\Sigma^{\rm{ex}}_t:=\{x: |x|\ge t-1\}$ denote the truncated flat time slice in the exterior region and $\mathscr{l}_{[2,t]}:=\{(\tau,x): |x|=\tau-1, 2\le \tau\le t\}$ denote the portion of the boundary of the light cone in the time strip $[2,t]$. We introduce two weight functions $t^{-\dz}$ and $\omega(z):=(2+z)^{1+\lz}$ ($\dz\ge 0, \lz>0$ are constants), and define the weighted energy functionals for the Klein-Gordon field $v$ and the Dirac field $\psi$ as 
\be
E^{{\rm{ex}},\lz}_{1,\dz}(t,v):=\int_{\Sigma^{\rm{ex}}_t}t^{-\dz}\omega(r-t)\big(|\partial v|^2+v^2\big){\rm{d}}x,\quad\quad E^{{\rm{ex}},\lz}_{D}(t,\psi):=\int_{\Sigma^{\rm{ex}}_t}\omega(r-t)|\psi|^2{\rm{d}}x
\ee
respectively (see \eqref{dexfs} and  \eqref{dexfD}), where $r=|x|$. In deriving estimates of $E^{{\rm{ex}},\lz}_{1,\dz}(t,v)$ and $E^{{\rm{ex}},\lz}_{D}(t,\psi)$, we can also bound space-time $L^2_{tx}$ norms and integrals on the boundary $\mathscr{l}_{[2,t]}$, i.e., we obtain estimates of $\|v\|_{Y^{{\rm{ex}},\lz}_{1,\dz,t}}+R_1(t,v)$ and $\|\psi\|_{Y^{{\rm{ex}},\lz}_{D,t}}+R_D(t,\psi)$, where 
\beq\label{s1: exnorm}
\|v\|_{Y^{{\rm{ex}},\lz}_{1,\dz,t}}:=[E^{{\rm{ex}},\lz}_{1,\dz}(t,v)]^{1/2}+\l(\int_2^t\!\!\int_{\Sigma^{\rm{ex}}_\tau}\tau^{-\dz}\omega'(r-\tau)\big(|Gv|^2+v^2\big){\rm{d}}x{\rm{d}}\tau\r)^{1/2},
\eeq
\beq\label{s1: 'exnorm}
\|\psi\|_{Y^{{\rm{ex}},\lz}_{D,t}}:=[E^{{\rm{ex}},\lz}_D(t,\psi)]^{1/2}+\l(\int_2^t\int_{\Sigma^{\rm{ex}}_\tau}\omega'(r-\tau)|[\psi]_{-}|^2{\rm{d}}x{\rm{d}}\tau\r)^{1/2},
\eeq
\beq\label{s1: exnorm'}
R_1(t,v):=\l(\int_{\mathscr{l}_{[2,t]}}\!\!\tau^{-\dz}\big(|Gv|^2+v^2\big){\rm{d}}\sz\r)^{1/2},\quad\quad R_D(t,\psi):=\l(\int_{\mathscr{l}_{[2,t]}}|[\psi]_{-}|^2{\rm{d}}\sz\r)^{1/2},
\eeq
where $|Gv|=\l(\sum_{a=1}^2|G_av|^2\r)^{1/2}$, $G_a=\partial_a+(x_a/r)\partial_t, a\in\{1,2\}$ denote the good derivatives. The energy functionals on the truncated hyperboloids $\mathscr{H}^{\rm{ex}}_s$ and $\mathscr{H}^{\rm{in}}_s$ (see \eqref{s2: mscHex} and \eqref{s2: mscHin}) are controlled in terms of the energy on flat time slices in the exterior region above.

The bootstrap assumption in the exterior region is the bound for the following quantity
\beq\label{s1: bsN}
\sum_{|I|\le N}\big\{\|\hat{\Gamma}^I\psi\|_{Y^{{\rm{ex}},\lz}_{D,t}}+\|\Gamma^Iv\|_{Y^{{\rm{ex}},\lz}_{1,\dz,t}}\big\}+\sum_{|I|\le N-1}\|\Gamma^Iv\|_{Y^{{\rm{ex}},\lz}_{1,0,t}}, 
\eeq
where $N\ge 6$ is an integer, $\lz>0, 0<\dz\ll 1$, $I$ denotes a multi-index and $\hat{\Gamma}$ denotes the set of compatible vector fields. To close the estimate of the top order energy in \eqref{s1: bsN}, we use a special structure generated in the energy estimate of $\psi$ (see Proposition \ref{s2: exfsD}) and the estimates of the space-time $L^2_{tx}$ norms
\be\label{s1: L^2_tx}
\l(\int_2^t\big\|(2+r-\tau)^{\f{\lz}{2}}\big(\tau^{-\f{\dz}{2}}|\Gamma^Iv|+|[\hat{\Gamma}^I\psi]_-|\big)\big\|_{L^2_x(\Sigma^{\rm{ex}}_\tau)}^2{\rm{d}}\tau\r)^{1/2}, \quad|I|\le N
\ee
(see \eqref{s1: exnorm} and \eqref{s1: 'exnorm}) together with weighted Sobolev inequality on truncated flat time slices in the exterior region. For the lower order case, we conduct a nonlinear transformation $\tilde{v}=v-\psi^*\gz^0\psi$ (which was also used in \cite{DW, DLMY}). 

The refined exterior energy estimates of the solution obtained above imply corresponding estimates on truncated exterior hyperboloids. To derive (unweighted) lower order energy estimates of $v$ on exterior hyperboloids, we need the improved pointwise estimate of $[\psi]_-$, and this is obtained by the weighted energy estimates of the solution $\Psi$ to the equation $\Box\Psi=F_\psi$, and weighted Sobolev and Hardy inequalities on flat time slices in the exterior region.

The (unweighted) lower order energy estimates of the solution $(\psi,v)$ on truncated exterior hyperboloids, combined with the bootstrap estimates in the interior region, give pointwise decay estimates for the solution on truncated interior hyperboloids, via Sobolev inequality on hyperboloids. We also point out that in deriving exterior energy estimates of the solution $(\psi,v)$ and the functions in the nonlinear transformations, we also obtain refined estimates of integrals on the boundary of the light cone (see \eqref{s1: exnorm'}). These are used to close the bootstrap estimates in the interior region.

Next we consider \eqref{s2: nDKG}-\eqref{s2: ninitial} with $n=3$ and $M\in[0,1]$. In this case, the uniform (in the mass parameter) weighted space-time $L^2_{tx}$ estimates, and the faster pointwise decay of the solution in terms of these $L^2_{tx}$ estimates, give the boundedness of the top order energy of the solution $(\psi, v)$ on flat time slices in the exterior region, and therefore on truncated exterior hyperboloids. Hence we obtain global existence which is uniform in the Dirac mass parameter $M$. By exploiting the hidden structure of the Dirac equation, we obtain pointwise decay with explicit dependence on the mass parameter $M$.

Finally, we also prove linear scattering of the solutions in all the cases under consideration. These scattering results cannot be obtained from the pointwise decay estimates directly since the Sobolev norms of the nonlinearities are not integrable in time. We prove a technical lemma which gives a sufficient condition for the linear scattering of the solutions; see Lemma \ref{Sca}.

In the sequel, we use $C$ to denote a universal constant whose value may change from line to line. As usual, $A\lesssim B$ means that $A\le CB$ for some constant $C$. Given a vector or a scalar $w$, we use Japanese bracket to denote $\langle w\rangle:=(1+|w|^2)^{1/2}$.

Next we state our main results. 

We first present the results in two space dimensions. Consider 
\beq\label{s1: DKG}
\l\{\begin{split}
i\gz^\mu\partial_\mu\psi&=F_\psi:=v\psi, \ \   \quad\quad (t,x)\in[2,\infty)\times\mathbb{R}^2,\\
-\Box v+v&=F_v:=\psi^*\gz^0\psi, \quad (t,x)\in[2,\infty)\times\mathbb{R}^2
\end{split}\r.
\eeq
with prescribed initial data at $t=2$
\beq\label{s1: initial}
(\psi,v,\partial_tv)|_{t=2}=(\psi_0,v_0,v_1).
\eeq

\begin{thm}\label{thm1}
Let $N\in\mathbb{N}$ with $N\ge 6$. Then there exists a constant $\ez_0>0$ such that for any $0<\ez<\ez_0$ and all initial data $(\psi_0,v_0,v_1)$ satisfying the smallness condition 
\beq\label{s1: inismall}
\sum_{k=0}^{N+1}\|\langle |x|\rangle^{N+1}\nabla^kv_0\|_{L^2(\mathbb{R}^2)}+\sum_{k=0}^N\big\|\langle|x|\rangle^{k+1}|\nabla^k\psi_0|+\langle|x|\rangle^{N+1}|\nabla^kv_1|\big\|_{L^2(\mathbb{R}^2)}\le\ez,
\eeq
the Cauchy problem \eqref{s1: DKG}-\eqref{s1: initial} admits a unique global-in-time solution $(\psi,v)$, which enjoys the following sharp pointwise decay estimates
\ben
|\psi(t,x)|&\lesssim&\ez\langle t+|x|\rangle^{-\f{1}{2}}\langle t-|x|\rangle^{-\f{1}{2}},\\
|v(t,x)|&\lesssim&\ez\langle t+|x|\rangle^{-1}.
\een
In addition, the solution $(\psi,v)$ scatters to a free solution in $\mathcal{X}_{N-1}(\mathbb{R}^2):=H^{N-1}(\mathbb{R}^2)\times H^{N}(\mathbb{R}^2)\times H^{N-1}(\mathbb{R}^2)$, i.e., there exists $(\psi^*_{0},v^*_0,v^*_1)\in \mathcal{X}_{N-1}(\mathbb{R}^2)$ such that
\ben
\lim_{t\to+\infty}\|(\psi, v,\partial_tv)-(\psi^*, v^*,\partial_tv^*)\|_{\mathcal{X}_{N-1}(\mathbb{R}^2)}=0,
\een
where $(\psi^*,v^*)$ is the solution to 
\be
\l\{\begin{split}
i\gz^\mu\partial_\mu\psi^*&=0,\quad (t,x)\in[2,\infty)\times\mathbb{R}^2,\\
-\Box v^*+v^*&=0,\quad (t,x)\in[2,\infty)\times\mathbb{R}^2
\end{split}\r.
\ee
with the initial data $(\psi^*,v^*,\partial_tv^*)|_{t=2}=(\psi^*_{0},v^*_0,v^*_1)$.
\end{thm}

We next state the results in three space dimensions. Consider
\beq\label{s5: 3DDKG}
\l\{\begin{split}
i\gz^\mu\partial_\mu\psi+M\psi&=F_\psi:=v\psi, \ \  \quad\quad (t,x)\in[2,\infty)\times\mathbb{R}^3,\\
-\Box v+v&=F_v:=\psi^*\gz^0\psi, \quad (t,x)\in[2,\infty)\times\mathbb{R}^3
\end{split}\r.
\eeq
with prescribed initial data at $t=2$
\beq\label{s5: 3Dinitial}
(\psi,v,\partial_tv)|_{t=2}=(\psi_0,v_0,v_1).
\eeq

\begin{thm}\label{thm2}
Let $n=3$, $M\in[0,1]$ and $N\in\mathbb{N}$ with $N\ge 5$. Then there exists a constant $\ez_0>0$, which is independent of the mass parameter $M$, such that for any $0<\ez<\ez_0$ and all initial data $(\psi_0,v_0,v_1)$ satisfying the smallness condition 
\beq\label{s1: 3Dinismall}
\sum_{k=0}^{N+1}\|\langle |x|\rangle^{N+1}\nabla^kv_0\|_{L^2(\mathbb{R}^3)}+\sum_{k=0}^N\big\|\langle|x|\rangle^{N+1}\big(|\nabla^k\psi_0|+|\nabla^kv_1|\big)\big\|_{L^2(\mathbb{R}^3)}\le\ez,
\eeq
the Cauchy problem \eqref{s5: 3DDKG}-\eqref{s5: 3Dinitial} admits a unique global-in-time solution $(\psi,v)$, which enjoys the following pointwise decay estimates
\ben
|\psi(t,x)|&\lesssim&\f{\ez}{\langle t+|x|\rangle\langle t-|x|\rangle^{\f{1}{2}}+M^2\langle t+|x|\rangle^{\f{3}{2}}},\\
|v(t,x)|&\lesssim&\ez\langle t+|x|\rangle^{-\f{3}{2}}.
\een
In addition, the solution $(\psi,v)$ scatters to a free solution in $\mathcal{X}_N(\mathbb{R}^3):=H^N(\mathbb{R}^3)\times H^{N+1}(\mathbb{R}^3)\times H^{N}(\mathbb{R}^3)$, i.e., there exists $(\psi^*_{0},v^*_0,v^*_1)\in \mathcal{X}_N(\mathbb{R}^3)$ such that
\ben
\lim_{t\to+\infty}\|(\psi, v,\partial_tv)-(\psi^*, v^*,\partial_tv^*)\|_{\mathcal{X}_N(\mathbb{R}^3)}=0,
\een
where $(\psi^*,v^*)$ is the solution to
\be
\l\{\begin{split}
i\gz^\mu\partial_\mu\psi^*+M\psi^*&=0,\quad (t,x)\in[2,\infty)\times\mathbb{R}^3,\\
-\Box v^*+v^*&=0,\quad (t,x)\in[2,\infty)\times\mathbb{R}^3
\end{split}\r.
\ee
with the initial data $(\psi^*,v^*,\partial_tv^*)|_{t=2}=(\psi^*_{0},v^*_0,v^*_1)$.
\end{thm}

The organization of this paper is as follows. In Section \ref{sP}, we introduce some notations, present energy and Sobolev estimates both on flat time slices in the exterior region and on truncated hyperboloids. In Sections \ref{s3} - \ref{s4}, we prove the global existence result in Theorem \ref{thm1}. In Section \ref{s5}, we show the uniform global existence result in Theorem \ref{thm2}. Section \ref{s6} is devoted to proving the scattering results in Theorems \ref{thm1} and \ref{thm2}.

\section{Preliminaries}\label{sP}

In this section, we assume that the spatial dimension $n=2$ or $3$. Unless otherwise specified, the definitions and conclusions hold for both $n=2$ and $n=3$. We denote 
\beq\label{s2: defN_0}
N_0:=2^{[(n+1)/2]},
\eeq
i.e. $N_0=2$ for $n=2$ and $N_0=4$ for $n=3$, where $[s]$ denotes the largest integer which is smaller than or equal to $s$, for $s\in\mathbb{R}$. 

\subsection{Notations}\label{sno}
We work in the $(1+n)$ dimensional space-time $\mathbb{R}^{1+n}$ with Minkowski metric $g=(-1,1,\cdots,1)$, which is used to raise or lower indices. The space indices are denoted by Roman letters $a,b,\cdots\in\{1,2,\cdots,n\}$, and the space-time indices are denoted by Greek letters $\mu,\nu,\az,\bz,\cdots\in\{0,1,2,\cdots,n\}$. Einstein summation convention for repeated upper and lower indices is adopted throughout the paper. We denote a point in $\mathbb{R}^{1+n}$ by $(t,x)=(x_0,x_1,x_2,\cdots,x_n)$ with $t=x_0,x=(x_1,x_2,\cdots,x_n),x^a=x_a,a=1,2,\cdots,n$, and its spatial radius is denoted by $r:=|x|=\sqrt{x_1^2+x_2^2+\cdots+x_n^2}$. The following vector fields will be used frequently in the analysis (see Klainerman \cite{K86}):
\begin{itemize}
\item[(i)] Translations: $\partial_\az:=\partial_{x_\az}$, for $\az\in\{0,1,2,\cdots,n\}$.
\item[(ii)] Lorentz boosts: $L_a:=x_a\partial_t+t\partial_a$, for $a\in\{1,2,\cdots,n\}$.
\item[(iii)] Rotations: $\Omega_{ab}:=x_a\partial_b-x_b\partial_a$, for $1\le a<b\le n$.
\item[(iv)] Scaling: $L_0=t\partial_t+x^a\partial_a$.
\end{itemize}
We also use the modified Lorentz boosts and rotations introduced by Bachelot \cite{Ba},
\beq\label{s2: hatLO}
\hat{L}_a:=L_a-\f{1}{2}\gz^0\gz^a,\quad\quad\hat{\Omega}_{ab}:=\Omega_{ab}-\f{1}{2}\gz^a\gz^b,
\eeq
which enjoy the following commutative property, i.e.
\be
[\hat{L}_a,i\gz^\mu\partial_\mu]=[\hat{\Omega}_{ab},i\gz^\mu\partial_\mu]=0,
\ee
where the commutator $[A,B]$ is defined as 
$$[A,B]:=AB-BA.$$

We decompose the space-time $\mathbb{R}^{1+n}$ into two regions:
$$\mathrm{interior\ region\ }\quad \mathscr{D}^{\rm{in}}:=\{(t,x): t\ge 2, r<t-1\},$$
$$\mathrm{exterior\ region\ }\quad \mathscr{D}^{\rm{ex}}:=\{(t,x): t\ge 2, r\ge t-1\}.$$
We denote the constant time slices which foliate $\mathscr{D}^{\rm{ex}}$ as
$$\Sigma^{\rm{ex}}_t:=\{x\in\mathbb{R}^n: |x|\ge t-1\}.$$
The portion of exterior region in the time strip $[2,T]$ for any fixed time $T$ is denoted by
$$\mathscr{D}^{\rm{ex}}_T:=\{(t,x)\in\mathscr{D}^{\rm{ex}}: 2\le t\le T\}.$$
For $2\le T_1<T_2$, we denote the portion of the boundary of the light cone in the time interval $[T_1,T_2]$ by 
\beq\label{s2: l[2,t]}
\mathscr{l}_{[T_1,T_2]}:=\{(t,x): r=t-1, T_1\le t\le T_2\}.
\eeq
For any $s\ge 2$, we use
\beq\label{s2: mscH}
\mathscr{H}_s:=\{(t,x): t^2=s^2+|x|^2\}
\eeq
to denote the hyperboloid at hyperbolic time $s$. We denote by $\mathscr{H}^{\rm{in}}_s$ (resp. $\mathscr{H}^{\rm{ex}}_s$) the portion of $\mathscr{H}_s$ contained in the interior region $\mathscr{D}^{\rm{in}}$ (resp. in the exterior region $\mathscr{D}^{\rm{ex}}$), i.e.,
\beq\label{s2: mscHin}
\mathscr{H}^{\rm{in}}_s:=\{(t,x)\in\mathscr{H}_s: |x|<(s^2-1)/2\},
\eeq
\beq\label{s2: mscHex}
\mathscr{H}^{\rm{ex}}_s:=\{(t,x)\in\mathscr{H}_s: |x|\ge(s^2-1)/2\}.
\eeq
We set 
\beq\label{s2: tr(s)}
t(s):=\f{s^2+1}{2},\quad\quad r(s):=\f{s^2-1}{2}. 
\eeq
In addition, we denote by $\mathscr{H}^{\rm{in}}_{[2,s]}$ the hyperbolic interior region limited by the hyperboloids $\mathscr{H}_2$ and $\mathscr{H}_s$, and by $\mathscr{H}^{\rm{ex}}_{[2,s]}$ the portion of exterior region below $\mathscr{H}^{\rm{ex}}_{s}$, i.e.,
$$\mathscr{H}^{\rm{in}}_{[2,s]}:=\{(t,x)\in\mathscr{D}^{\rm{in}}: 2^2\le t^2-|x|^2\le s^2\},$$
$$\mathscr{H}^{\rm{ex}}_{[2,s]}:=\{(t,x)\in\mathscr{D}^{\rm{ex}}: t^2-|x|^2\le s^2\}.$$
For any sufficiently smooth function $g=g(t,x)$ defined on a hyperboloid $\mathscr{H}_s$, we denote 
$$\int_{\mathscr{H}_s}g{\rm{d}}x:=\int_{\mathbb{R}^n}g(\sqrt{s^2+|x|^2},x){\rm{d}}x,$$
$$\int_{\mathscr{H}^{\rm{in}}_s}g{\rm{d}}x:=\int_{|x|<r(s)}g(\sqrt{s^2+|x|^2},x){\rm{d}}x,\quad\quad \int_{\mathscr{H}^{\rm{ex}}_s}g{\rm{d}}x:=\int_{|x|\ge r(s)}g(\sqrt{s^2+|x|^2},x){\rm{d}}x,$$
where $r(s)=\f{s^2-1}{2}$ is as in \eqref{s2: tr(s)}. For any $1\le p<\infty$, we denote
\beq\label{s2: L^pH_s}
\|g\|_{L^p(\mathscr{H}_s)}:=\l(\int_{\mathscr{H}_s}|g|^p{\rm{d}}x\r)^{\f{1}{p}}=\l(\int_{\mathbb{R}^n}|g(\sqrt{s^2+|x|^2},x)|^p{\rm{d}}x\r)^{\f{1}{p}},
\eeq
$$\|g\|_{L^p(\mathscr{H}^{\rm{in}}_s)}:=\l(\int_{\mathscr{H}^{\rm{in}}_s}|g|^p{\rm{d}}x\r)^{\f{1}{p}},\quad\quad \|g\|_{L^p(\mathscr{H}^{\rm{ex}}_s)}:=\l(\int_{\mathscr{H}^{\rm{ex}}_s}|g|^p{\rm{d}}x\r)^{\f{1}{p}}.$$

For any sufficiently smooth function $u$, we denote 
\beq\label{s2: defNa}
|\partial u|=\big(|\partial_tu|^2+|\nabla u|^2\big)^{\f{1}{2}},\quad\quad|\nabla u|:=\l(\sum_{a=1}^n|\partial_au|^2\r)^{\f{1}{2}}.
\eeq
\beq\label{s2: defLO} 
|Lu|:=\l(\sum_{a=1}^n|L_au|^2\r)^{\f{1}{2}},\quad\quad |\Omega u|:=\l(\sum_{1\le a<b\le n}|\Omega_{ab}u|^2\r)^{\f{1}{2}}.
\eeq
We also denote 
\beq\label{s2: G_a}
G_a=\partial_a+\f{x_a}{r}\partial_t, \quad\quad |Gu|=\l(\sum_{a=1}^n|G_au|^2\r)^{\f{1}{2}}.
\eeq
and 
\beq\label{s2: barpar_a}
\bar{\partial}_a:=t^{-1}L_a=\partial_a+\f{x_a}{t}\partial_t.
\eeq
Let 
\beq\label{s2: defn_0}
n_0=\f{1}{2}(n^2+3n+2),
\eeq
i.e., $n_0=6$ for $n=2$ and $n_0=10$ for $n=3$. We define the ordered sets
\beq\label{s2: gamma_k}
\begin{split}
&\{\Gamma_k\}_{k=1}^{n_0}:=\big\{(\partial_\az)_{0\le\az\le n}, (L_a)_{1\le a\le n}, (\Omega_{ab})_{1\le a<b\le n}\big\}\\
&\{\hat{\Gamma}_k\}_{k=1}^{n_0}:=\big\{(\partial_\az)_{0\le\az\le n}, (\hat{L}_a)_{1\le a\le n},(\hat{\Omega}_{ab})_{1\le a<b\le n}\big\}.
\end{split}
\eeq
For any multi-index $I=(i_1,\cdots,i_{n_0})\in\mathbb{N}^{n_0}$ of length $|I|=\sum_{k=1}^{n_0}i_k$, we denote
\beq\label{s2: Ga^I}
\begin{split}
&\Gamma^I=\prod_{k=1}^{n_0}\Gamma_k^{i_k},\quad\mathrm{where}\quad\Gamma=(\Gamma_1,\dots,\Gamma_{n_0}),\\
&\hat{\Gamma}^I=\prod_{k=1}^{n_0}\hat{\Gamma}_k^{i_k},\quad\mathrm{where}\quad\hat{\Gamma}=(\hat{\Gamma}_1,\dots,\hat{\Gamma}_{n_0}).
\end{split}
\eeq
For any $I=(i_0,i_1,\cdots,i_n)\in\mathbb{N}^{1+n}$, $J=(j_1,\cdots,j_n)\in\mathbb{N}^n$, $K=(k^{ab})_{1\le a<b\le n}\in\mathbb{N}^{n_1}$, where $n_1:=(n^2-n)/2$, let 
\beq\label{s2: pOL^J}
\begin{split}
&\partial^I:=\prod_{\az=0}^{n}\partial_\az^{i_\az},\quad\quad L^J:=\prod_{a=1}^{n}L_a^{j_a},\quad\quad\Omega^K:=\prod_{1\le a<b\le n}\Omega_{ab}^{k^{ab}}\\
&\hat{L}^J:=\prod_{a=1}^{n}\hat{L}_a^{j_a},\quad\quad\hat{\Omega}^K:=\prod_{1\le a<b\le n}\hat{\Omega}_{ab}^{k^{ab}}.
\end{split}
\eeq

Given an index set $\Theta$, a finite set $\{P_{\!\tz}: \tz\in\Theta\}$ of linear operators and a linear operator $Q$, we write
\be
Q=\underset{\tz\in\Theta}{\sum\nolimits^{\prime}}P_{\!\tz}
\ee
if there exist some constants $c_{\!\tz}$ such that we have
\be
Q=\sum_{\tz\in\Theta}c_{\!\tz} P_{\!\tz}.
\ee
The following relations on commutators between the vector fields in $\{\Gamma_k\}_{k=1}^{n_0}$ are well-known (see for example \cite{Sog}):
\beqn\label{s2: commutator}
&&[\partial_\az, L_a]=\underset{0\le\bz\le n}{\sum\nolimits^{\prime}}\partial_\bz,\quad\quad [\partial_\az, \Omega_{ab}]=\underset{0\le\bz\le n}{\sum\nolimits^{\prime}}\partial_\bz,\quad\quad [\partial_\az, L_0]=\partial_\az,\quad\quad 0\le\az\le n, 1\le a<b\le n,\nonumber\\
&&[Z_k,Z_l]=\underset{Z_j\in V}{\sum\nolimits^{\prime}}Z_j,\quad\quad Z_k, Z_l\in V:=\{(L_a)_{1\le a\le n}, (\Omega_{ab})_{1\le a<b\le n}\},\nonumber\\
&&[L_0,\Omega_{ab}]=[L_0,L_c]=0,\quad\quad 1\le a<b\le n, 1\le c\le n.
\eeqn

The following proposition is a direct result from the definitions \eqref{s2: hatLO} and \eqref{s2: Ga^I}.

\begin{prop}\label{hatgz}
Let $N_0, n_0$ be as in \eqref{s2: defN_0} and \eqref{s2: defn_0} respectively. For any sufficiently smooth vector field $\phi=\phi(t,x):\mathbb{R}^{1+n}\to\mathbb{C}^{N_0}$, and any $I,J\in\mathbb{N}^{n_0}$, we have
\be
\hat{\Gamma}^I\phi=\Gamma^I\phi+\sum_{|I'|<|I|}c_{I'}\Gamma^{I'}\phi,\quad\quad\Gamma^J\phi=\hat{\Gamma}^J\phi+\sum_{|J'|<|J|}c_{J'}\hat{\Gamma}^{J'}\phi
\ee
for some $N_0\times N_0$ constant matrices $c_{I'}, c_{J'}$.
\end{prop}

\subsection{Energy estimates for wave and Klein-Gordon equations}

For any fixed constants $0\le m\le 1$, $\lz>0$ and $\dz\ge 0$, we define
\beq\label{dexfs}
E^{{\rm{ex}},\lz}_{m,\dz}(t,u):=\int_{\Sigma^{\rm{ex}}_t}t^{-\dz}(2+r-t)^{1+\lz}(|\partial u|^2+m^2u^2)(t,x){\rm{d}}x.
\eeq
When $\dz=0$, we denote for simplicity
\beq\label{dexfsaz} 
E^{{\rm{ex}},\lz}_m(t,u):=\int_{\Sigma^{\rm{ex}}_t}(2+r-t)^{1+\lz}(|\partial u|^2+m^2u^2)(t,x){\rm{d}}x.
\eeq
We also denote 
\beq\label{dexfsdz} 
E^{{\rm{ex}}}_{m,\dz}(t,u):=\int_{\Sigma^{\rm{ex}}_t}t^{-\dz}(|\partial u|^2+m^2u^2)(t,x){\rm{d}}x.
\eeq
In particular, we define the unweighted energy
\beq\label{dexfs0} 
E^{{\rm{ex}}}_m(t,u):=E^{{\rm{ex}}}_{m,0}(t,u)=\int_{\Sigma^{\rm{ex}}_t}(|\partial u|^2+m^2u^2)(t,x){\rm{d}}x.
\eeq

We have the following weighted energy estimate in the exterior region.

\begin{prop}\label{s2: exfs}
Let $\lz>0$ and $\dz\ge 0$. For any $t\ge 2$ and any sufficiently smooth function $u$ which is defined in the region $\mathscr{D}^{\rm{ex}}_t$ and decays sufficiently fast at space infinity, we have
\ben
\|u\|_{Y^{{\rm{ex}},\lz}_{m,\dz,t}}+R_m(t,u)\lesssim\|u\|_{Y^{{\rm{ex}},\lz}_{m,\dz,2}}+\int_2^t\|\tau^{-\f{\dz}{2}}(2+r-\tau)^{\f{1+\lz}{2}}(-\Box u+m^2u)\|_{L^2_x(\Sigma^{\rm{ex}}_\tau)}{\rm{d}}\tau,
\een
where 
\beq\label{s2: exnorm}
\begin{split}
\|u\|_{Y^{{\rm{ex}},\lz}_{m,\dz,t}}:&=[E^{{\rm{ex}},\lz}_{m,\dz}(t,u)]^{\f{1}{2}}+\l(\int_2^t\int_{\Sigma^{\rm{ex}}_\tau}\tau^{-\dz}(2+r-\tau)^\lz\big(|Gu|^2+m^2u^2\big){\rm{d}}x{\rm{d}}\tau\r)^{\f{1}{2}},\\
R_m(t,u):&=\l(\int_{\mathscr{l}_{[2,t]}}\tau^{-\dz}\big(|Gu|^2+m^2u^2\big){\rm{d}}\sz\r)^{\f{1}{2}},
\end{split}
\eeq
and we recall the definitions \eqref{s2: G_a} and \eqref{s2: l[2,t]}.

\begin{proof}
Let $\omega$ be a positive weight. Denote $F_u:=-\Box u+m^2u$. Multiplying on both sides of this equality by $t^{-\dz}\omega(r-t)\partial_t u$ and by direct calculation, we obtain 
\beqn\label{s2: omd_tubo}
\partial_t\big(t^{-\dz}\omega(|\partial u|^2+m^2u^2)\big)-2\partial_a\big(t^{-\dz}\omega\partial_tu\partial^au\big)&+&t^{-\dz}\omega'\big(|Gu|^2+m^2u^2\big)+\dz t^{-\dz-1}\omega(|\partial u|^2+m^2u^2)\nonumber\\
&=&2t^{-\dz}\omega\partial_tu F_u.
\eeqn
We note that the upward unit normal to $\mathscr{l}_{[2,t]}$ is $2^{-\f{1}{2}}(1,-x/r)$. Taking $\omega(z)=(2+z)^{1+\lz}$, and integrating \eqref{s2: omd_tubo} over the region $\mathscr{D}^{\rm{ex}}_{t}$, we have
\beqn\label{s2: exfse}
E^{{\rm{ex}},\lz}_{m,\dz}(t,u)&+&(1+\lz)\int_2^t\int_{\Sigma^{\rm{ex}}_\tau}\tau^{-\dz}(2+r-\tau)^\lz\big(|Gu|^2+m^2u^2\big){\rm{d}}x{\rm{d}}\tau+\int_{\mathscr{l}_{[2,t]}}\tau^{-\dz}\big(|Gu|^2+m^2u^2\big)2^{-\f{1}{2}}{\rm{d}}\sz\nonumber\\
&+&\dz\int_2^t\int_{\Sigma^{\rm{ex}}_\tau} \tau^{-\dz-1}(2+r-\tau)^{1+\lz}(|\partial u|^2+m^2u^2){\rm{d}}x{\rm{d}}\tau\nonumber\\
&=&E^{{\rm{ex}},\lz}_{m,\dz}(2,u)+2\int_2^t\int_{\Sigma^{\rm{ex}}_\tau}\tau^{-\dz}(2+r-\tau)^{1+\lz}\partial_tuF_u {\rm{d}}x{\rm{d}}\tau.
\eeqn
Differentiating \eqref{s2: exfse} with respect to $t$ and using H\"older inequality, we derive
\be
\begin{split}
\partial_t\Bigg\{E^{{\rm{ex}},\lz}_{m,\dz}(t,u)&+(1+\lz)\int_2^t\int_{\Sigma^{\rm{ex}}_\tau}\tau^{-\dz}(2+r-\tau)^\lz\big(|Gu|^2+m^2u^2\big){\rm{d}}x{\rm{d}}\tau+\int_{\mathscr{l}_{[2,t]}}\tau^{-\dz}\big(|Gu|^2+m^2u^2\big)2^{-\f{1}{2}}{\rm{d}}\sz\Bigg\}\\
&\le 2[E^{{\rm{ex}},\lz}_{m,\dz}(t,u)]^{\f{1}{2}}\cdot\|t^{-\f{\dz}{2}}(2+r-t)^{\f{1+\lz}{2}}F_u\|_{L^2_x(\Sigma^{\rm{ex}}_t)},
\end{split}
\ee
which implies the conclusion of the proposition.
\end{proof}
\end{prop}

\begin{rem}\label{RemK}
Taking $\omega\equiv 1,\dz=0$ in \eqref{s2: omd_tubo} and integrating over $[2,t]\times\mathbb{R}^n$, we obtain the standard energy inequality for wave and Klein-Gordon equations, i.e.,
\be
\|\big(|\partial u|+m|u|\big)(t)\|_{L^2_x(\mathbb{R}^n)}\lesssim\|\big(|\partial u|+m|u|\big)(2)\|_{L^2_x(\mathbb{R}^n)}+\int_2^t\|(-\Box u+m^2u)(\tau)\|_{L^2_x(\mathbb{R}^n)}{\rm{d}}\tau,
\ee
for any sufficiently smooth function $u$ which decays sufficiently fast at space infinity.
\end{rem}

Let $\dz\ge 0$ and $s\ge 2$. The energy functionals on the truncated hyperboloids $\mathscr{H}^{\rm{in}}_s$ and $\mathscr{H}^{\rm{ex}}_s$ are defined as
\beq\label{dinhs}
\mathcal{E}^{{\rm{in}},h}_{m,\dz}(s,u):=\int_{\mathscr{H}^{\rm{in}}_s}t^{-\dz}e^h_m[u]{\rm{d}}x,
\eeq
\beq\label{dexhs}
\mathcal{E}^{{\rm{ex}},h}_{m,\dz}(s,u):=\int_{\mathscr{H}^{\rm{ex}}_s}t^{-\dz}e^h_m[u]{\rm{d}}x
\eeq
respectively, where   
\beq\label{e^h_m}
\begin{split}
e^h_m[u]:&=|\partial u|^2+m^2u^2+2\f{x^a}{t}\partial_tu\partial_au=\f{s^2}{t^2}|\partial_tu|^2+t^{-2}|Lu|^2+m^2u^2\\
&=\f{s^2}{t^2}|\nabla u|^2+t^{-2}|L_0u|^2+t^{-2}|\Omega u|^2+m^2u^2
\end{split}
\eeq
 and we recall the definitions \eqref{s2: defNa} and \eqref{s2: defLO}. We denote for simplicity 
\beq\label{dhs0}
\mathcal{E}^{{\rm{in}},h}_{m}(s,u):=\mathcal{E}^{{\rm{in}},h}_{m,0}(s,u),\quad\quad \mathcal{E}^{{\rm{ex}},h}_m(s,u):=\mathcal{E}^{{\rm{ex}},h}_{m,0}(s,u).
\eeq

The energy estimate on the truncated exterior hyperboloids $\mathscr{H}^{\rm{ex}}_s$ is stated as follows.

\begin{prop}\label{s2: exhs}
Let $\dz\ge 0$. For any $s\in[2,\infty)$ and any sufficiently smooth function $u$ which is defined in $\mathscr{H}^{\rm{ex}}_{[2,s]}$ and decays sufficiently fast at space infinity, we have
\ben
&&\mathcal{E}^{{\rm{ex}},h}_{m,\dz}(s,u)+\int_{\mathscr{l}_{[2,\,t(s)]}}t^{-\dz}\big(|Gu|^2+m^2u^2\big){\rm{d}}\sz\\
&\lesssim&E^{\rm{ex}}_{m,\dz}(2,u)+\sup_{t\in[2,\infty)}[E^{\rm{ex}}_{m,\dz}(t,u)]^{\f{1}{2}}\cdot\int_2^\infty\|t^{-\f{\dz}{2}}(-\Box u+m^2u)\|_{L^2_x(\Sigma^s_t)}{\rm{d}}t,
\een
where $t(s)=\f{s^2+1}{2}$ is as in \eqref{s2: tr(s)}, we recall \eqref{s2: l[2,t]} and \eqref{dexfsdz} for the definitions of $\mathscr{l}_{[2,\,t(s)]}$ and $E^{\rm{ex}}_{m,\dz}(t,u)$, and
\beq\label{s2: Sigma^s_t}
\Sigma^s_t:=\l\{\begin{split}
\big\{x\in\mathbb{R}^n: |x|\ge t-1\big\},\quad& t\le t(s);\\
\big\{x\in\mathbb{R}^n: |x|\ge\sqrt{t^2-s^2}\big\}, \quad& t> t(s).
\end{split}\r.
\eeq

\begin{proof}
Recall the definition \eqref{s2: mscH}. We note that on $\mathscr{H}_{s}$ we have ${\bf n}{\rm{d}}\sz=(1,-x/t){\rm{d}}x$, where ${\bf n}$ and ${\rm{d}}\sz$ denote the upward unit normal and the volume element of $\mathscr{H}_s$ respectively. Taking $\omega\equiv 1$ in \eqref{s2: omd_tubo} and integrating over the region $\mathscr{H}^{\rm{ex}}_{[2,s]}$, we obtain
\ben
\mathcal{E}^{{\rm{ex}},h}_{m,\dz}(s,u)&+&\int_{\mathscr{l}_{[2,\,t(s)]}}t^{-\dz}\big(|Gu|^2+m^2u^2\big)2^{-\f{1}{2}}{\rm{d}}\sz+\dz\int_{\mathscr{H}^{\rm{ex}}_{[2,s]}}t^{-\dz-1}(|\partial u|^2+m^2u^2){\rm{d}}x{\rm{d}}t\\
&=&E^{\rm{ex}}_{m,\dz}(2,u)+2\int_2^\infty\int_{\Sigma^s_t}t^{-\dz}\partial_tu(-\Box u+m^2u){\rm{d}}x{\rm{d}}t.
\een
We note that $\Sigma^s_t\subset\Sigma^{\rm{ex}}_t$, hence the conclusion follows.
\end{proof}
\end{prop}

On the truncated interior hyperboloids $\mathscr{H}^{\rm{in}}_s$, we have the following energy estimate.

\begin{prop}\label{s2: inhs}
Let $\dz\ge 0$. For any $s\ge 2$ and any sufficiently smooth function $u$ defined in $\mathscr{H}^{\rm{in}}_{[2,s]}$, we have
\ben
\mathcal{E}^{{\rm{in}},h}_{m,\dz}(s,u)&\lesssim&\mathcal{E}^{{\rm{in}},h}_{m,\dz}(2,u)+\int_{\mathscr{l}_{[t(2),\,t(s)]}}t^{-\dz}\big(|Gu|^2+m^2u^2\big){\rm{d}}\sz\\
&+&\int_2^s[\mathcal{E}^{{\rm{in}},h}_{m,\dz}(\tau,u)]^{\f{1}{2}}\cdot\|t^{-\f{\dz}{2}}(-\Box u+m^2u)\|_{L^2(\mathscr{H}^{\rm{in}}_\tau)}{\rm{d}}\tau,
\een
where $t(s)=\f{s^2+1}{2}$ is as in \eqref{s2: tr(s)} and $t(2)=\f{5}{2}$\,, and we recall \eqref{s2: l[2,t]} for the definition of $\mathscr{l}_{[t(2),\,t(s)]}$.

\begin{proof}
Taking $\omega\equiv 1$ in \eqref{s2: omd_tubo}, integrating over the region $\mathscr{H}^{\rm{in}}_{[2,s]}$, and using the transformation $(t,x)\to (\tau,x)$, where $\tau=\sqrt{t^2-|x|^2}$, we obtain
\ben
\mathcal{E}^{{\rm{in}},h}_{m,\dz}(s,u)+\dz\int_{\mathscr{H}^{\rm{in}}_{[2,s]}}t^{-\dz-1}(|\partial u|^2+m^2u^2){\rm{d}}x{\rm{d}}t&=&\mathcal{E}^{{\rm{in}},h}_{m,\dz}(2,u)+\int_{\mathscr{l}_{[t(2),\,t(s)]}}t^{-\dz}\big(|Gu|^2+m^2u^2\big)2^{-\f{1}{2}}{\rm{d}}\sz\\
&+&2\int_2^s\int_{\mathscr{H}^{{\rm{in}}}_\tau}t^{-\dz}\partial_tu(-\Box u+m^2u)\f{\tau}{t}{\rm{d}}x{\rm{d}}\tau.
\een
Using H\"older inequality, the conclusion follows.
\end{proof}
\end{prop}

\subsection{Energy estimates for the Dirac equation}

Let $\psi$ be any $\mathbb{C}^{N_0}$-valued function defined on $\Sigma^{\rm{ex}}_t$, where $N_0$ is as in \eqref{s2: defN_0}. For any fixed constant $\lz>0$, we define
\beq\label{dexfD}
E^{{\rm{ex}},\lz}_D(t,\psi):=\int_{\Sigma^{\rm{ex}}_t}(2+r-t)^{1+\lz}|\psi|^2(t,x){\rm{d}}x.
\eeq
For the unweighted case we denote
\beq\label{dexfsD0} 
E^{{\rm{ex}}}_D(t,\psi):=\int_{\Sigma^{\rm{ex}}_t}|\psi|^2(t,x){\rm{d}}x.
\eeq
We also denote 
\beq\label{s2: psi-}
[\psi]_{-}:=\psi-\f{x_a}{r}\gz^0\gz^a\psi,\quad\quad [\psi]_{+}:=\psi+\f{x_a}{r}\gz^0\gz^a\psi.
\eeq

We next state the energy estimates for the Dirac equation on flat time slices in the exterior region.

\begin{prop}\label{s2: exfsD}
Let $\lz>0$. For any $M\in\mathbb{R}$, any $t\in[2,\infty)$ and any sufficiently smooth $\mathbb{C}^{N_0}$-valued function $\psi$ which is defined in $\mathscr{D}^{\rm{ex}}_t$ and decays sufficiently fast at space infinity, we have
\ben
\|\psi\|_{Y^{{\rm{ex}},\lz}_{D,t}}^2+(R_D(t,\psi))^2&\lesssim&\|\psi\|_{Y^{{\rm{ex}},\lz}_{D,2}}^2+\int_2^t\|(2+r-\tau)^{1+\lz}\psi^*\gz^0(i\gz^\mu\partial_\mu\psi+M\psi)\|_{L^1_x(\Sigma^{\rm{ex}}_\tau)}{\rm{d}}\tau,\\
\|\psi\|_{Y^{{\rm{ex}},\lz}_{D,t}}+R_D(t,\psi)&\lesssim&\|\psi\|_{Y^{{\rm{ex}},\lz}_{D,2}}+\int_2^t\|(2+r-\tau)^{\f{1+\lz}{2}}(i\gz^\mu\partial_\mu\psi+M\psi)\|_{L^2_x(\Sigma^{\rm{ex}}_\tau)}{\rm{d}}\tau,
\een
where
\beq\label{s2: exnormD}
\begin{split}
\|\psi\|_{Y^{{\rm{ex}},\lz}_{D,t}}:&=[E^{{\rm{ex}},\lz}_D(t,\psi)]^{\f{1}{2}}+\l(\int_2^t\int_{\Sigma^{\rm{ex}}_\tau}(2+r-\tau)^\lz|[\psi]_{-}|^2{\rm{d}}x{\rm{d}}\tau\r)^{\f{1}{2}},\\
R_D(t,\psi):&=\l(\int_{\mathscr{l}_{[2,t]}}|[\psi]_{-}|^2{\rm{d}}\sz\r)^{\f{1}{2}}.
\end{split}
\eeq

\begin{proof}
Let $\omega$ be a positive weight. Set $F_\psi:=i\gz^\mu\partial_\mu\psi+M\psi$. Multiplying on both sides of this identity by $-i\omega(r-t)\psi^*\gz^0$, we obtain
\be
\omega\psi^*\partial_t\psi+\omega\psi^*\gz^0\gz^a\partial_a\psi-iM\omega\psi^*\gz^0\psi=-i\omega\psi^*\gz^0F_\psi.
\ee
Taking the complex conjugate of the above equality, we have
\be
\omega\partial_t\psi^*\psi+\omega\partial_a\psi^*\gz^0\gz^a\psi+iM\omega\psi^*\gz^0\psi=i\omega F_\psi^*\gz^0\psi.
\ee
Adding the last two equalities, we derive
\beq\label{s2: opsi*D}
\partial_t\big(\omega\psi^*\psi\big)+\partial_a\big(\omega\psi^*\gz^0\gz^a\psi\big)+\f{1}{2}\omega'|[\psi]_{-}|^2=-i\omega(\psi^*\gz^0F_\psi-F_\psi^*\gz^0\psi),
\eeq
where we use the relation
\be
\psi^*\psi-\f{x_a}{r}\psi^*\gz^0\gz^a\psi=\f{1}{2}\l|\psi-\f{x_a}{r}\gz^0\gz^a\psi\r|^2.
\ee
Let $\omega(z):=(2+z)^{1+\lz}$ and integrate \eqref{s2: opsi*D} over the region $\mathscr{D}^{\rm{ex}}_t$, then we obtain
\ben
E^{{\rm{ex}},\lz}_D(t,\psi)&+&\f{1+\lz}{2}\int_2^t\int_{\Sigma^{\rm{ex}}_\tau}(2+r-\tau)^\lz|[\psi]_{-}|^2{\rm{d}}x{\rm{d}}\tau+\f{1}{2}\int_{\mathscr{l}_{[2,t]}}|[\psi]_{-}|^2\,2^{-\f{1}{2}}{\rm{d}}\sz\\
&=&E^{{\rm{ex}},\lz}_D(2,\psi)-i\int_2^t\int_{\Sigma^{\rm{ex}}_\tau}(2+r-\tau)^{1+\lz}(\psi^*\gz^0F_\psi-F_\psi^*\gz^0\psi){\rm{d}}x{\rm{d}}\tau.
\een
Differentiating this equality with respect to $t$ and using H\"older inequality, we obtain
\ben
\partial_t\Bigg\{E^{{\rm{ex}},\lz}_D(t,\psi)&+&\f{1+\lz}{2}\int_2^t\int_{\Sigma^{\rm{ex}}_\tau}(2+r-\tau)^\lz|[\psi]_{-}|^2{\rm{d}}x{\rm{d}}\tau+\f{1}{2}\int_{\mathscr{l}_{[2,t]}}|[\psi]_{-}|^2\,2^{-\f{1}{2}}{\rm{d}}\sz\Bigg\}\\
&\lesssim&[E^{{\rm{ex}},\lz}_D(t,\psi)]^{\f{1}{2}}\cdot\|(2+r-t)^{\f{1+\lz}{2}}F_\psi\|_{L^2_x(\Sigma^{\rm{ex}}_t)}.
\een
The conclusions of the proposition follow.
\end{proof}
\end{prop}

\begin{rem}\label{RemD}
Taking $\omega\equiv 1$ in \eqref{s2: opsi*D} and integrating over $[2,t]\times\mathbb{R}^n$, we obtain the standard energy inequality for the Dirac equation, i.e.,
\be
\|\psi(t)\|_{L^2_x(\mathbb{R}^n)}\lesssim\|\psi(2)\|_{L^2_x(\mathbb{R}^n)}+\int_2^t\|(i\gz^\mu\partial_\mu\psi+M\psi)(\tau)\|_{L^2_x(\mathbb{R}^n)}{\rm{d}}\tau,
\ee
for any sufficiently smooth function $\psi=\psi(t,x): [2,\infty)\times\mathbb{R}^{n}\to\mathbb{C}^{N_0}$ which decays sufficiently fast at space infinity.
\end{rem}

For any $\mathbb{C}^{N_0}$-valued function $\psi=\psi(t,x)$, we denote 
\beq\label{s2: psi)-}
(\psi)_{-}:=\psi-\f{x_a}{t}\gz^0\gz^a\psi,\quad\quad (\psi)_{+}:=\psi+\f{x_a}{t}\gz^0\gz^a\psi.
\eeq
The energy functionals for the Dirac field on the truncated hyperboloids $\mathscr{H}^{\rm{in}}_s$ and $\mathscr{H}^{\rm{ex}}_s$ are defined as
\beq\label{dinhsD}
\mathcal{E}^{{\rm{in}},h}_D(s,\psi):=\int_{\mathscr{H}^{\rm{in}}_s}\big\{|(\psi)_{-}|^2+|(s/t)\psi|^2\big\}{\rm{d}}x,
\eeq
\beq\label{dexhsD}
\mathcal{E}^{{\rm{ex}},h}_D(s,\psi):=\int_{\mathscr{H}^{\rm{ex}}_s}\big\{|(\psi)_{-}|^2+|(s/t)\psi|^2\big\}{\rm{d}}x.
\eeq

We have the following energy estimate for the Dirac equation on the exterior hyperboloids $\mathscr{H}^{\rm{ex}}_{s}$.

\begin{prop}\label{s2: exhsD}
For any $M\in\mathbb{R}$, any $s\in[2,\infty)$ and any sufficiently smooth $\mathbb{C}^{N_0}$-valued function $\psi$ which is defined in $\mathscr{H}^{\rm{ex}}_{[2,s]}$ and decays sufficiently fast at space infinity, we have
\ben
\mathcal{E}^{{\rm{ex}},h}_D(s,\psi)+\int_{\mathscr{l}_{[2,\,t(s)]}}|[\psi]_{-}|^2{\rm{d}}\sz
&\lesssim& E^{\rm{ex}}_D(2,\psi)+\int_2^\infty\|\psi^*\gz^0(i\gz^\mu\partial_\mu\psi+M\psi)\|_{L^1_x(\Sigma^s_t)}{\rm{d}}t\\
&\lesssim&E^{\rm{ex}}_D(2,\psi)+\sup_{t\in[2,\infty)}[E^{\rm{ex}}_D(t,\psi)]^{\f{1}{2}}\cdot\int_2^\infty\|i\gz^\mu\partial_\mu\psi+M\psi\|_{L^2_x(\Sigma^s_t)}{\rm{d}}t,
\een
where $t(s)=\f{s^2+1}{2}$ is as in \eqref{s2: tr(s)}, $\Sigma^s_t$ is as in \eqref{s2: Sigma^s_t}, and we recall \eqref{s2: l[2,t]}, \eqref{dexfsD0} and \eqref{s2: psi-}.

\begin{proof}
Taking $\omega\equiv 1$ in \eqref{s2: opsi*D} and integrating over the region $\mathscr{H}^{\rm{ex}}_{[2,s]}$, we obtain
\ben
\f{1}{2}\mathcal{E}^{{\rm{ex}},h}_D(s,\psi)+\f{1}{2}\int_{\mathscr{l}_{[2,\,t(s)]}}|[\psi]_{-}|^2\ 2^{-\f{1}{2}}{\rm{d}}\sz\le E^{\rm{ex}}_D(2,\psi)+2\int_2^\infty\int_{\Sigma^s_t}|\psi^*\gz^0F_\psi|{\rm{d}}x{\rm{d}}t,
\een
where $F_\psi:=i\gz^\mu\partial_\mu\psi+M\psi$ and we use the following identity on $\mathscr{H}_s$:
\be
\psi^*\psi-\f{x_a}{t}\psi^*\gz^0\gz^a\psi=\f{1}{2}\l(\l|\psi-\f{x_a}{t}\gz^0\gz^a\psi\r|^2+\f{s^2}{t^2}|\psi|^2\r).
\ee
The proof is done.
\end{proof}
\end{prop}

The energy estimate for the Dirac equation on the interior hyperboloids is stated as follows.

\begin{prop}\label{s2: inhsD}
For any $M\in\mathbb{R}$, any $s\ge 2$ and any sufficiently smooth $\mathbb{C}^{N_0}$-valued function $\psi$ defined in $\mathscr{H}^{\rm{in}}_{[2,s]}$, we have
\ben
\mathcal{E}^{{\rm{in}},h}_D(s,\psi)\lesssim \mathcal{E}^{{\rm{in}},h}_D(2,\psi)+\int_{\mathscr{l}_{[t(2),\,t(s)]}}|[\psi]_{-}|^2{\rm{d}}\sz
+\int_2^s[\mathcal{E}^{{\rm{in}},h}_D(\tau,\psi)]^{\f{1}{2}}\cdot\|i\gz^\mu\partial_\mu\psi+M\psi\|_{L^2(\mathscr{H}^{\rm{in}}_\tau)}{\rm{d}}\tau,
\een
where $t(s)=\f{s^2+1}{2}$ is as in \eqref{s2: tr(s)} and $t(2)=\f{5}{2}$\,, and we recall \eqref{s2: l[2,t]} for the definition of $\mathscr{l}_{[t(2),\,t(s)]}$.

\begin{proof}
Taking $\omega\equiv 1$ in \eqref{s2: opsi*D}, integrating over the region $\mathscr{H}^{\rm{in}}_{[2,s]}$ and using the transformation $(t,x)\to (\tau,x)$, where $\tau=\sqrt{t^2-|x|^2}$, we obtain
\ben
\f{1}{2}\mathcal{E}^{{\rm{in}},h}_D(s,\psi)=\f{1}{2}\mathcal{E}^{{\rm{in}},h}_D(2,\psi)+\f{1}{2}\int_{\mathscr{l}_{[t(2),\,t(s)]}}|[\psi]_{-}|^2\ 2^{-\f{1}{2}}{\rm{d}}\sz -i\int_2^s\int_{\mathscr{H}^{\rm{in}}_\tau}(\psi^*\gz^0F_\psi-F_\psi^*\gz^0\psi)\f{\tau}{t}{\rm{d}}x{\rm{d}}\tau,
\een
where $F_\psi:=i\gz^\mu\partial_\mu\psi+M\psi$. Using H\"older inequality, the conclusion of the proposition follows.
\end{proof}
\end{prop}

\subsection{Sobolev and Hardy inequalities}

We give the following weighted Sobolev inequality in the exterior region.

\begin{lem}\label{s2: Sobex}
Let $\Lz\in\mathbb{R}$. For any $t\in[2,\infty)$ and any sufficiently smooth function $w$ defined on $\Sigma^{\rm{ex}}_t$, we have
\beqn\label{s2: Soeq0}
&&\sup_{\Sigma^{\rm{ex}}_t}(2+r-t)^\Lz r^{n-1}|w(t,x)|^2\lesssim\sum_{|J|\le n-1}\int_{\Sigma^{\rm{ex}}_t}\big\{(2+r-t)^{\Lz+1}|\partial_r\Omega^Jw|^2+(2+r-t)^{\Lz-1}|\Omega^Jw|^2\big\}{\rm{d}}x,\nonumber\\
&&\sup_{\Sigma^{\rm{ex}}_t}(2+r-t)^\Lz r^{n-1}|w(t,x)|^2\lesssim\sum_{|J|\le n-1}\int_{\Sigma^{\rm{ex}}_t}(2+r-t)^{\Lz}\big\{|\partial_r\Omega^Jw|^2+|\Omega^Jw|^2\big\}{\rm{d}}x,
\eeqn
where we recall \eqref{s2: pOL^J}.
\begin{proof}
We only give the proof for $n=2$. The proof for $n=3$ is similar (using spherical coordinates instead) and was given in \cite[Lemmas 4.1 and 4.2]{HS}. Let $(r,\tz)$ be the polar coordinates in $\mathbb{R}^2$. Then 
\beq\label{s2: Sexeq0}
\sup_{\mathbb{S}^1}|w(t,r,\tz)|^2\lesssim\sum_{0\le l\le 1}\int_{\mathbb{S}^1}|\partial^l_\tz w(t,r,\tz)|^2{\rm{d}}\tz.
\eeq
For any fixed $t$ we denote $\zeta(r):=2+r-t$. Let $v$ be any smooth function which is compactly supported in $x$, then in the region $\{r\ge t-1\}$ we have
\ben
\partial_r\big((\zeta(r))^\Lz rv^2\big)\ge\Lz(\zeta(r))^{\Lz-1}r v^2+2(\zeta(r))^{\Lz}rv\partial_rv.
\een
Hence for $x\in\Sigma^{\rm{ex}}_t$ we have
\beq\label{s2: HS1}
(\zeta(r))^{\Lz} rv^2(t,x)\lesssim\int_r^\infty\rz\big\{(\zeta(\rz))^{\Lz-1}v^2+(\zeta(\rz))^{\Lz+1}|\partial_\rz v|^2\big\}{\rm{d}}\rz.
\eeq
If $w$ is compactly supported in $x$, we take $v=\partial^l_\tz w(t,r,\tz) (l=0,1)$ in \eqref{s2: HS1}, integrate over $\mathbb{S}^1$ and use \eqref{s2: Sexeq0}, and then obtain
\beqn\label{s2: Sexeq1}
&&\sup_{\Sigma^{\rm{ex}}_t}(\zeta(r))^\Lz r |w(t,x)|^2\lesssim\sup_{r\ge t-1}\sum_{0\le l\le 1}\int_{\mathbb{S}^1}(\zeta(r))^\Lz r |\partial^l_\tz w(t,r,\tz)|^2{\rm{d}}\tz\nonumber\\
&\lesssim&\sum_{0\le l\le 1}\int_{\mathbb{S}^1}\int_{t-1}^\infty r \big\{(\zeta(r))^{\Lz-1}|\partial^l_\tz w(t,r,\tz)|^2+(\zeta(r))^{\Lz+1}|\partial_r \partial^l_\tz w(t,r,\tz)|^2\big\}{\rm{d}}r{\rm{d}}\tz\nonumber\\
&\lesssim&\sum_{0\le l\le 1}\int_{\Sigma^{\rm{ex}}_t}\big\{(\zeta(r))^{\Lz-1}|\partial^l_\tz w(t,x)|^2+(\zeta(r))^{\Lz+1}|\partial_r\partial^l_\tz w(t,x)|^2\big\}{\rm{d}}x.
\eeqn
In the general case where $w$ is not compactly supported in $x$, we choose a cut-off function $\chi\in C^\infty_0(\mathbb{R})$ and apply \eqref{s2: Sexeq1} to $\chi(\ez r)w$ for any $\ez>0$, then we derive
\be
\begin{split}
&\quad\sup_{\Sigma^{\rm{ex}}_t}(\zeta(r))^\Lz r |\chi(\ez r)w(t,x)|^2\\
&\lesssim \sum_{0\le l\le 1}\int_{\Sigma^{\rm{ex}}_t}\l\{(\zeta(r))^{\Lz-1}|\chi(\ez r)\partial^l_\tz w|^2+(\zeta(r))^{\Lz+1}\l(|\chi(\ez r)\partial_r\partial^l_\tz w|^2+|\ez\chi'(\ez r)\partial^l_\tz w|^2\r)\r\}{\rm{d}}x.
\end{split}
\ee
Note that on $\Sigma^{\rm{ex}}_t$ we have $|\ez\chi'(\ez r)|^2\le(\zeta(r))^{2}|\ez\chi'(\ez r)|^2\le|\ez r\chi'(\ez r)|^2\lesssim 1$. Let $\ez\to 0$, and we obtain
\be
\sup_{\Sigma^{\rm{ex}}_t}(\zeta(r))^\Lz r |w(t,x)|^2 \lesssim\sum_{0\le l\le 1}\int_{\Sigma^{\rm{ex}}_t}\big\{(\zeta(r))^{\Lz-1}|\partial^l_\tz w|^2+(\zeta(r))^{\Lz+1}|\partial_r\partial^l_\tz w|^2\big\}{\rm{d}}x.
\ee
Taking the same weight $(\zeta(\rz))^{\Lz}$ on the right hand side of \eqref{s2: HS1} and arguing as above, we also obtain
\be
\sup_{\Sigma^{\rm{ex}}_t}(\zeta(r))^{\Lz} r |w(t,x)|^2 \lesssim\sum_{0\le l\le 1}\int_{\Sigma^{\rm{ex}}_t}(\zeta(r))^{\Lz}\l\{|\partial^l_\tz w|^2+|\partial_r\partial^l_\tz w|^2\r\}{\rm{d}}x.
\ee
We note that $\partial_\tz=x_1\partial_2-x_2\partial_1=\Omega_{12}$, hence the conclusions of the lemma follow.
\end{proof}
\end{lem}

\begin{lem}\label{s2: Hardyex}
Let $\Lz>-1$. For any $t\in [2,\infty)$ and any sufficiently smooth function $w$ which is defined on $\Sigma^{\rm{ex}}_t$ and is compactly supported in $x$, we have
\be
\int_{\Sigma^{\rm{ex}}_t}(2+r-t)^\Lz w^2{\rm{d}}x\lesssim\int_{\Sigma^{\rm{ex}}_t}(2+r-t)^{\Lz+2}|\partial_rw|^2{\rm{d}}x.
\ee
\begin{proof}
The proof is similar to \cite[Lemma 4.5]{HS}. Using that
\be
\partial_r\big((2+r-t)^{\Lz+1}r^{n-1}\big)\ge(\Lz+1)(2+r-t)^\Lz r^{n-1},
\ee
we obtain
\ben
\int_{\Sigma^{\rm{ex}}_t}(2+r-t)^\Lz w^2(t,x){\rm{d}}x&\le&\f{1}{\Lz+1}\int_{\Sigma^{\rm{ex}}_t}\f{x^a}{r^n}\partial_a\big((2+r-t)^{\Lz+1}r^{n-1}\big)w^2(t,x){\rm{d}}x\\
&=&-\f{1}{\Lz+1}\int_{|x|=t-1}w^2(t,x){\rm{d}}\sz(x)-\f{2}{\Lz+1}\int_{\Sigma^{\rm{ex}}_t}(2+r-t)^{\Lz+1}w\partial_rw{\rm{d}}x\\
&\le&\f{2}{\Lz+1}\l(\int_{\Sigma^{\rm{ex}}_t}(2+r-t)^{\Lz}w^2{\rm{d}}x\r)^{\f{1}{2}}\cdot\l(\int_{\Sigma^{\rm{ex}}_t}(2+r-t)^{\Lz+2}|\partial_rw|^2{\rm{d}}x\r)^{\f{1}{2}}.
\een
The conclusion of the lemma follows.
\end{proof}
\end{lem}

We denote $\{Z_k\}_{k=1}^{n_2}:=\{(L_a)_{1\le a\le n},(\Omega_{ab})_{1\le a<b\le n}\}$, where $n_2:=(n^2+n)/2$. For any multi-index $J=(j_1,\cdots,j_{n_2})\in\mathbb{N}^{n_2}$ of length $|J|=\sum_{k=1}^{n_2}j_k$, we denote 
\beq\label{s2: Z^J}
Z^J=\prod_{k=1}^{n_2}Z_k^{j_k}.
\eeq 

\begin{lem}\label{s2: exL_auL^2}
Let $\Lz>-1$. For any $t\in[2,\infty)$ and any sufficiently smooth function $w$ defined on $\Sigma^{\rm{ex}}_t$, we have
\be
\begin{split}
&\int_{\Sigma^{\rm{ex}}_t}(2+r-t)^\Lz|Z_kw|^2{\rm{d}}x\lesssim\sum_{|J|\le 1}\int_{\Sigma^{\rm{ex}}_t}(2+r-t)^{\Lz+2}|\partial Z^{J}w|^2{\rm{d}}x,\quad k\in\{1,\cdots,n_2\},\\
&\int_{\Sigma^{\rm{ex}}_t}(2+r-t)^\Lz|L_0w|^2{\rm{d}}x\lesssim\sum_{|J|\le 1}\int_{\Sigma^{\rm{ex}}_t}(2+r-t)^{\Lz+2}|\partial Z^Jw|^2{\rm{d}}x.
\end{split}
\ee
\begin{proof}
The second inequality follows from the first one, using the relation 
\beq\label{L_0}
L_0=(t-r)\partial_t+(r-t)\partial_r+(x^a/r)L_a
\eeq
and the fact that 
\beq\label{s2: Edecay0}
|t-r|\lesssim 2+r-t\quad\mathrm{in}\ \mathscr{D}^{\rm{ex}}.
\eeq 
The proof of the first inequality is as in \cite[Corollary 4.6]{HS}. We choose a cut-off function $\chi\in C^\infty_0(\mathbb{R})$, apply Lemma \ref{s2: Hardyex} to $\chi(\ez r)Z_kw$ for any $\ez>0$, and obtain
\ben
\int_{\Sigma^{\rm{ex}}_t}(2+r-t)^\Lz|\chi(\ez r)Z_kw|^2{\rm{d}}x&\lesssim&\int_{\Sigma^{\rm{ex}}_t}(2+r-t)^{\Lz+2}\l\{|\chi(\ez r)\partial_rZ_kw|^2+|\ez\chi'(\ez r)Z_kw|^2\r\}{\rm{d}}x\\
&\lesssim&\int_{\Sigma^{\rm{ex}}_t}(2+r-t)^{\Lz+2}\l\{|\partial_rZ_kw|^2+|\partial w|^2\r\}{\rm{d}}x,
\een
where we use that $|Z_kw|\lesssim r|\partial w|$ on $\Sigma^{\rm{ex}}_t$ and that $|\ez r\chi'(\ez r)|^2\lesssim 1$. Let $\ez\to 0$, and the proof is done.
\end{proof}
\end{lem}

\begin{lem}\label{s2: exL_auL^infty}
Let $\lz>0$. For any $t\in[2,\infty)$ and any sufficiently smooth function $w$ defined on $\Sigma^{\rm{ex}}_t$, we have
\be
\begin{split}
&\sup_{\Sigma^{\rm{ex}}_t}(2+r-t)^\lz r^{n-1}|Z_kw|^2\lesssim\sum_{|K|\le n}\int_{\Sigma^{\rm{ex}}_t}(2+r-t)^{\lz+1}|\partial Z^{K}w|^2{\rm{d}}x,\quad k\in\{1,\cdots,n_2\},\\
&\sup_{\Sigma^{\rm{ex}}_t}(2+r-t)^{\lz-1} r^{n-1}|L_0w|^2\lesssim\sum_{|K|\le n+1}\int_{\Sigma^{\rm{ex}}_t}(2+r-t)^{\lz+1}|\partial \Gamma^{K}w|^2{\rm{d}}x,
\end{split}
\ee
where we recall \eqref{s2: Ga^I} and \eqref{s2: Z^J}.
\begin{proof}
For any fixed $t$ we denote $\zeta(r):=2+r-t$. By Lemmas \ref{s2: Sobex} and \ref{s2: exL_auL^2}, we have
\be
\begin{split}
\sup_{\Sigma^{\rm{ex}}_t}(\zeta(r))^\lz r^{n-1}|Z_kw|^2&\lesssim\sum_{|J|\le n-1}\int_{\Sigma^{\rm{ex}}_t}\l\{(\zeta(r))^{\lz+1}|\partial_r\Omega^JZ_kw|^2+(\zeta(r))^{\lz-1}|\Omega^JZ_kw|^2\r\}{\rm{d}}x\\
&\lesssim\sum_{|K|\le n}\int_{\Sigma^{\rm{ex}}_t}(\zeta(r))^{\lz+1}|\partial Z^{K}w|^2{\rm{d}}x,\\
\sup_{\Sigma^{\rm{ex}}_t}(\zeta(r))^{\lz-1} r^{n-1}|L_0w|^2&\lesssim\sum_{|J|\le n-1}\int_{\Sigma^{\rm{ex}}_t}(\zeta(r))^{\lz-1}\big\{|\partial_r\Omega^JL_0w|^2+|\Omega^JL_0w|^2\big\}{\rm{d}}x\\
&\lesssim\sum_{|J|\le n-1}\int_{\Sigma^{\rm{ex}}_t}(\zeta(r))^{\lz-1}\big\{|L_0\partial\Omega^Jw|^2+|\partial\Omega^Jw|^2+|L_0\Omega^Jw|^2\big\}{\rm{d}}x\\
&\lesssim\sum_{|K|\le n+1}\int_{\Sigma^{\rm{ex}}_t}(\zeta(r))^{\lz+1}|\partial \Gamma^{K}w|^2{\rm{d}}x,
\end{split}
\ee
where we use the estimates on commutators in \eqref{s2: commutator}. The proof is done.
\end{proof}
\end{lem}

We give the following Sobolev inequality on hyperboloids; see \cite{LM18}.

\begin{lem}\label{s4: Sobin}
For any sufficiently smooth function $w$ defined in the cone $\mathscr{C}:=\{r<t\}$ and all $(t,x)\in\mathscr{C}$, we have
\be
|w(t,x)|^2\lesssim t^{-n}\sum_{|J|\le 2}\int_{B(x,t/3)}|L^Jw\big(\sqrt{s^2+|y|^2},y\big)|^2{\rm{d}}y,
\ee
where $s:=\sqrt{t^2-|x|^2}$, $B(x,t/3)\subset\mathbb{R}^n$ denotes the ball centered at $x$ with radius $t/3$, and we recall \eqref{s2: pOL^J}.
\end{lem}

\begin{lem}\label{sC: DL2.3'}
For any sufficiently smooth scalar function $w$ and $\mathbb{C}^{N_0}$-valued function $\phi$ (recall \eqref{s2: defN_0}), the following estimates hold in the region $\mathscr{C}:=\{r<t\}$:
\be
\begin{split}
&\sum_{|K|\le 2}|L^K\big((s/t)w\big)|\lesssim\sum_{|K|\le 2}|(s/t)L^Kw|,\\
&\sum_{|K|\le 2}|\hat{L}^K(\phi)_{-}|\lesssim\sum_{|K|\le 2}|(\hat{L}^K\phi)_{-}|,
\end{split}
\ee
where we recall \eqref{s2: psi)-}, \eqref{s2: hatLO} and \eqref{s2: pOL^J}.
\begin{proof}
The first inequality follows from the identities below (recall \eqref{s2: barpar_a}):
\be
L_a\l(\f{s}{t}\r)=t\bar{\partial}_a\l(\f{s}{t}\r)=-\f{s}{t}\f{x_a}{t},\quad\quad L_bL_a\l(\f{s}{t}\r)=-t\bar{\partial}_b\l(\f{s}{t}\f{x_a}{t}\r)=\l(2\f{x_ax_b}{t^2}-\dz_{ab}\r)\f{s}{t}.
\ee
The proof of the second inequality can be found in \cite[Proposition 4.2]{DW}. We sketch it below. Let ${\boldsymbol{I}}:=I_{N_0}$ denote the $N_0\times N_0$ identity matrix and $\xi_a=\xi_a(t,x):=(x_a/t)$ for $1\le a\le n$. By direct computation,
\be
\begin{split}
[\hat{L}_b,{\boldsymbol{I}}-\xi_a\gz^0\gz^a]&=-(\gz^0\gz^b+\xi_b)({\boldsymbol{I}}-\xi_a\gz^0\gz^a),\\
[\hat{L}_c,\gz^0\gz^b+\xi_b]&=-\xi_b\xi_c-\gz^b\gz^c.
\end{split}
\ee
It follows that
\be
\begin{split}
\hat{L}_c[\hat{L}_b,{\boldsymbol{I}}-\xi_a\gz^0\gz^a]&=-\hat{L}_c(\gz^0\gz^b+\xi_b)({\boldsymbol{I}}-\xi_a\gz^0\gz^a)\\
&=-(\gz^0\gz^b+\xi_b)\hat{L}_c({\boldsymbol{I}}-\xi_a\gz^0\gz^a)+(\xi_b\xi_c+\gz^b\gz^c)({\boldsymbol{I}}-\xi_a\gz^0\gz^a)\\
&=-(\gz^0\gz^b+\xi_b)\big\{({\boldsymbol{I}}-\xi_a\gz^0\gz^a)\hat{L}_c-(\gz^0\gz^c+\xi_c)({\boldsymbol{I}}-\xi_a\gz^0\gz^a)\big\}\\
&\quad+(\xi_b\xi_c+\gz^b\gz^c)({\boldsymbol{I}}-\xi_a\gz^0\gz^a)\\
&=-(\gz^0\gz^b+\xi_b)({\boldsymbol{I}}-\xi_a\gz^0\gz^a)\hat{L}_c+(\xi_b\gz^0\gz^c+\xi_c\gz^0\gz^b+2\xi_b\xi_c)({\boldsymbol{I}}-\xi_a\gz^0\gz^a),\\
[\hat{L}_c\hat{L}_b,{\boldsymbol{I}}-\xi_a\gz^0\gz^a]
&=\hat{L}_c[\hat{L}_b,{\boldsymbol{I}}-\xi_a\gz^0\gz^a]+[\hat{L}_c,{\boldsymbol{I}}-\xi_a\gz^0\gz^a]\hat{L}_b\\
&=-(\gz^0\gz^b+\xi_b)({\boldsymbol{I}}-\xi_a\gz^0\gz^a)\hat{L}_c-(\gz^0\gz^c+\xi_c)({\boldsymbol{I}}-\xi_a\gz^0\gz^a)\hat{L}_b\\
&\quad+(\xi_b\gz^0\gz^c+\xi_c\gz^0\gz^b+2\xi_b\xi_c)({\boldsymbol{I}}-\xi_a\gz^0\gz^a).
\end{split}
\ee 
The proof is done.
\end{proof}
\end{lem}

\subsection{Estimates on good derivatives and the $\gz^0$-structure}

\begin{lem}\label{sC: DMY2.1}
For any sufficiently smooth function $u$ defined in $[2,\infty)\times\mathbb{R}^{n}$, we have
\be
|(t-r)\partial u|+(t+r)|G u|\lesssim |L_0u|+\sum_{|I|=1}|Z^Iu|,
\ee
where we recall \eqref{s2: G_a} and \eqref{s2: Z^J}.
\begin{proof}
We first recall the estimate of $|\partial u|$ below (see for example \cite{Sog})
\be\label{paru}
|(t-r)\partial u|\lesssim|L_0u|+\sum_{|I|=1}|Z^Iu|.
\ee
The estimate of $|Gu|$ follows from this and the identities below:
\be\label{G_au}
G_au=\f{1}{r}\l(L_a u+(r-t)\partial_a u\r)=\f{1}{t}\l(L_a u-\f{x_a}{r}(r-t)\partial_t u\r).
\ee
The proof is done.
\end{proof}
\end{lem}

\begin{lem}\label{sC: 2.2-4}
The following statements hold:
\begin{itemize}
\item[$i)$] Let $u: [2,\infty)\times\mathbb{R}^{n}\to\mathbb{R}$, $\phi,\tilde{\phi}: [2,\infty)\times\mathbb{R}^{n}\to\mathbb{C}^{N_0}$ be sufficiently smooth functions. Then
\beq\label{s2: 2.2-41}
\hat{\Gamma}^I(u\phi)=\sum_{I_1+I_2=I}(\Gamma^{I_1}u)(\hat{\Gamma}^{I_2}\phi),
\eeq
\beq\label{s2: 2.2-42}
\Gamma^I(\phi^*\gz^0\tilde{\phi})=\sum_{I_1+I_2=I}(\hat{\Gamma}^{I_1}\phi)^*\gz^0(\hat{\Gamma}^{I_2}\tilde{\phi}),
\eeq
\beq\label{s2: 2.2-43'}
|\phi^{*}\gz^0\tilde{\phi}|\lesssim|[\phi]_{-}|\cdot|\tilde{\phi}|+|\phi|\cdot|[\tilde{\phi}]_{-}|,
\eeq
\beq\label{s2: 2.2-43}
|\phi^{*}\gz^0\tilde{\phi}|\lesssim|(\phi)_{-}|\cdot|\tilde{\phi}|+|\phi|\cdot|(\tilde{\phi})_{-}|+\f{|t^2-r^2|}{t^2}|\phi|\cdot|\tilde{\phi}|\quad\quad\mathrm{in}\ \{r<t\},
\eeq
for any multi-index $I\in\mathbb{N}^{n_0}$, where we recall \eqref{s2: Ga^I}, \eqref{s2: psi-} and \eqref{s2: psi)-}.
\item[$ii)$] For any sufficiently smooth function $\phi: [2,\infty)\times\mathbb{R}^{n}\to \mathbb{C}^{N_0}$, let $\tilde{\phi}:=i\gz^\mu\partial_\mu\phi$. Then 
\beq\label{s2: 2.2-44}
|[\tilde{\phi}]_{-}|\lesssim|G\phi|,
\eeq
where we recall \eqref{s2: G_a}.
\end{itemize}

\begin{proof}
$i)$ Recall the definition \eqref{s2: gamma_k}. We have
\be
\hat{L}_a(u\phi)=(L_au)\phi+u(\hat{L}_a\phi),\quad\quad \mathrm{for\ }1\le a\le n
\ee
and similarly for $\hat{\Omega}_{ab}, 1\le a<b\le n$. We also have 
\be
\begin{split}
L_a(\phi^*\gz^0\tilde{\phi})=(L_a\phi)^*\gz^0\tilde{\phi}+\phi^*\gz^0L_a\tilde{\phi}&=(\hat{L}_a\phi)^*\gz^0\tilde{\phi}+\f{1}{2}\phi^*\gz^0\gz^a\gz^0\tilde{\phi}+\phi^*\gz^0\hat{L}_a\tilde{\phi}+\f{1}{2}\phi^*\gz^0\gz^0\gz^a\tilde{\phi}\\
&=(\hat{L}_a\phi)^*\gz^0\tilde{\phi}+\phi^*\gz^0\hat{L}_a\tilde{\phi},\quad\quad \mathrm{for\ }1\le a\le n,\\
\Omega_{ab}(\phi^*\gz^0\tilde{\phi})=(\Omega_{ab}\phi)^*\gz^0\tilde{\phi}+\phi^*\gz^0\Omega_{ab}\tilde{\phi}&=(\hat{\Omega}_{ab}\phi)^*\gz^0\tilde{\phi}+\f{1}{2}\phi^*\gz^b\gz^a\gz^0\tilde{\phi}+\phi^*\gz^0\hat{\Omega}_{ab}\tilde{\phi}+\f{1}{2}\phi^*\gz^0\gz^a\gz^b\tilde{\phi}\\
&=(\hat{\Omega}_{ab}\phi)^*\gz^0\tilde{\phi}+\phi^*\gz^0\hat{\Omega}_{ab}\tilde{\phi},\quad\quad \mathrm{for\ }1\le a<b\le n.
\end{split}
\ee
By induction, we obtain \eqref{s2: 2.2-41} and \eqref{s2: 2.2-42}. Next, we write
\be
\phi=\f{(\phi)_{-}+(\phi)_{+}}{2}=\f{[\phi]_{-}+[\phi]_{+}}{2},\quad\quad \tilde{\phi}=\f{(\tilde{\phi})_{-}+(\tilde{\phi})_{+}}{2}=\f{[\tilde{\phi}]_{-}+[\tilde{\phi}]_{+}}{2}.
\ee
Denote by ${\boldsymbol{I}}:=I_{N_0}$ the $N_0\times N_0$ identity matrix, and let $\xi_a:=x_a/t$, $\eta_a:=x_a/r$. By direct computation,
\be\label{sC: A-}
({\boldsymbol{I}}-\xi_b\gz^0\gz^b)\gz^0({\boldsymbol{I}}-\xi_a\gz^0\gz^a)=\gz^0-\xi_b\gz^0\gz^b\gz^0-\xi_a\gz^0\gz^0\gz^a+\xi_a\xi_b\gz^0\gz^b\gz^0\gz^0\gz^a=(1-r^2/t^2)\gz^0
\ee
and similarly for $({\boldsymbol{I}}+\xi_b\gz^0\gz^b)\gz^0({\boldsymbol{I}}+\xi_a\gz^0\gz^a)$. We also have
\beq\label{sC: T-}
({\boldsymbol{I}}-\eta_b\gz^0\gz^b)\gz^0({\boldsymbol{I}}-\eta_a\gz^0\gz^a)=\gz^0-\eta_b\gz^0\gz^b\gz^0-\eta_a\gz^0\gz^0\gz^a+\eta_a\eta_b\gz^0\gz^b\gz^0\gz^0\gz^a=0
\eeq
and similarly for $({\boldsymbol{I}}+\eta_b\gz^0\gz^b)\gz^0({\boldsymbol{I}}+\eta_a\gz^0\gz^a)$. It follows that
\be
\phi^{*}\gz^0\tilde{\phi}=\f{1}{4}\big\{(\phi)^*_{-}\gz^0(\tilde{\phi})_{+}+(\phi)^*_{+}\gz^0(\tilde{\phi})_{-}+2(1-r^2/t^2)\phi^*\gz^0\tilde{\phi}\big\}=\f{1}{4}\big\{[\phi]^*_{-}\gz^0[\tilde{\phi}]_{+}+[\phi]^*_{+}\gz^0[\tilde{\phi}]_{-}\big\}.
\ee
Using the definitions \eqref{s2: psi-} and \eqref{s2: psi)-}, we obtain \eqref{s2: 2.2-43'} and \eqref{s2: 2.2-43}.

$ii)$ Using the relation $\partial_a=G_a-(x_a/r)\partial_t$, we can write
\be
\gz^0\partial_t+\gz^a\partial_a=\gz^0\l({\boldsymbol{I}}-\f{x_a}{r}\gz^0\gz^a\r)\partial_t+\gz^aG_a.
\ee
Using \eqref{sC: T-}, we obtain 
\be
[\tilde{\phi}]_{-}=i\l({\boldsymbol{I}}-\f{x_b}{r}\gz^0\gz^b\r)\gz^\mu\partial_\mu\phi=i\l({\boldsymbol{I}}-\f{x_b}{r}\gz^0\gz^b\r)\gz^aG_a\phi.
\ee
The proof is done.
\end{proof}
\end{lem}

\subsection{Nonlinear transforms}

The following lemma is obtained by straightforward computation; see also \cite{DW, DLMY}. 

\begin{lem}\label{s2: nlt}
Let $(\psi,v)$ be the solution to \eqref{s2: nDKG}-\eqref{s2: ninitial} with $n=2$ and $M=0$. Then the following statements hold:
\begin{itemize}
\item[$i)$] Let $\tilde{\psi}:=\psi-i\gz^\mu\partial_\mu(v\psi)$. Then $\tilde{\psi}$ solves the equation
\beq\label{s2: tildeG}
i\gz^\mu\partial_\mu\tilde{\psi}=\tilde{F}_\psi:=(\psi^*\gz^0\psi)\psi-iv\gz^\mu\partial_\mu(v\psi)-2\partial_\az v\partial^\az\psi.
\eeq
\item[$ii)$] Let $\tilde{v}:=v-\psi^*\gz^0\psi$. Then $\tilde{v}$ solves the equation
\beq\label{s2: tildeF}
-\Box\tilde{v}+\tilde{v}=\tilde{F}_v:=-i\partial_\mu(v\psi^*)\gz^0\gz^\mu\psi+i\psi^*\gz^0\gz^\mu\partial_\mu(v\psi)+2\partial_\az\psi^*\gz^0\partial^\az\psi.
\eeq
\item[$iii)$] Let $\Psi$ be the solution to
\beq\label{s2: boxPsi}
\Box\Psi=i\gz^\mu\partial_\mu\psi=v\psi,\quad\quad(\Psi,\partial_t\Psi)|_{t=2}=(0,-i\gz^0\psi_0).
\eeq
Then we have $\psi=i\gz^\mu\partial_\mu\Psi$. Let $\tilde{\Psi}:=\Psi-v\psi$. Then $\tilde{\Psi}$ solves the equation
\beq\label{s2: boxtiPsi}
\Box\tilde{\Psi}=\tilde{F}_\psi=(\psi^*\gz^0\psi)\psi-iv\gz^\mu\partial_\mu(v\psi)-2\partial_\az v\partial^\az\psi.
\eeq
\end{itemize}
\end{lem}

We next give estimates of the null forms in the above lemma.

\begin{lem}\label{s2: Q_0}
For any sufficiently smooth functions $u$ and $\phi$, we have
\ben
|\partial_\az u\partial^\az \phi|&\lesssim&t^{-1}\sum_{a=1}^n\big\{|L_a u|\cdot|\partial\phi|+|\partial u|\cdot|L_a\phi|+|t-r|\cdot|\partial u|\cdot|\partial\phi|\big\},\\
|\Gamma_k\l(\partial_\az u\partial^\az \phi\r)|&\lesssim&|(\partial_\az \Gamma_k u)(\partial^\az \phi)|+|(\partial_\az u)(\partial^\az \Gamma_k\phi)|,\quad k=1,\cdots,n_0,
\een
where we recall \eqref{s2: gamma_k}.
\begin{proof}
The second inequality is well-known; see for example \cite{Sog}. For the first one, using the relation $\partial_a=\bar{\partial}_a-(x_a/t)\partial_t$ where $\bar{\partial}_a=t^{-1}L_a$, we have
\ben
-\partial_\az u\partial^\az \phi&=&\f{1}{t}\partial_tu\l(t\partial_t\phi+x^a\partial_a\phi-x^a\partial_a\phi\r)-\l(\bar{\partial}_au-\f{x_a}{t}\partial_tu\r)\partial^a\phi\\
&=&t^{-1}\partial_tu\l(t\partial_t\phi+r\partial_r\phi\r)-\bar{\partial}_au\partial^a\phi\\
&=&t^{-1}\partial_tu\big((t-r)\partial_t\phi+(r-t)\partial_r\phi+(x^a/r)L_a\phi)\big)-\bar{\partial}_au\partial^a\phi\\
&=&\f{t-r}{t}\partial_tu\partial_t\phi+\f{r-t}{t}\partial_tu\partial_r\phi+\f{x^a}{r}\partial_tu\bar{\partial}_a\phi-\bar{\partial}_au\partial^a\phi.
\een
The proof is done.
\end{proof}
\end{lem}

\begin{lem}\label{s2: Edecay}
For any sufficiently smooth scalar function $v$ and $\mathbb{C}^{N_0}$-valued vector field $\psi$, we have
\beqn\label{s2: parpsi}
&&|v|\lesssim\f{2+r-t}{t}\sum_{|I|\le 1}|\partial\Gamma^Iv|+|-\Box v+v|\quad\quad\mathrm{in}\ \mathscr{D}^{\rm{ex}}\cap\{r\le 3t\},\nonumber\\
&&|\partial\psi|\lesssim\f{1}{t-r}\sum_{|I|=1}|\Gamma^I\psi|+\f{t}{t-r}|i\gz^\mu\partial_\mu\psi|\quad\quad\mathrm{in}\ \mathscr{D}^{\rm{in}}.
\eeqn

\begin{proof}
We write the d'Alembert operator $-\Box$ as
\be
-\Box=\f{(t-r)(t+r)}{t^2}\partial_t\partial_t+\f{x^a}{t^2}\partial_tL_a-\f{1}{t}\partial^aL_a+\f{n}{t}\partial_t-\f{x^a}{t^2}\partial_a.
\ee
Then the estimate of $v$ follows from \eqref{s2: Edecay0}. The estimate of $|\partial\psi|$ was given in \cite[Lemma 4.3]{DW}. Using the relation $\partial_a=t^{-1}L_a-(x_a/t)\partial_t$, we rewrite the identity $F_\psi:=i\gz^\mu\partial_\mu\psi$ as
\be
i\l(\gz^0-\f{x_a}{t}\gz^a\r)\partial_t\psi=-i\f{1}{t}\gz^aL_a\psi+F_\psi.
\ee
Multiplying on both sides of this equation by $-i\l(\gz^0-\f{x_b}{t}\gz^b\r)$ yields
\be
\f{t^2-r^2}{t^2}|\partial_t\psi|\lesssim \f{1}{t}\sum_{a=1}^n|L_a\psi|+|F_\psi|\quad\quad\mathrm{in}\ \mathscr{D}^{\rm{in}},
\ee
which implies
\be
|\partial_t\psi|\lesssim\f{1}{t-r}\sum_{a=1}^n|L_a\psi|+\f{t}{t-r}|F_\psi|\quad\quad\mathrm{in}\ \mathscr{D}^{\rm{in}}.
\ee
Using the relation $\partial_a=t^{-1}L_a-(x_a/t)\partial_t$ again, we obtain the same estimate for $\partial_a\psi$.
\end{proof}
\end{lem}

\subsection{A scattering lemma}

In this subsection, we prove a technical lemma which will be used in Section \ref{s6} to prove the scattering results in Theorems \ref{thm1} and \ref{thm2}.

Let $(\psi,v)$ be the global solution to \eqref{s2: nDKG}-\eqref{s2: ninitial}. We denote $\vec{v}=(v,\partial_tv)'=\l(\begin{array}{c}
v\\
\partial_tv
\end{array}\r)$, where $(a_1,a_2)'$ denotes the transpose of a vector $\vec{a}=(a_1,a_2)\in\mathbb{R}^2$. Let
\beq\label{s5: vecF}
\vec{F}_v=(0,F_v)',\quad\quad\vec{v}_0=(v_0,v_1)'.
\eeq 
By the linear theory of Dirac and Klein-Gordon equations, we have
\beq\label{s5: psi=S}
\psi=\mathcal{S}(t-2)\psi_0-i\int_{2}^t\mathcal{S}(t-\tau)\gz^0F_\psi(\tau)\rm{d}\tau,
\eeq
\beq\label{s5: vecv=}
\vec{v}=\mathcal{\tilde{S}}(t-2)\vec{v}_0+\int_2^t\mathcal{\tilde{S}}(t-\tau)\vec{F}_v(\tau){\rm{d}}\tau,
\eeq 
where $\mathcal{S}(t):=e^{it(i\gz^0\gz^a\partial_a+M\gz^0)}$ is the propagator for the linear Dirac equation, and
\be
\mathcal{\tilde{S}}(t)=\l(\begin{array}{cc}
\cos(t\langle\nabla\rangle)&\f{\sin(t\langle\nabla\rangle)}{\langle\nabla\rangle}\\
-\langle\nabla\rangle\sin(t\langle\nabla\rangle)&\cos(t\langle\nabla\rangle)
\end{array}\r).
\ee
Let $l\in\mathbb{N}$. We denote ${\bf H}^l(\mathbb{R}^n):=H^{l+1}(\mathbb{R}^n)\times H^l(\mathbb{R}^n)$.

\begin{lem}\label{Sca}
The following statements hold:
\begin{itemize}
\item[i)] Let $l\in\mathbb{N}$ and $\dz>0$. For any scalar function $g$ and $\mathbb{R}^2$-valued function $\vec{f}(\tau,x)=(f_1,f_2)'$ which are defined in $[2,\infty)\times\mathbb{R}^n$ and satisfy $g(\tau,\cdot)\in H^l(\mathbb{R}^n)$, $\vec{f}(\tau,\cdot)\in{\bf{H}}^l(\mathbb{R}^n)$ for any fixed $\tau\in[2,\infty)$, and any $t\ge 2$, $4\le T_1<T_2<\infty$, we have
\be
\begin{split}
&\quad\quad\l\|\int_{T_1}^{T_2}\mathcal{S}(t-\tau)(g(\tau)){\rm{d}}\tau\r\|_{H^{l}(\mathbb{R}^n)}\\
&\lesssim T_1^{-\f{\dz}{2}}\sum_{k=0}^{l}\Bigg\{\l(\int_{T_1}^{T_2}\|\nabla^kg\|_{L^2_x(\Sigma^{\rm{ex}}_\tau)}^2 \cdot\tau^{1+\dz}{\rm{d}}\tau\r)^{\f{1}{2}}+\l(\int_{T_1^{\f{1}{2}}}^{T_2}\|\nabla^kg\|_{L^2(\mathscr{H}^{\rm{in}}_s)}^2\cdot s^{1+2\dz}{\rm{d}}s\r)^{\f{1}{2}}\Bigg\},\\
&\quad\quad\l\|\int_{T_1}^{T_2}\mathcal{\tilde{S}}(t-\tau)(\vec{f}(\tau)){\rm{d}}\tau\r\|_{{\bf{H}}^l(\mathbb{R}^n)}\\
&\lesssim T_1^{-\f{\dz}{2}}\sum_{k=0}^l\Bigg\{\l(\int_{T_1}^{T_2}\|\,|\vec{\nabla}^k\vec{f}|+|f_1|\,\|_{L^2_x(\Sigma^{\rm{ex}}_\tau)}^2 \cdot\tau^{1+\dz}{\rm{d}}\tau\r)^{\f{1}{2}}+\l(\int_{T_1^{\f{1}{2}}}^{T_2}\|\,|\vec{\nabla}^k\vec{f}|+|f_1|\,\|_{L^2(\mathscr{H}^{\rm{in}}_s)}^2\cdot s^{1+2\dz}{\rm{d}}s\r)^{\f{1}{2}}\Bigg\},
\end{split}
\ee
where $|\vec{\nabla}^k\vec{f}|:=|\nabla^{k+1}f_1|+|\nabla^kf_2|$ and $\nabla=(\partial_1,\cdots,\partial_n)$.

\item[ii)] Let $l\in\mathbb{N}$ and $\psi,\vec{v},F_\psi,\vec{F}_v,\psi_0,\vec{v}_0$ be as in \eqref{s5: vecF}, \eqref{s5: psi=S} and \eqref{s5: vecv=} with $\psi_0\in H^{l}(\mathbb{R}^n), \vec{v}_0\in{\bf{H}}^{l}(\mathbb{R}^n)$ and $F_\psi(\tau),F_v(\tau)\in H^l(\mathbb{R}^n)$ for any fixed $\tau\in[2,\infty)$. If for some $\dz>0$, it holds that
\be
A:=\int_2^4\|F(\tau,\cdot)\|_{L^2_x(\mathbb{R}^n)}{\rm{d}}\tau+\l(\int_{4}^{+\infty}\|F(\tau,\cdot)\|_{L^2_x(\Sigma^{\rm{ex}}_\tau)}^2 \cdot\tau^{1+\dz}{\rm{d}}\tau+\int_{2}^{+\infty}\|F\|_{L^2(\mathscr{H}^{\rm{in}}_s)}^2\cdot s^{1+2\dz}{\rm{d}}s\r)^{\f{1}{2}}<\infty,
\ee
where $F:=\sum_{k=0}^{l}\big(|\nabla^{k}F_\psi|+|\nabla^kF_v|\big)$, then the solution $(\psi,v)$ scatters to a free solution in $H^{l}(\mathbb{R}^n)\times {\bf{H}}^{l}(\mathbb{R}^n)$, i.e., there exist $\psi^*_{0}\in H^{l}(\mathbb{R}^n)$ and $\vec{v}^*_0=(v^*_0,v^*_1)'\in{\bf{H}}^{l}(\mathbb{R}^n)$ such that
\ben
\lim_{t\to+\infty}\|\psi-\psi^*\|_{H^{l}(\mathbb{R}^n)}=0,\quad\quad\lim_{t\to+\infty}\|\vec{v}-\vec{v}^*\|_{{\bf{H}}^{l}(\mathbb{R}^n)}=0,
\een
where $\vec{v}^*=(v^*,\partial_tv^*)'$, and $(\psi^*,v^*)$ is the solution to  
\be
\l\{\begin{split}
i\gz^\mu\partial_\mu\psi^*+M\psi^*&=0,\quad (t,x)\in[2,\infty)\times\mathbb{R}^n,\\
-\Box v^*+v^*&=0,\quad (t,x)\in[2,\infty)\times\mathbb{R}^n,
\end{split}\r.\quad\quad(\psi^*,v^*,\partial_tv^*)|_{t=2}=(\psi^*_{0},v^*_0,v^*_1).
\ee
\end{itemize}

\begin{proof}
$i)$ We only need to consider the case $l=0$. For any fixed $t\ge 2$, let $\vec{U}(\tau):=\mathcal{\tilde{S}}(t-\tau)(\vec{f}(\tau))=(U_1,U_2)'$. For any $4\le T_1<T_2<\infty$, by standard energy inequalities, we have
\ben
&&\l\|\int_{T_1}^{T_2}\vec{U}(\tau){\rm{d}}\tau\r\|_{{\bf{H}}^0(\mathbb{R}^n)}\lesssim\l(\int_{\mathbb{R}^n}\Bigg\{\l|\int_{T_1}^{T_2}\nabla U_1(\tau){\rm{d}}\tau\r|^2+\l|\int_{T_1}^{T_2}U_2(\tau){\rm{d}}\tau\r|^2+\l|\int_{T_1}^{T_2} U_1(\tau){\rm{d}}\tau\r|^2\Bigg\}{\rm{d}}x\r)^{\f{1}{2}}\\
&\lesssim&\l(\int_{\mathbb{R}^n}\l\{\int_{T_1}^{T_2}\big(|\nabla U_1|+|U_2|+|U_1|\big)^2(\tau,x)\cdot \tau^{1+\dz}{\rm{d}}\tau\r\}\cdot\int_{T_1}^{T_2}\tau^{-(1+\dz)}{\rm{d}}\tau{\rm{d}}x\r)^{\f{1}{2}}\\
&\lesssim&T_1^{-\f{\dz}{2}}\l(\int_{T_1}^{T_2}\int_{\mathbb{R}^n}\big(|\nabla U_1|+|U_2|+|U_1|\big)^2(\tau,x)\cdot \tau^{1+\dz}{\rm{d}}x{\rm{d}}\tau\r)^{\f{1}{2}}\\
&\lesssim&II:=T_1^{-\f{\dz}{2}}\l(\int_{T_1}^{T_2}\l(\int_{r\ge \tau-1}+\int_{r<\tau-1}\r)\big(|\nabla f_1|+|f_2|+|f_1|\big)^2(\tau,x)\cdot \tau^{1+\dz}{\rm{d}}x{\rm{d}}\tau\r)^{\f{1}{2}}.
\een
Denote $\tilde{f}:=|\nabla f_1|+|f_2|+|f_1|$. By a change of variables $(\tau,x)\to(s,x)$ with $s=\sqrt{\tau^2-|x|^2}$, we obtain
\beqn\label{s5: II}
II&\lesssim&T_1^{-\f{\dz}{2}}\l(\int_{T_1}^{T_2}\int_{\Sigma^{\rm{ex}}_\tau}|\tilde{f}|^2(\tau,x)\cdot \tau^{1+\dz}{\rm{d}}x{\rm{d}}\tau\r)^{\f{1}{2}}+T_1^{-\f{\dz}{2}}\l(\int_{T_1^{\f{1}{2}}}^{T_2}\int_{r<\f{s^2-1}{2}}|\tilde{f}|^2(\sqrt{s^2+|x|^2},x)\cdot \tau^{1+\dz}\f{s}{\tau}{\rm{d}}x{\rm{d}}s\r)^{\f{1}{2}}\nonumber\\
&\lesssim&T_1^{-\f{\dz}{2}}\Bigg\{\l(\int_{T_1}^{T_2}\|\tilde{f}\|_{L^2_x(\Sigma^{\rm{ex}}_\tau)}^2 \cdot\tau^{1+\dz}{\rm{d}}\tau\r)^{\f{1}{2}}+\l(\int_{T_1^{\f{1}{2}}}^{T_2}\|\tilde{f}\|_{L^2(\mathscr{H}^{\rm{in}}_s)}^2\cdot s^{1+2\dz}{\rm{d}}s\r)^{\f{1}{2}}\Bigg\}.
\eeqn
Similarly,
\ben
\l\|\int_{T_1}^{T_2}\mathcal{S}(t-\tau)(g(\tau)){\rm{d}}\tau\r\|_{L^2(\mathbb{R}^n)}\lesssim T_1^{-\f{\dz}{2}}\Bigg\{\l(\int_{T_1}^{T_2}\|g\|_{L^2_x(\Sigma^{\rm{ex}}_\tau)}^2 \cdot\tau^{1+\dz}{\rm{d}}\tau\r)^{\f{1}{2}}+\l(\int_{T_1^{\f{1}{2}}}^{T_2}\|g\|_{L^2(\mathscr{H}^{\rm{in}}_s)}^2\cdot s^{1+2\dz}{\rm{d}}s\r)^{\f{1}{2}}\Bigg\}.
\een
$ii)$ Let 
\ben
\psi^*_{0}:&=&\psi_0-i\int_2^{+\infty}\mathcal{S}(2-\tau)\gz^0F_\psi(\tau){\rm{d}}\tau,\quad\quad \psi^*=\mathcal{S}(t-2)\psi^*_{0},\\
\vec{v}^*_0:&=&\vec{v}_0+\int_2^{+\infty}\mathcal{\tilde{S}}(2-\tau)\vec{F}_v(\tau){\rm{d}}\tau,\quad\quad \vec{v}^*=\mathcal{\tilde{S}}(t-2)\vec{v}^*_0.
\een
For any $4\le T_1<T_2<\infty$, by $i)$, we have
\be
\begin{split}
&\quad\quad\l\|\int_{T_1}^{T_2}\mathcal{S}(2-\tau)\gz^0F_\psi(\tau){\rm{d}}\tau\r\|_{H^{l}(\mathbb{R}^n)}+\l\|\int_{T_1}^{T_2}\mathcal{\tilde{S}}(2-\tau)\vec{F}_v(\tau){\rm{d}}\tau\r\|_{{\bf{H}}^{l}(\mathbb{R}^n)}\\
&\lesssim T_1^{-\f{\dz}{2}}\Bigg\{\l(\int_{T_1}^{T_2}\|F\|_{L^2_x(\Sigma^{\rm{ex}}_\tau)}^2 \cdot\tau^{1+\dz}{\rm{d}}\tau\r)^{\f{1}{2}}+\l(\int_{T_1^{\f{1}{2}}}^{T_2}\|F\|_{L^2(\mathscr{H}^{\rm{in}}_s)}^2\cdot s^{1+2\dz}{\rm{d}}s\r)^{\f{1}{2}}\Bigg\}\lesssim AT_1^{-\f{\dz}{2}}\to 0\quad\mathrm{as}\ \  T_1\to+\infty,
\end{split}
\ee
where $F:=\sum_{k=0}^{l}\big(|\nabla^{k}F_\psi|+|\nabla^kF_v|\big)$. Hence, $\psi^*_{0}$ and $\vec{v}^*_0$ are well-defined in $H^{l}(\mathbb{R}^n)$ and ${\bf{H}}^{l}(\mathbb{R}^n)$ respectively. Similarly, using $i)$ again, we have
\ben
&&\|\psi-\psi^*\|_{H^{l}(\mathbb{R}^n)}=\l\|\int_t^{+\infty}\mathcal{S}(t-\tau)\gz^0F_\psi(\tau){\rm{d}}\tau\r\|_{H^{l}(\mathbb{R}^n)}\lesssim At^{-\f{\dz}{2}}\to 0\quad\mathrm{as}\ \ t\to+\infty,\\
&&\|\vec{v}-\vec{v}^*\|_{{\bf{H}}^{l}(\mathbb{R}^n)}=\l\|\int_t^{+\infty}\mathcal{\tilde{S}}(t-\tau)\vec{F}_v(\tau){\rm{d}}\tau\r\|_{{\bf{H}}^{l}(\mathbb{R}^n)}\lesssim At^{-\f{\dz}{2}}\to 0\quad\mathrm{as}\ \ t\to+\infty.
\een
The proof is done.
\end{proof}
\end{lem}

\section{$n=2$: Global existence and estimates in the exterior region}\label{s3}

In this section, we prove global existence of the solution to \eqref{s1: DKG}-\eqref{s1: initial} in the exterior region.

\subsection{Bootstrap setting}\label{s3.1}

Fix $N\ge 6$, $\lz>0$ and $0<\dz\ll \lz_0$, where
\beq\label{s3: deflz_0}
\lz_0:=\f{1}{2}\min\{\lz,1\}.
\eeq
Let $C_1\gg 1$ and $0<\ez\ll C_1^{-1}$ be chosen later. We introduce the following bootstrap setting for the solution $(\psi,v)$ to \eqref{s1: DKG}-\eqref{s1: initial}:
\beq\label{s3: bspsiv}
\sum_{|I|\le N}\big\{\|\hat{\Gamma}^I\psi\|_{Y^{{\rm{ex}},\lz}_{D,t}}+\|\Gamma^Iv\|_{Y^{{\rm{ex}},\lz}_{1,\dz,t}}\big\}+\sum_{|I|\le N-1}\|\Gamma^Iv\|_{Y^{{\rm{ex}},\lz}_{1,0,t}}\le C_1\ez,
\eeq
where (see \eqref{s2: exnorm} and \eqref{s2: exnormD})
\beq\label{s3: Yhgzpsi0}
\|\hat{\Gamma}^I\psi\|_{Y^{{\rm{ex}},\lz}_{D,t}}:=[E^{{\rm{ex}},\lz}_D(t,\hat{\Gamma}^I\psi)]^{\f{1}{2}}+\l(\int_2^t\int_{\Sigma^{\rm{ex}}_\tau}(2+r-\tau)^\lz|[\hat{\Gamma}^I\psi]_{-}|^2{\rm{d}}x{\rm{d}}\tau\r)^{\f{1}{2}},
\eeq
\beq\label{s3: Ygzvdz}
\|\Gamma^Iv\|_{Y^{{\rm{ex}},\lz}_{1,\dz,t}}:=[E^{{\rm{ex}},\lz}_{1,\dz}(t,\Gamma^Iv)]^{\f{1}{2}}+\l(\int_2^t\int_{\Sigma^{\rm{ex}}_\tau}\tau^{-\dz}(2+r-\tau)^\lz\big(|G\Gamma^Iv|^2+|\Gamma^Iv|^2\big){\rm{d}}x{\rm{d}}\tau\r)^{\f{1}{2}},
\eeq
\beqn\label{s3: Ygzv0}
\|\Gamma^Iv\|_{Y^{{\rm{ex}},\lz}_{1,0,t}}:=[E^{{\rm{ex}},\lz}_{1}(t,\Gamma^Iv)]^{\f{1}{2}}+\l(\int_2^t\int_{\Sigma^{\rm{ex}}_\tau}(2+r-\tau)^\lz\big(|G\Gamma^Iv|^2+|\Gamma^Iv|^2\big){\rm{d}}x{\rm{d}}\tau\r)^{\f{1}{2}},
\eeqn
and we recall \eqref{s2: G_a} and \eqref{s2: psi-} for the definitions of $|G\Gamma^I v|$ and $[\hat{\Gamma}^I\psi]_{-}$. We also recall \eqref{dexfD}, \eqref{dexfs} and \eqref{dexfsaz} for the definitions of the energy functionals above. Let
\beq\label{s3: maxt}
T^*:=\sup\{T>2: \eqref{s3: bspsiv}\ \mathrm{holds\ for}\ t\in[2,T]\}.
\eeq

\begin{prop}\label{s3: maxT}
There exist some constants $C_1>0$ sufficiently large and $0<\ez_0\ll C_1^{-1}$ sufficiently small such that, for any $0<\ez<\ez_0$, if $(\psi,v)$ is a solution to \eqref{s1: DKG}-\eqref{s1: initial} in a time interval $[2,T]$ and satisfies \eqref{s3: bspsiv}, then for $t\in[2,T]$ we have
\be\label{s3: imbsex}
\sum_{|I|\le N}\big\{\|\hat{\Gamma}^I\psi\|_{Y^{{\rm{ex}},\lz}_{D,t}}+\|\Gamma^Iv\|_{Y^{{\rm{ex}},\lz}_{1,\dz,t}}\big\}+\sum_{|I|\le N-1}\|\Gamma^Iv\|_{Y^{{\rm{ex}},\lz}_{1,0,t}}\le \f{1}{2}C_1\ez.
\ee
\end{prop}

In the above proposition $T$ is arbitrary, hence the solution $(\psi,v)$ exists globally in time in $\mathscr{D}^{\rm{ex}}$ (i.e., $T^*=\infty$ where $T^*$ is as in \eqref{s3: maxt}) and satisfies \eqref{s3: bspsiv} for all $t\in[2,\infty)$. Below we give the proof of Proposition \ref{s3: maxT}. In the sequel, the implied constants in $\lesssim$ do not depend on the constants $C_1$ and $\ez$ appearing in the bootstrap assumption \eqref{s3: bspsiv}. Let 
\beqn\label{s3: l(tau)}
l(\tau):&=&(C_1\ez)^{-2}\Bigg\{\sum_{|I|\le N}\int_{\Sigma^{\rm{ex}}_\tau}(2+r-\tau)^\lz\big(|[\hat{\Gamma}^I\psi]_{-}|^2+\tau^{-\dz}|\Gamma^Iv|^2\big){\rm{d}}x\nonumber\\
&&\quad\quad\quad+\sum_{|I|\le N-1}\int_{\Sigma^{\rm{ex}}_\tau}(2+r-\tau)^\lz|\Gamma^Iv|^2{\rm{d}}x\Bigg\}.
\eeqn
Then 
\beq\label{s3: l(tau)L1}
l(\tau)\in L^1([2,t]),\quad\quad \|l\|_{L^1([2,t])}\lesssim 1
\eeq
uniformly in $t\in[2,T]$. By \eqref{s3: bspsiv}, we obtain the following $L^2$-type estimates for the solution $(\psi,v)$ on the time interval $[2,T]$:
\beq\label{s3: L^2psiv}
\l\{\begin{split}
&\sum_{|I|\le N}\|(2+r-t)^{\f{1+\lz}{2}}\big\{|\hat{\Gamma}^I\psi|+t^{-\f{\dz}{2}}\big(|\partial\Gamma^Iv|+|\Gamma^Iv|\big)\big\}\|_{L^2_x(\Sigma^{\rm{ex}}_t)}\lesssim C_1\ez,\\
&\sum_{|I|\le N-1}\|(2+r-t)^{\f{1+\lz}{2}}\big(|\partial\Gamma^Iv|+|\Gamma^Iv|\big)\|_{L^2_x(\Sigma^{\rm{ex}}_t)}\lesssim C_1\ez,\\
&\sum_{|I|\le N}\|(2+r-\tau)^{\f{\lz}{2}}\big(|[\hat{\Gamma}^I\psi]_{-}|+\tau^{-\f{\dz}{2}}|\Gamma^Iv|\big)\|_{L^2_x(\Sigma^{\rm{ex}}_\tau)}\lesssim C_1\ez\sqrt{l(\tau)},\quad\tau\in[2,t],\\
&\sum_{|I|\le N-1}\|(2+r-\tau)^{\f{\lz}{2}}\Gamma^Iv\|_{L^2_x(\Sigma^{\rm{ex}}_\tau)}\lesssim C_1\ez\sqrt{l(\tau)},\quad\tau\in[2,t].
\end{split}\r.
\eeq
We recall \eqref{s2: hatLO} and \eqref{s2: pOL^J}. We claim that $[\partial_r,(x_a/r)\gz^0\gz^a]=[\hat{\Omega}_{12},(x_a/r)\gz^0\gz^a]=0$. Indeed, we have
\be
\begin{split}
\l[\f{x^b}{r}\partial_b,\f{x_a}{r}\gz^0\gz^a\r]&=\f{x^b}{r}\l(\f{\dz_{ab}}{r}-\f{x_ax_b}{r^3}\r)\gz^0\gz^a=0,\\
\l[\hat{\Omega}_{12},\f{x_a}{r}\gz^0\gz^a\r]&=\l[x_1\partial_2-x_2\partial_1-\f{1}{2}\gz^1\gz^2,\f{x_a}{r}\gz^0\gz^a\r]\\
&=x_1\l(\f{\dz_{a2}}{r}-\f{x_ax_2}{r^3}\r)\gz^0\gz^a-x_2\l(\f{\dz_{a1}}{r}-\f{x_ax_1}{r^3}\r)\gz^0\gz^a-\f{1}{2}\f{x_a}{r}\l(\gz^1\gz^2\gz^0\gz^a-\gz^0\gz^a\gz^1\gz^2\r)\\
&=\f{x_1}{r}\gz^0\gz^2-\f{x_2}{r}\gz^0\gz^1-\f{1}{2}\f{x_a}{r}\l(\gz^0\gz^1(-2\dz^{a}_2-\gz^a\gz^2)-\gz^0(-2\dz^{a}_1-\gz^1\gz^a)\gz^2\r)=0,
\end{split}
\ee
hence the claim holds. Then by \eqref{s3: L^2psiv}, we have 
\beqn\label{s3: L^2psi-}
&&\sum_{|I|\le N-2, |J|\le 1}\big\{\|(2+r-\tau)^{\f{\lz}{2}}\hat{\Omega}^J[\hat{\Gamma}^I\psi]_{-}\|_{L^2_x(\Sigma^{\rm{ex}}_\tau)}+\|(2+r-\tau)^{\f{\lz}{2}}\partial_r\hat{\Omega}^J[\hat{\Gamma}^I\psi]_{-}\|_{L^2_x(\Sigma^{\rm{ex}}_\tau)}\big\}\nonumber\\
&\lesssim&\sum_{|I|\le N-2, |J|\le 1}\big\{\|(2+r-\tau)^{\f{\lz}{2}}[\hat{\Omega}^J\hat{\Gamma}^I\psi]_{-}\|_{L^2_x(\Sigma^{\rm{ex}}_\tau)}+\|(2+r-\tau)^{\f{\lz}{2}}[\partial_r\hat{\Omega}^J\hat{\Gamma}^I\psi]_{-}\|_{L^2_x(\Sigma^{\rm{ex}}_\tau)}\big\}\nonumber\\
&\lesssim&\sum_{|I|\le N}\|(2+r-\tau)^{\f{\lz}{2}}[\hat{\Gamma}^I\psi]_{-}\|_{L^2_x(\Sigma^{\rm{ex}}_\tau)}\lesssim C_1\ez\sqrt{l(\tau)},\quad\tau\in[2,t].
\eeqn
Using \eqref{s3: L^2psiv}, \eqref{s3: L^2psi-} and \eqref{s2: Soeq0} in Lemma \ref{s2: Sobex}, we obtain the following pointwise estimates for the solution $(\psi,v)$ on $[2,T]$:
\beq\label{s3: pointpsiv}
\l\{\begin{split}
&\sum_{|I|\le N-3}\sup_{\Sigma^{\rm{ex}}_t}\ (2+r-t)^{\f{1+\lz}{2}} r^{\f{1}{2}}\big(|\hat{\Gamma}^I\psi|+|\partial\Gamma^Iv|+|\Gamma^Iv|\big)\lesssim C_1\ez,\\
&\sum_{|I|\le N-3}\sup_{\Sigma^{\rm{ex}}_\tau}\ (2+r-\tau)^{\f{\lz}{2}} r^{\f{1}{2}}\big(|[\hat{\Gamma}^I\psi]_{-}|+|\Gamma^Iv|\big)\lesssim C_1\ez\sqrt{l(\tau)},\quad\tau\in[2,t].
\end{split}\r.
\eeq
We note that on $\Sigma^{\rm{ex}}_t$ we have
\beq\label{s3: 2+r-t}
(2+r-t)^{\f{1-\lz}{2}}\lesssim\max\{r^{\f{1-\lz}{2}},1\}.
\eeq
By \eqref{s3: pointpsiv}, we obtain 
\beq\label{s3: pointpsiv0*}
\l\{\begin{split}
&\sum_{|I|\le N-3}\sup_{\Sigma^{\rm{ex}}_t}(2+r-t) r^{\lz_0}\big(|\hat{\Gamma}^I\psi|+|\partial\Gamma^Iv|+|\Gamma^Iv|\big)\lesssim C_1\ez,\\
&\sum_{|I|\le N-3}\sup_{\Sigma^{\rm{ex}}_\tau}(2+r-\tau)^{\f{1}{2}}r^{\lz_0}\big(|[\hat{\Gamma}^I\psi]_{-}|+|\Gamma^Iv|\big)\lesssim C_1\ez\sqrt{l(\tau)},\quad \tau\in[2,t],
\end{split}\r.
\eeq
where $\lz_0$ is as in \eqref{s3: deflz_0}. By Lemma \ref{s2: Edecay}, \eqref{s3: pointpsiv} and \eqref{s3: pointpsiv0*}, on $\Sigma^{\rm{ex}}_t\cap\{r\le 3t\}$, we have
\beqn\label{s3: pointve*}
\sum_{|I|\le N-4}|\Gamma^Iv|&\lesssim&\f{2+r-t}{t}\sum_{|I|\le N-3}|\partial\Gamma^Iv|+\sum_{|I|\le N-4}|\Gamma^I(\psi^*\gz^0\psi)|\nonumber\\
&\lesssim&\f{2+r-t}{t}\sum_{|I|\le N-3}|\partial\Gamma^Iv|+\sum_{|I_1|+|I_2|\le N-4}|\hat{\Gamma}^{I_1}\psi|\cdot|\hat{\Gamma}^{I_2}\psi|\nonumber\\
&\lesssim&C_1\ez r^{-1-\lz_0}+(C_1\ez)^2(2+r-t)^{-\lz-1}r^{-1}\lesssim C_1\ez r^{-1}.
\eeqn
By \eqref{s2: 2.2-42} and \eqref{s2: 2.2-43'} in Lemma \ref{sC: 2.2-4}, for any $\tau\in[2,t]$, on $\Sigma^{\rm{ex}}_\tau\cap\{r\le 3\tau\}$ we have 
\beqn\label{s3: pointve}
\sum_{|I|\le N-4}|\Gamma^Iv|&\lesssim&\f{2+r-\tau}{\tau}\sum_{|I|\le N-3}|\partial\Gamma^Iv|+\sum_{|I|\le N-4}|\Gamma^I(\psi^*\gz^0\psi)|\nonumber\\
&\lesssim&\f{2+r-\tau}{\tau}\sum_{|I|\le N-3}|\partial\Gamma^Iv|+\sum_{|I_1|+|I_2|\le N-4}|[\hat{\Gamma}^{I_1}\psi]_{-}|\cdot|\hat{\Gamma}^{I_2}\psi|\nonumber\\
&\lesssim&C_1\ez  r^{-1-\lz_0}+(C_1\ez)^2(2+r-\tau)^{-\lz-\f{1}{2}} \sqrt{l(\tau)}r^{-1}\nonumber\\
&\lesssim&C_1\ez \big\{r^{-1-\lz_0}+(2+r-\tau)^{-\f{1}{2}} \sqrt{l(\tau)}r^{-1}\big\}.
\eeqn
On the other hand, on $\Sigma^{\rm{ex}}_\tau\cap\{r\ge 3\tau\}$ we have $2+r-\tau\sim r$, hence by \eqref{s3: pointpsiv}, we have
\beq\label{s3: pointvf}
|\Gamma^Iv|\lesssim C_1\ez r^{-1-\f{\lz}{2}},\quad |I|\le N-4.
\eeq
Combining \eqref{s3: pointve*}, \eqref{s3: pointve} and \eqref{s3: pointvf}, we obtain
\beq\label{s3: exvextra*}
|\Gamma^Iv|\lesssim C_1\ez r^{-1}\quad\quad\mathrm{on}\ \ \Sigma^{\rm{ex}}_t,\ \mathrm{for}\ |I|\le N-4,
\eeq
\beq\label{s3: exvextra}
|\Gamma^Iv|\lesssim C_1\ez \big\{r^{-1-\lz_0}+(2+r-\tau)^{-\f{1}{2}} \sqrt{l(\tau)}r^{-1}\big\}\quad\quad\mathrm{on}\ \ \Sigma^{\rm{ex}}_\tau, \ \mathrm{for}\ \tau\in[2,t],\ |I|\le N-4.
\eeq

\subsection{Improved estimates for the solution $(\psi,v)$ in the exterior region}\label{s3.2}

${\bf Step\ 1.}$ First, we refine the energy estimate of $\psi$. Acting the vector fields $\hat{\Gamma}^I, |I|\le N$ on both sides of the first equation in \eqref{s1: DKG} and applying Proposition \ref{s2: exfsD}, we obtain
\be
\|\hat{\Gamma}^I\psi\|_{Y^{{\rm{ex}},\lz}_{D,t}}^2\lesssim\|\hat{\Gamma}^I\psi\|_{Y^{{\rm{ex}},\lz}_{D,2}}^2+\int_2^t\|(2+r-\tau)^{1+\lz}(\hat{\Gamma}^I\psi)^*\gz^0\hat{\Gamma}^I(v\psi)\|_{L^1_x(\Sigma^{\rm{ex}}_\tau)}{\rm{d}}\tau,
\ee
where we recall \eqref{s3: Yhgzpsi0}. For $|I|\le N$ and $\tau\in[2,t]$, by \eqref{s2: 2.2-43'} and \eqref{s2: 2.2-41} in Lemma \ref{sC: 2.2-4}, we have
\beqn\label{s3: psiNL^2}
&&\|(2+r-\tau)^{1+\lz}(\hat{\Gamma}^I\psi)^*\gz^0\hat{\Gamma}^I(v\psi)\|_{L^1_x(\Sigma^{\rm{ex}}_\tau)}\nonumber\\
&\lesssim&\sum_{|I_1|+|I_2|\le N}\|(2+r-\tau)^{1+\lz}\big\{|[\hat{\Gamma}^I\psi]_{-}|\cdot|\Gamma^{I_1}v|\cdot|\hat{\Gamma}^{I_2}\psi|+|\hat{\Gamma}^I\psi|\cdot|\Gamma^{I_1}v|\cdot|[\hat{\Gamma}^{I_2}\psi]_{-}|\big\}\|_{L^1_x(\Sigma^{\rm{ex}}_\tau)}\nonumber\\
&\lesssim&\sum_{\substack{|I_1|\le N-4\\|I'|, |I_2|\le N}}\|(2+r-\tau)^{\f{\lz}{2}}[\hat{\Gamma}^{I'}\psi]_{-}\|_{L^2_x(\Sigma^{\rm{ex}}_\tau)}\cdot\|(2+r-\tau)^{\f{1}{2}}\Gamma^{I_1}v\|_{L^\infty_x(\Sigma^{\rm{ex}}_\tau)}\cdot\|(2+r-\tau)^{\f{1+\lz}{2}}\hat{\Gamma}^{I_2}\psi\|_{L^2_x(\Sigma^{\rm{ex}}_\tau)}\nonumber\\
&+&\sum_{\substack{|I_2|\le N-3\\|I'|,|I_1|\le N}}\|(2+r-\tau)^{\f{\lz}{2}}[\hat{\Gamma}^{I'}\psi]_{-}\|_{L^2_x(\Sigma^{\rm{ex}}_\tau)}\cdot\|(2+r-\tau)^{\f{\lz}{2}}\Gamma^{I_1}v\|_{L^2_x(\Sigma^{\rm{ex}}_\tau)}\cdot\|(2+r-\tau)\hat{\Gamma}^{I_2}\psi\|_{L^\infty_x(\Sigma^{\rm{ex}}_\tau)}\nonumber\\
&+&\sum_{\substack{|I_2|\le N-3\\|I'|,|I_1|\le N}}\|(2+r-\tau)^{\f{1+\lz}{2}}\hat{\Gamma}^{I'}\psi\|_{L^2_x(\Sigma^{\rm{ex}}_\tau)}\cdot\|(2+r-\tau)^{\f{\lz}{2}}\Gamma^{I_1}v\|_{L^2_x(\Sigma^{\rm{ex}}_\tau)}\cdot\|(2+r-\tau)^{\f{1}{2}}[\hat{\Gamma}^{I_2}\psi]_{-}\|_{L^\infty_x(\Sigma^{\rm{ex}}_\tau)}\nonumber\\
&\lesssim&(C_1\ez)^3\sqrt{l(\tau)}\big\{\tau^{-\f{1}{2}-\lz_0}+\sqrt{l(\tau)}\tau^{-1}\big\}+(C_1\ez)^3 l(\tau)\tau^{-\lz_0+\f{\dz}{2}},
\eeqn
where we use \eqref{s3: L^2psiv}, \eqref{s3: exvextra} and \eqref{s3: pointpsiv0*}. Recall that $0<\dz\ll \lz_0$. Hence by \eqref{s3: l(tau)L1}, we have
\beq\label{s3: S1psi}
\sum_{|I|\le N}\|\hat{\Gamma}^I\psi\|_{Y^{{\rm{ex}},\lz}_{D,t}}\lesssim \ez+(C_1\ez)^{\f{3}{2}}.
\eeq

${\bf Step\ 2.}$ We now refine the highest order energy estimate of $v$. Acting the vector fields $\Gamma^I, |I|\le N$ on both sides of the second equation in \eqref{s1: DKG} and applying Proposition \ref{s2: exfs}, we obtain
\be
\|\Gamma^Iv\|_{Y^{{\rm{ex}},\lz}_{1,\dz,t}}\lesssim\|\Gamma^Iv\|_{Y^{{\rm{ex}},\lz}_{1,\dz,2}}+\int_2^t\|\tau^{-\f{\dz}{2}}(2+r-\tau)^{\f{1+\lz}{2}}\Gamma^I(\psi^*\gz^0\psi)\|_{L^2_x(\Sigma^{\rm{ex}}_\tau)}{\rm{d}}\tau,
\ee
where we recall \eqref{s3: Ygzvdz}. For $|I|\le N$ and $\tau\in[2,t]$, by \eqref{s2: 2.2-42} and \eqref{s2: 2.2-43'} in Lemma \ref{sC: 2.2-4}, we have
\beqn\label{s3: vNL^2}
&&\!\!\!\|\tau^{-\f{\dz}{2}}(2+r-\tau)^{\f{1+\lz}{2}}\Gamma^I(\psi^*\gz^0\psi)\|_{{L^2_x(\Sigma^{\rm{ex}}_\tau)}}\lesssim\sum_{|I_1|+|I_2|\le N}\|\tau^{-\f{\dz}{2}}(2+r-\tau)^{\f{1+\lz}{2}}|[\hat{\Gamma}^{I_1}\psi]_{-}|\cdot|\hat{\Gamma}^{I_2}\psi|\,\|_{L^2_x(\Sigma^{\rm{ex}}_\tau)}\nonumber\\
&\lesssim&\sum_{\substack{|I_1|\le N-3\\|I_2|\le N}}\|\tau^{-\f{\dz}{2}}[\hat{\Gamma}^{I_1}\psi]_{-}\|_{L^\infty_x(\Sigma^{\rm{ex}}_\tau)}\cdot\|(2+r-\tau)^{\f{1+\lz}{2}}\hat{\Gamma}^{I_2}\psi\|_{L^2_x(\Sigma^{\rm{ex}}_\tau)}\nonumber\\
&+&\sum_{\substack{|I_2|\le N-3\\|I_1|\le N}}\|[\hat{\Gamma}^{I_1}\psi]_{-}\|_{L^2_x(\Sigma^{\rm{ex}}_\tau)}\cdot\|\tau^{-\f{\dz}{2}}(2+r-\tau)^{\f{1+\lz}{2}}\hat{\Gamma}^{I_2}\psi\|_{L^\infty_x(\Sigma^{\rm{ex}}_\tau)}\lesssim(C_1\ez)^2\sqrt{l(\tau)}\tau^{-\f{1}{2}-\f{\dz}{2}},
\eeqn
where we use \eqref{s3: L^2psiv} and \eqref{s3: pointpsiv}. Hence \eqref{s3: l(tau)L1} gives
\beq\label{s3: S2v}
\sum_{|I|\le N}\|\Gamma^Iv\|_{Y^{{\rm{ex}},\lz}_{1,\dz,t}}\lesssim\ez+(C_1\ez)^2.
\eeq

${\bf Step\ 3.}$ Next we refine the lower order energy estimates of $v$. Let $\tilde{v}:=v-\psi^*\gz^0\psi$. For $|I|\le N-1$, by $ii)$ in Lemma \ref{s2: nlt} and Proposition \ref{s2: exfs}, we obtain 
\be
\|\Gamma^I\tilde{v}\|_{Y^{{\rm{ex}},\lz}_{1,0,t}}\lesssim\|\Gamma^I\tilde{v}\|_{Y^{{\rm{ex}},\lz}_{1,0,2}}+\int_2^t\|(2+r-\tau)^{\f{1+\lz}{2}}\Gamma^I\tilde{F}_v\|_{L^2_x(\Sigma^{\rm{ex}}_\tau)}{\rm{d}}\tau,
\ee
where 
\beq\label{s3: exnvN-3}
\|\Gamma^I\tilde{v}\|_{Y^{{\rm{ex}},\lz}_{1,0,t}}=[E^{{\rm{ex}},\lz}_1(t,\Gamma^I\tilde{v})]^{\f{1}{2}}+\l(\int_2^t\int_{\Sigma^{\rm{ex}}_\tau}(2+r-\tau)^\lz\big(|G\Gamma^I\tilde{v}|^2+|\Gamma^I\tilde{v}|^2\big){\rm{d}}x{\rm{d}}\tau\r)^{\f{1}{2}}
\eeq
and (see \eqref{s2: tildeF})
\beq\label{s3: tildeF}
\tilde{F}_v:=-i\partial_\mu(v\psi^*)\gz^0\gz^\mu\psi+i\psi^*\gz^0\gz^\mu\partial_\mu(v\psi)+2\partial_\az\psi^*\gz^0\partial^\az\psi.
\eeq
Here we recall \eqref{dexfsaz} for the definition of $E^{{\rm{ex}},\lz}_1(t,\Gamma^I\tilde{v})$. For $|I|\le N-1$ and $\tau\in[2,t]$, by \eqref{s3: L^2psiv} and \eqref{s3: pointpsiv}, we have
\beqn\label{s3: tildevF1}
&&\|(2+r-\tau)^{\f{1+\lz}{2}}\Gamma^I\big(\partial_\mu(v\psi^*)\gz^0\gz^\mu\psi\big)\|_{L^2_x(\Sigma^{\rm{ex}}_\tau)}\nonumber\\
&\lesssim&\sum_{|I_1|+|I_2|+|I_3|\le N}\|(2+r-\tau)^{\f{1+\lz}{2}}|\Gamma^{I_1}v|\cdot|\Gamma^{I_2}\psi|\cdot|\Gamma^{I_3}\psi|\,\|_{L^2_x(\Sigma^{\rm{ex}}_\tau)}\nonumber\\
&\lesssim&\sum_{\substack{|I_1|,|I_2|\le N-3\\|I_3|\le N}}\|\Gamma^{I_1}v\|_{L^\infty_x(\Sigma^{\rm{ex}}_\tau)}\cdot\|\Gamma^{I_2}\psi\|_{L^\infty_x(\Sigma^{\rm{ex}}_\tau)}\cdot\|(2+r-\tau)^{\f{1+\lz}{2}}\Gamma^{I_3}\psi\|_{L^2_x(\Sigma^{\rm{ex}}_\tau)}\nonumber\\
&+&\sum_{\substack{|I_2|,|I_3|\le N-3\\|I_1|\le N}}\|\Gamma^{I_1}v\|_{L^2_x(\Sigma^{\rm{ex}}_\tau)}\cdot\|\Gamma^{I_2}\psi\|_{L^\infty_x(\Sigma^{\rm{ex}}_\tau)}\cdot\|(2+r-\tau)^{\f{1+\lz}{2}}\Gamma^{I_3}\psi\|_{L^\infty_x(\Sigma^{\rm{ex}}_\tau)}\nonumber\\
&\lesssim&(C_1\ez)^3\sqrt{l(\tau)}\tau^{-1+\f{\dz}{2}}.
\eeqn
For $|I|\le N-1$ and $\tau\in[2,t]$, by Lemma \ref{s2: Q_0}, we obtain
\beqn\label{s3: tildvF2}
&&\!\!\!\!\|(2+r-\tau)^{\f{1+\lz}{2}}\Gamma^I\big(\partial_\az\psi^*\gz^0\partial^\az\psi\big)\|_{L^2_x(\Sigma^{\rm{ex}}_\tau)}\lesssim\sum_{|I_1|+|I_2|\le N-1}\|(2+r-\tau)^{\f{1+\lz}{2}}(\partial_\az\Gamma^{I_1}\psi)^*\gz^0\partial^\az\Gamma^{I_2}\psi\|_{L^2_x(\Sigma^{\rm{ex}}_\tau)}\nonumber\\
&\lesssim&\sum_{|I_1|+|I_2|\le N-1}\sum_{a=1}^2\|\tau^{-1}(2+r-\tau)^{\f{1+\lz}{2}}\big(|L_a\Gamma^{I_1}\psi|\cdot|\partial\Gamma^{I_2}\psi|+|\tau-r|\cdot|\partial\Gamma^{I_1}\psi|\cdot|\partial\Gamma^{I_2}\psi|\big)\|_{L^2_x(\Sigma^{\rm{ex}}_\tau)}\nonumber\\
&\lesssim&\sum_{\substack{|I_1|+|I_2|\le N-1\\|J|=|J'|=1}}\|\tau^{-1}(2+r-\tau)^{1+\f{1+\lz}{2}}|\Gamma^J\Gamma^{I_1}\psi|\cdot|\Gamma^{J'}\Gamma^{I_2}\psi|\,\|_{L^2_x(\Sigma^{\rm{ex}}_\tau)}\nonumber\\
&\lesssim&\sum_{\substack{|I_1|\le N-4\\|J|=1,|I_2|\le N}}\|\tau^{-1}(2+r-\tau)\Gamma^J\Gamma^{I_1}\psi\|_{L^\infty_x(\Sigma^{\rm{ex}}_\tau)}\cdot\|(2+r-\tau)^{\f{1+\lz}{2}}\Gamma^{I_2}\psi\|_{L^2_x(\Sigma^{\rm{ex}}_\tau)}\lesssim(C_1\ez)^2\tau^{-1-\lz_0},
\eeqn
where we use \eqref{s2: Edecay0}, \eqref{s3: L^2psiv} and \eqref{s3: pointpsiv0*}. The estimates \eqref{s3: tildevF1} and \eqref{s3: tildvF2} yield 
\beq\label{s3: stildevF}
\|(2+r-\tau)^{\f{1+\lz}{2}}\Gamma^I\tilde{F}_v\|_{L^2_x(\Sigma^{\rm{ex}}_\tau)}\lesssim(C_1\ez)^2\big\{\sqrt{l(\tau)}\tau^{-1+\f{\dz}{2}}+\tau^{-1-\lz_0}\big\},\quad |I|\le N-1,\tau\in[2,t],
\eeq
which together with \eqref{s3: l(tau)L1} implies
\beq\label{s3: S3tildev}
\sum_{|I|\le N-1}\|\Gamma^I\tilde{v}\|_{Y^{{\rm{ex}},\lz}_{1,0,t}}\lesssim\ez+(C_1\ez)^2.
\eeq
We recall that $\tilde{v}=v-\psi^*\gz^0\psi$. For $|I|\le N-1$, using \eqref{dexfsaz},  \eqref{s3: L^2psiv} and \eqref{s3: pointpsiv}, we have
\ben
[E^{{\rm{ex}},\lz}_1(t,\Gamma^I(\psi^*\gz^0\psi))]^{\f{1}{2}}&\lesssim&\sum_{|I'|\le N}\|(2+r-t)^{\f{1+\lz}{2}}\Gamma^{I'}(\psi^*\gz^0\psi)\|_{L^2_x(\Sigma^{\rm{ex}}_t)}\\
&\lesssim&\sum_{|I_1|\le N-3,|I_2|\le N}\|\hat{\Gamma}^{I_1}\psi\|_{L^\infty_x(\Sigma^{\rm{ex}}_t)}\cdot\|(2+r-t)^{\f{1+\lz}{2}}\hat{\Gamma}^{I_2}\psi\|_{L^2_x(\Sigma^{\rm{ex}}_t)}\lesssim (C_1\ez)^2.
\een
By \eqref{s2: 2.2-42} and \eqref{s2: 2.2-43'} in Lemma \ref{sC: 2.2-4}, we also have
\ben
&&\sum_{|I'|\le N}\l(\int_2^t\|(2+r-\tau)^{\f{\lz}{2}}\Gamma^{I'}(\psi^*\gz^0\psi)\|_{L^2_x(\Sigma^{\rm{ex}}_\tau)}^2{\rm{d}}\tau\r)^{\f{1}{2}}\nonumber\\
&\lesssim&\sum_{|I_1|+|I_2|\le N}\l(\int_2^t\|(2+r-\tau)^{\f{\lz}{2}}|[\hat{\Gamma}^{I_1}\psi]_{-}|\cdot|\hat{\Gamma}^{I_2}\psi|\,\|_{L^2_x(\Sigma^{\rm{ex}}_\tau)}^2{\rm{d}}\tau\r)^{\f{1}{2}}\nonumber\\
&\lesssim&\sum_{|I_1|\le N-3,|I_2|\le N}\l(\int_2^t\|(2+r-\tau)^{\f{\lz}{2}}[\hat{\Gamma}^{I_1}\psi]_{-}\|_{L^\infty_x(\Sigma^{\rm{ex}}_\tau)}^2\cdot\|\hat{\Gamma}^{I_2}\psi\|_{L^2_x(\Sigma^{\rm{ex}}_\tau)}^2{\rm{d}}\tau\r)^{\f{1}{2}}\nonumber\\
&+&\sum_{|I_2|\le N-3,|I_1|\le N}\l(\int_2^t\|(2+r-\tau)^{\f{\lz}{2}}[\hat{\Gamma}^{I_1}\psi]_{-}\|_{L^2_x(\Sigma^{\rm{ex}}_\tau)}^2\cdot\|\hat{\Gamma}^{I_2}\psi\|_{L^\infty_x(\Sigma^{\rm{ex}}_\tau)}^2{\rm{d}}\tau\r)^{\f{1}{2}}\lesssim(C_1\ez)^2.
\een
Recall \eqref{s3: Ygzv0} and \eqref{s3: exnvN-3}. The last two estimates together with \eqref{s3: S3tildev} give
\beq\label{s3: gzv0im}
\sum_{|I|\le N-1}\|\Gamma^Iv\|_{Y^{{\rm{ex}},\lz}_{1,0,t}}\lesssim\ez+(C_1\ez)^2.
\eeq

Combining \eqref{s3: S1psi}, \eqref{s3: S2v} and \eqref{s3: gzv0im}, we have strictly improved the bootstrap estimate \eqref{s3: bspsiv} for $C_1$ sufficiently large and $0<\ez\ll C_1^{-1}$ sufficiently small. Hence the proof of Proposition \ref{s3: maxT} is completed. We conclude that the solution $(\psi,v)$ exists globally in the exterior region $\mathscr{D}^{\rm{ex}}$ and satisfies \eqref{s3: bspsiv} for all $t\in[2,\infty)$.

We next give an improved pointwise estimate of $[\hat{\Gamma}^I\psi]_{-}$ for $|I|\le N-4$. Let $\Psi$ solve \eqref{s2: boxPsi} and $\tilde{\Psi}:=\Psi-v\psi$. For $|I|\le N-1$, by $iii)$ in Lemma \ref{s2: nlt} and Proposition \ref{s2: exfs}, we obtain
\be
[E^{{\rm{ex}},\lz}_0(t,\Gamma^I\tilde{\Psi})]^{\f{1}{2}}\lesssim [E^{{\rm{ex}},\lz}_0(2,\Gamma^I\tilde{\Psi})]^{\f{1}{2}}+\int_2^t\|(2+r-\tau)^{\f{1+\lz}{2}}\Gamma^I\tilde{F}_\psi\|_{L^2_x(\Sigma^{\rm{ex}}_\tau)}{\rm{d}}\tau
\ee
for all $t\in[2,\infty)$, where we recall \eqref{dexfsaz} for the definition of $E^{{\rm{ex}},\lz}_0(t,\Gamma^I\tilde{\Psi})$ and (see \eqref{s2: boxtiPsi})
\beq\label{s3: tildeG}
\tilde{F}_\psi:=(\psi^*\gz^0\psi)\psi-iv\gz^\mu\partial_\mu(v\psi)-2\partial_\az v\partial^\az\psi
\eeq
For $|I|\le N-1$ and $\tau\in[2,t]$, by  \eqref{s2: 2.2-41}, \eqref{s2: 2.2-42} and \eqref{s2: 2.2-43'} in Lemma \ref{sC: 2.2-4}, we have
\beqn\label{s3: tildevG1}
&&\|(2+r-\tau)^{\f{1+\lz}{2}}\Gamma^I\big((\psi^*\gz^0\psi)\psi\big)\|_{L^2_x(\Sigma^{\rm{ex}}_\tau)}\nonumber\\
&\lesssim&\sum_{|I_1|+|I_2|+|I_3|\le N-1}\|(2+r-\tau)^{\f{1+\lz}{2}}|[\hat{\Gamma}^{I_1}\psi]_{-}|\cdot|\hat{\Gamma}^{I_2}\psi|\cdot|\hat{\Gamma}^{I_3}\psi|\,\|_{L^2_x(\Sigma^{\rm{ex}}_\tau)}\nonumber\\
&\lesssim&\sum_{\substack{|I_1|,|I_2|\le N-3\\|I_3|\le N-1}}\|[\hat{\Gamma}^{I_1}\psi]_{-}\|_{L^\infty_x(\Sigma^{\rm{ex}}_\tau)}\cdot\|\hat{\Gamma}^{I_2}\psi\|_{L^\infty_x(\Sigma^{\rm{ex}}_\tau)}\cdot\|(2+r-\tau)^{\f{1+\lz}{2}}\hat{\Gamma}^{I_3}\psi\|_{L^2_x(\Sigma^{\rm{ex}}_\tau)}\nonumber\\
&+&\sum_{\substack{|I_2|,|I_3|\le N-3\\|I_1|\le N-1}}\|[\hat{\Gamma}^{I_1}\psi]_{-}\|_{L^2_x(\Sigma^{\rm{ex}}_\tau)}\cdot\|(2+r-\tau)^{\f{1+\lz}{2}}\hat{\Gamma}^{I_2}\psi\|_{L^\infty_x(\Sigma^{\rm{ex}}_\tau)}\cdot\|\hat{\Gamma}^{I_3}\psi\|_{L^\infty_x(\Sigma^{\rm{ex}}_\tau)}\nonumber\\
&\lesssim&(C_1\ez)^3\sqrt{l(\tau)}\tau^{-1}
\eeqn
and
\beqn\label{s3: tildevG2}
&&\!\!\!\!\|(2+r-\tau)^{\f{1+\lz}{2}}\Gamma^I\big(v\gz^\mu\partial_\mu(v\psi)\big)\|_{L^2_x(\Sigma^{\rm{ex}}_\tau)}\lesssim\sum_{|I_1|+|I_2|+|I_3|\le N}\!\|(2+r-\tau)^{\f{1+\lz}{2}}|\Gamma^{I_1}v|\cdot|\Gamma^{I_2}v|\cdot|\Gamma^{I_3}\psi|\,\|_{L^2_x(\Sigma^{\rm{ex}}_\tau)}\nonumber\\
&\lesssim&\sum_{\substack{|I_1|,|I_2|\le N-3\\|I_3|\le N}}\|\Gamma^{I_1}v\|_{L^\infty_x(\Sigma^{\rm{ex}}_\tau)}\cdot\|\Gamma^{I_2}v\|_{L^\infty_x(\Sigma^{\rm{ex}}_\tau)}\cdot\|(2+r-\tau)^{\f{1+\lz}{2}}\Gamma^{I_3}\psi\|_{L^2_x(\Sigma^{\rm{ex}}_\tau)}\nonumber\\
&+&\sum_{\substack{|I_2|,|I_3|\le N-3\\|I_1|\le N}}\|\Gamma^{I_1}v\|_{L^2_x(\Sigma^{\rm{ex}}_\tau)}\cdot\|\Gamma^{I_2}v\|_{L^\infty_x(\Sigma^{\rm{ex}}_\tau)}\cdot\|(2+r-\tau)^{\f{1+\lz}{2}}\Gamma^{I_3}\psi\|_{L^\infty_x(\Sigma^{\rm{ex}}_\tau)}\nonumber\\
&\lesssim&(C_1\ez)^3\sqrt{l(\tau)}\tau^{-1+\f{\dz}{2}},
\eeqn
where we use \eqref{s3: L^2psiv} and \eqref{s3: pointpsiv}. For $|I|\le N-1$ and $\tau\in[2,t]$, by Lemma \ref{s2: Q_0}, we have
\beqn\label{s3: tildevG3}
&&\|(2+r-\tau)^{\f{1+\lz}{2}}\Gamma^I\big(\partial_\az v\partial^\az\psi\big)\|_{L^2_x(\Sigma^{\rm{ex}}_\tau)}\nonumber\\
&\lesssim&\sum_{\substack{|I_1|+|I_2|\le N-1\\|J|=|J'|=1}}\|\tau^{-1}(2+r-\tau)^{\f{1+\lz}{2}}\big\{|\Gamma^J\Gamma^{I_1}v|\cdot|\Gamma^{J'}\Gamma^{I_2}\psi|+|\tau-r|\cdot|\partial\Gamma^{I_1}v|\cdot|\partial\Gamma^{I_2}\psi|\big\}\|_{L^2_x(\Sigma^{\rm{ex}}_\tau)}\nonumber\\
&\lesssim&\sum_{\substack{|I_1|+|I_2|\le N-1\\|J|=|J'|=1}}\|\tau^{-1}(2+r-\tau)^{1+\f{1+\lz}{2}}|\Gamma^J\Gamma^{I_1}v|\cdot|\Gamma^{J'}\Gamma^{I_2}\psi|\,\|_{L^2_x(\Sigma^{\rm{ex}}_\tau)}\nonumber\\
&\lesssim&\sum_{\substack{|I_1|\le N-4\\|J|=1,|I_2|\le N}}\|\tau^{-1}(2+r-\tau)\Gamma^J\Gamma^{I_1}v\|_{L^\infty_x(\Sigma^{\rm{ex}}_\tau)}\cdot\|(2+r-\tau)^{\f{1+\lz}{2}}\Gamma^{I_2}\psi\|_{L^2_x(\Sigma^{\rm{ex}}_\tau)}\nonumber\\
&+&\!\!\!\sum_{\substack{|I_2|\le N-4\\|J'|=1,|I_1|\le N}}\|(2+r-\tau)^{\f{1+\lz}{2}}\Gamma^{I_1}v\|_{L^2_x(\Sigma^{\rm{ex}}_\tau)}\cdot\|\tau^{-1}(2+r-\tau)\Gamma^{J'}\Gamma^{I_2}\psi\|_{L^\infty_x(\Sigma^{\rm{ex}}_\tau)}\nonumber\\
&\lesssim&(C_1\ez)^2\tau^{-1-\lz_0+\f{\dz}{2}},
\eeqn
where we use \eqref{s2: Edecay0}, \eqref{s3: L^2psiv} and \eqref{s3: pointpsiv0*}. Hence \eqref{s3: tildevG1}-\eqref{s3: tildevG3} yield
\beq\label{s3: stildeG}
\|(2+r-\tau)^{\f{1+\lz}{2}}\Gamma^I\tilde{F}_\psi\|_{L^2_x(\Sigma^{\rm{ex}}_\tau)}\lesssim(C_1\ez)^2\big\{\sqrt{l(\tau)}\tau^{-1+\f{\dz}{2}}+\tau^{-1-\lz_0+\f{\dz}{2}}\big\},\quad |I|\le N-1,\tau\in[2,t].
\eeq
Recall that $0<\dz\ll\lz_0$. Then the estimates \eqref{s3: stildeG} and \eqref{s3: l(tau)L1} imply that for all $t\in[2,\infty)$
\beq\label{s3: gztiPsi}
[E^{{\rm{ex}},\lz}_0(t,\Gamma^I\tilde{\Psi})]^{\f{1}{2}}\lesssim\ez+(C_1\ez)^2,\quad |I|\le N-1.
\eeq
We recall that $\tilde{\Psi}=\Psi-v\psi$. For $|I|\le N-1$, by \eqref{s3: L^2psiv}, \eqref{s3: pointpsiv} and the definition \eqref{dexfsaz}, we have
\ben
&&[E^{{\rm{ex}},\lz}_0(t,\Gamma^I(v\psi))]^{\f{1}{2}}\lesssim\sum_{|I'|\le N}\|(2+r-t)^{\f{1+\lz}{2}}\Gamma^{I'}(v\psi)\|_{L^2_x(\Sigma^{\rm{ex}}_t)}\\
&\lesssim&\sum_{\substack{|I_1|\le N-3\\|I_2|\le N}}\|\Gamma^{I_1}v\|_{L^\infty_x(\Sigma^{\rm{ex}}_t)}\cdot\|(2+r-t)^{\f{1+\lz}{2}}\Gamma^{I_2}\psi\|_{L^2_x(\Sigma^{\rm{ex}}_t)}\\
&+&\sum_{\substack{|I_2|\le N-3\\|I_1|\le N}}\|t^{-\f{\dz}{2}}(2+r-t)^{\f{1+\lz}{2}}\Gamma^{I_1}v\|_{L^2_x(\Sigma^{\rm{ex}}_t)}\cdot\|t^{\f{\dz}{2}}\Gamma^{I_2}\psi\|_{L^\infty_x(\Sigma^{\rm{ex}}_t)}
\lesssim (C_1\ez)^2.
\een
This together with \eqref{s3: gztiPsi} gives the following estimate for all $t\in[2,\infty)$:
\beq\label{s3: gzv0}
[E^{{\rm{ex}},\lz}_{0}(t,\Gamma^I\Psi)]^{\f{1}{2}}\lesssim\ez+(C_1\ez)^2,\quad |I|\le N-1.
\eeq
By \eqref{s3: deflz_0}, \eqref{s3: 2+r-t} and Lemma \ref{s2: exL_auL^infty}, we have
\beqn\label{s3: PsiL^inftyim}
&&\sum_{\substack{|I|\le N-4\\|J|=1}}\sup_{\Sigma^{\rm{ex}}_t}r^{\lz_0}\big(|L_0\Gamma^I\Psi|+|Z^J\Gamma^I\Psi|\big)\lesssim\sum_{\substack{|I|\le N-4\\|J|=1}}\sup_{\Sigma^{\rm{ex}}_t}(2+r-t)^{\f{\lz-1}{2}} r^{\f{1}{2}}\big(|L_0\Gamma^I\Psi|+|Z^J\Gamma^I\Psi|\big)\nonumber\\
&\lesssim&\sum_{|I'|\le N-1}\|(2+r-t)^{\f{1+\lz}{2}}\partial\Gamma^{I'}\Psi\|_{L^2_x(\Sigma^{\rm{ex}}_t)}\lesssim\ez+(C_1\ez)^2
\eeqn
for all $t\in[2,\infty)$, where we recall \eqref{s2: Z^J}. By $iii)$ in Lemma \ref{s2: nlt}, we have
\beq\label{s3: hapsi=Psi}
\hat{\Gamma}^I\psi=i\gz^\mu\partial_\mu\hat{\Gamma}^I\Psi.
\eeq
Then by \eqref{s2: 2.2-44} in Lemma \ref{sC: 2.2-4}, Lemma \ref{sC: DMY2.1} and \eqref{s3: PsiL^inftyim}, we have the following estimate in the region $\mathscr{D}^{\rm{ex}}$:
\beq\label{s3: psi-L^inftyim}
|[\hat{\Gamma}^I\psi]_{-}|\lesssim |G\hat{\Gamma}^I\Psi|\lesssim r^{-1}\sum_{|J|=1}\big(|L_0\hat{\Gamma}^I\Psi|+|Z^J\hat{\Gamma}^I\Psi|\big)\lesssim[\ez+(C_1\ez)^2]r^{-1-\lz_0},\quad |I|\le N-4.
\eeq
We recall \eqref{s2: psi-} and \eqref{s2: psi)-}. By \eqref{s3: hapsi=Psi}, Lemma \ref{sC: DMY2.1}, \eqref{s3: psi-L^inftyim} and \eqref{s3: PsiL^inftyim}, in the region $\mathscr{D}^{\rm{ex}}$, it holds that 
\beqn\label{s3: (psi)-L^inftyim}
|(\hat{\Gamma}^I\psi)_{-}|&=&\l|[\hat{\Gamma}^I\psi]_{-}+\f{x_a}{r}\f{t-r}{t}\gz^0\gz^a\hat{\Gamma}^I\psi\r|\lesssim|[\hat{\Gamma}^I\psi]_{-}|+\f{|t-r|}{t}|\partial\hat{\Gamma}^I\Psi|\nonumber\\
&\lesssim&|[\hat{\Gamma}^I\psi]_{-}|+t^{-1}\sum_{|J|=1}\big(|L_0\hat{\Gamma}^I\Psi|+|Z^J\hat{\Gamma}^I\Psi|\big)\nonumber\\
&\lesssim&[\ez+(C_1\ez)^2]r^{-\lz_0}t^{-1},\quad |I|\le N-4.
\eeqn

Finally, to close the bootstrap in the interior region in Section \ref{s4}, we give the following estimate of the function $\tilde{\psi}:=\psi-i\gz^\mu\partial_\mu(v\psi)$. For $|I|\le N-1$, by $i)$ in Lemma \ref{s2: nlt} and Proposition \ref{s2: exfsD}, we have
\beq\label{s3: S1tipsi}
\l(\int_{\mathscr{l}_{[2,t]}}|[\hat{\Gamma}^I\tilde{\psi}]_{-}|^2{\rm{d}}\sz\r)^{\f{1}{2}}\lesssim\|\hat{\Gamma}^I\tilde{\psi}\|_{Y^{{\rm{ex}},\lz}_{D,2}}+\int_2^t\|(2+r-\tau)^{\f{1+\lz}{2}}\hat{\Gamma}^I\tilde{F}_\psi\|_{L^2_x(\Sigma^{\rm{ex}}_\tau)}{\rm{d}}\tau\lesssim\ez+(C_1\ez)^2
\eeq
for all $t\in[2,\infty)$, where $\tilde{F}_\psi$ is as in \eqref{s3: tildeG} and we use \eqref{s3: stildeG}.

\subsection{Estimates of the solution $(\psi,v)$ on truncated exterior hyperboloids}\label{s3.3}

${\bf Step\ 1.}$ First, we derive the highest order energy estimates of $(\psi,v)$ on exterior hyperboloids. By Propositions \ref{s2: exhsD} and \ref{s2: exhs}, for $|I|\le N$ and any $s\in[2,\infty)$, we have
\beqn\label{s3: psiexh}
\mathcal{E}^{{\rm{ex}},h}_D(s,\hat{\Gamma}^I\psi)+\int_{\mathscr{l}_{[2,\,t(s)]}}|[\hat{\Gamma}^I\psi]_{-}|^2{\rm{d}}\sz&\lesssim& E^{\rm{ex}}_D(2,\hat{\Gamma}^I\psi)+\int_2^\infty\|(\hat{\Gamma}^I\psi)^*\gz^0\hat{\Gamma}^I(v\psi)\|_{L^1_x(\Sigma^{\rm{ex}}_t)}{\rm{d}}t\nonumber\\
&\lesssim&\ez^2+(C_1\ez)^3
\eeqn
and
\beqn\label{s3: vexh}
\!\!\!\!\!\!\!\!\!\!\!\!\!\!&&\mathcal{E}^{{\rm{ex}},h}_{1,\dz}(s,\Gamma^Iv)+\int_{\mathscr{l}_{[2,\,t(s)]}}t^{-\dz}\big(|G\Gamma^Iv|^2+|\Gamma^Iv|^2\big){\rm{d}}\sz\nonumber\\
\!\!\!\!\!\!\!\!\!\!\!\!\!\!&\lesssim& E^{\rm{ex}}_{1,\dz}(2,\Gamma^Iv)+\sup_{t\in[2,\infty)}[E^{\rm{ex}}_{1,\dz}(t,\Gamma^Iv)]^{\f{1}{2}}\cdot\int_2^\infty\|t^{-\f{\dz}{2}}\Gamma^I(\psi^*\gz^0\psi)\|_{L^2_x(\Sigma^{\rm{ex}}_t)}{\rm{d}}t\lesssim\ez^2+(C_1\ez)^3,
\eeqn
where $t(s)=\f{s^2+1}{2}$, we use \eqref{s3: psiNL^2}, \eqref{s3: vNL^2} and \eqref{s3: S2v} (recall \eqref{s3: Ygzvdz} and \eqref{dexfs}), and recall \eqref{dexhsD}, \eqref{dexfsD0}, \eqref{dexhs} and \eqref{dexfsdz} for the definitions of $\mathcal{E}^{{\rm{ex}},h}_D(s,\hat{\Gamma}^I\psi)$, $ E^{\rm{ex}}_D(2,\hat{\Gamma}^I\psi)$, $\mathcal{E}^{{\rm{ex}},h}_{1,\dz}(s,\Gamma^Iv)$ and $E^{\rm{ex}}_{1,\dz}(t,\Gamma^Iv)$.

${\bf Step\ 2.}$ Next we show lower order energy estimates of $v$ on exterior hyperboloids. Let $\tilde{v}=v-\psi^*\gz^0\psi$ be as in ${\bf Step\ 3}$ in Section \ref{s3.2}. By Proposition \ref{s2: exhs}, for $|I|\le N-1$ and all $s\in[2,\infty)$, we have
\beqn\label{s3: tivexh}
&&\mathcal{E}^{{\rm{ex}},h}_1(s,\Gamma^I\tilde{v})+\int_{\mathscr{l}_{[2,\,t(s)]}}\big(|G\Gamma^I\tilde{v}|^2+|\Gamma^I\tilde{v}|^2\big){\rm{d}}\sz\nonumber\\
&\lesssim& E^{\rm{ex}}_1(2,\Gamma^I\tilde{v})+\sup_{t\in[2,\infty)}[E^{\rm{ex}}_1(t,\Gamma^I\tilde{v})]^{\f{1}{2}}\cdot\int_2^\infty\|\Gamma^I\tilde{F}_v\|_{L^2_x(\Sigma^{\rm{ex}}_t)}{\rm{d}}t\lesssim\ez^2+(C_1\ez)^3,
\eeqn
where $t(s)=\f{s^2+1}{2}$, $\tilde{F}_v$ is as in \eqref{s3: tildeF}, and we use \eqref{s3: stildevF} and \eqref{s3: S3tildev} (recall \eqref{s3: exnvN-3} and \eqref{dexfsaz}). We also recall \eqref{dhs0} and \eqref{dexfs0} for the definitions of $\mathcal{E}^{{\rm{ex}},h}_1(s,\Gamma^I\tilde{v})$ and $E^{\rm{ex}}_1(t,\Gamma^I\tilde{v})$. For $|I|\le N-1$, by \eqref{dhs0}, and \eqref{s2: 2.2-42} and \eqref{s2: 2.2-43} in Lemma \ref{sC: 2.2-4}, we have
\ben
&&[\mathcal{E}^{{\rm{ex}},h}_1(s,\Gamma^I(\psi^*\gz^0\psi))]^{\f{1}{2}}\lesssim\sum_{|J|\le 1}\|\Gamma^J\Gamma^I(\psi^*\gz^0\psi)\|_{L^2(\mathscr{H}^{\rm{ex}}_s)}\\
&\lesssim&\sum_{|I_1|+|I_2|\le N}\|\,|(\hat{\Gamma}^{I_1}\psi)_{-}|\cdot|\hat{\Gamma}^{I_2}\psi|+(s^2/t^2)|\hat{\Gamma}^{I_1}\psi|\cdot\hat{\Gamma}^{I_2}\psi|\,\|_{L^2(\mathscr{H}^{\rm{ex}}_s)}\\
&\lesssim&\sum_{\substack{|I_1|\le N-4\\|I_2|\le N}}\|t(\hat{\Gamma}^{I_1}\psi)_{-}\|_{L^\infty(\mathscr{H}^{\rm{ex}}_s)}\cdot\|t^{-1}\hat{\Gamma}^{I_2}\psi\|_{L^2(\mathscr{H}^{\rm{ex}}_s)}+\sum_{\substack{|I_2|\le N-3\\|I_1|\le N}}\|(\hat{\Gamma}^{I_1}\psi)_{-}\|_{L^2(\mathscr{H}^{\rm{ex}}_s)}\cdot\|\hat{\Gamma}^{I_2}\psi\|_{L^\infty(\mathscr{H}^{\rm{ex}}_s)}\\
&+&\sum_{\substack{|I_1|\le N-3\\|I_2|\le N}}\|\hat{\Gamma}^{I_1}\psi\|_{L^\infty(\mathscr{H}^{\rm{ex}}_s)}\cdot\|(s/t)\hat{\Gamma}^{I_2}\psi\|_{L^2(\mathscr{H}^{\rm{ex}}_s)}\lesssim (C_1\ez)^2,
\een
where we use \eqref{s3: (psi)-L^inftyim}, \eqref{s3: psiexh} and \eqref{s3: pointpsiv}, and recall the definition \eqref{dexhsD}. This together with \eqref{s3: tivexh} gives
\beqn\label{s3: vexhN-3}
\mathcal{E}^{{\rm{ex}},h}_1(s,\Gamma^Iv)\lesssim\ez^2+(C_1\ez)^3,\quad |I|\le N-1.
\eeqn

\section{$n=2$: Global existence and estimates in the interior region}\label{s4}

In this section, we prove global existence of the solution to \eqref{s1: DKG}-\eqref{s1: initial} in the interior region.

\subsection{Bootstrap assumptions}

Let $N\ge 6$, $0<\dz\ll \lz_0$ and $C_1\gg1$ be as in Section \ref{s3.1}. Let $C_2\gg C_1$ and $0<\ez\ll C_2^{-1}$ be chosen later. The bootstrap assumptions in the interior region are the following energy bounds for the solution $(\psi,v)$ to \eqref{s1: DKG}-\eqref{s1: initial}:
\beq\label{s4: bspsiv}
\sum_{|I|\le N} \big\{[\mathcal{E}^{{\rm{in}},h}_D(s,\hat{\Gamma}^I\psi)]^{\f{1}{2}}+[\mathcal{E}^{{\rm{in}},h}_1(s,\Gamma^Iv)]^{\f{1}{2}}\big\}\le C_2\ez s^\dz,
\eeq
\beq\label{s4: bspsivN-3}
\sum_{|I|\le N-1} \big\{[\mathcal{E}^{{\rm{in}},h}_D(s,\hat{\Gamma}^I\psi)]^{\f{1}{2}}+[\mathcal{E}^{{\rm{in}},h}_1(s,\Gamma^Iv)]^{\f{1}{2}}\big\}\le C_2\ez,
\eeq
where we recall \eqref{dinhsD} and \eqref{dhs0} for the definitions of the energy functionals above. Let
\beq\label{s4: s*}
\tilde{T}^*:=\sup\{\tilde{T}>2: \eqref{s4: bspsiv}-\eqref{s4: bspsivN-3}\ \mathrm{hold\ for\ }s\in[2,\tilde{T}]\}.
\eeq

\begin{prop}\label{s4: maxtildeT}
There exist some constants $C_2>0$ sufficiently large and $0<\ez_1\ll C_2^{-1}$ sufficiently small such that for any $0<\ez<\ez_1$, if $(\psi,v)$ is a solution to \eqref{s1: DKG}-\eqref{s1: initial} and satisfies the exterior estimate \eqref{s3: bspsiv} globally in time as well as the interior bounds \eqref{s4: bspsiv}-\eqref{s4: bspsivN-3} for all $s\in[2,\tilde{T}]$, then we have the following improved interior estimates for all $s\in[2,\tilde{T}]$:
\be
\begin{split}
&\sum_{|I|\le N} \big\{[\mathcal{E}^{{\rm{in}},h}_D(s,\hat{\Gamma}^I\psi)]^{\f{1}{2}}+[\mathcal{E}^{{\rm{in}},h}_1(s,\Gamma^Iv)]^{\f{1}{2}}\big\}\le \f{1}{2}C_2\ez s^\dz,\\
&\sum_{|I|\le N-1} \big\{[\mathcal{E}^{{\rm{in}},h}_D(s,\hat{\Gamma}^I\psi)]^{\f{1}{2}}+[\mathcal{E}^{{\rm{in}},h}_1(s,\Gamma^Iv)]^{\f{1}{2}}\big\}\le \f{1}{2}C_2\ez.
\end{split}
\ee
\end{prop}

In the above proposition the hyperbolic time $\tilde{T}$ is arbitrary, hence the solution $(\psi,v)$ exists globally in the interior region $\mathscr{D}^{\rm{in}}$ (i.e., $\tilde{T}^*=\infty$ where $\tilde{T}^*$ is as in \eqref{s4: s*}) and satisfies the estimates \eqref{s4: bspsiv}-\eqref{s4: bspsivN-3} for all $s\in[2,\infty)$. Below we give the proof of Proposition \ref{s4: maxtildeT}. In the sequel, the implied constants in $\lesssim$ do not depend on the constants $C_2$ and $\ez$ appearing in the bootstrap assumptions \eqref{s4: bspsiv}-\eqref{s4: bspsivN-3}.

We recall the definitions \eqref{dinhsD}, \eqref{dexhsD}, \eqref{dhs0}, \eqref{s2: L^pH_s} and that $\mathscr{H}_s=\mathscr{H}^{\rm{in}}_s\cup\mathscr{H}^{\rm{ex}}_s$ by \eqref{s2: mscH}-\eqref{s2: mscHex}. By \eqref{s3: psiexh}, \eqref{s3: vexhN-3},  \eqref{s4: bspsiv} and \eqref{s4: bspsivN-3}, we obtain the following $L^2$-type estimates for the solution $(\psi,v)$, for $s\in[2,\tilde{T}]$:
\beq\label{s4: psivL^2}
\sum_{|I|\le N}\|\,|(\hat{\Gamma}^I\psi)_{-}|+|(s/t)\hat{\Gamma}^I\psi|+|(s/t)\partial\Gamma^Iv|+|\Gamma^Iv|\,\|_{L^2(\mathscr{H}^{\rm{in}}_s)}\lesssim C_2\ez s^\dz,
\eeq
\beq\label{s4: psivL^2N-3}
\sum_{|I|\le N-1}\|\,|(\hat{\Gamma}^I\psi)_{-}|+|(s/t)\hat{\Gamma}^I\psi|+|(s/t)\partial\Gamma^Iv|+|\Gamma^Iv|\,\|_{L^2(\mathscr{H}_s)}\lesssim C_2\ez.
\eeq
By \eqref{s4: psivL^2N-3}, Lemmas \ref{s4: Sobin} and \ref{sC: DL2.3'}, the following pointwise estimates hold for $s\in[2,\tilde{T}]$:
\beq\label{s4: psivL^infty}
\sum_{|I|\le N-3}\|t\big\{|(\hat{\Gamma}^I\psi)_{-}|+|(s/t)\hat{\Gamma}^I\psi|+|(s/t)\partial\Gamma^Iv|+|\Gamma^Iv|\big\}\|_{L^\infty(\mathscr{H}^{\rm{in}}_s)}\lesssim C_2\ez.
\eeq

\subsection{Improved estimates for the solution $(\psi,v)$ in the interior region}\label{ss4.2}

${\bf Step\ 1.}$ First, we refine the highest order energy estimates of $\psi$ and $v$. By Proposition \ref{s2: inhsD}, for $|I|\le N$, we have
\beqn\label{s4: S1gzpsi}
\mathcal{E}^{{\rm{in}},h}_D(s,\hat{\Gamma}^I\psi)&\lesssim&\mathcal{E}^{{\rm{in}},h}_D(2,\hat{\Gamma}^I\psi)+\int_{\mathscr{l}_{[t(2),\,t(s)]}}|[\hat{\Gamma}^I\psi]_{-}|^2{\rm{d}}\sz\nonumber\\
&+&\int_2^s[\mathcal{E}^{{\rm{in}},h}_D(\tau,\hat{\Gamma}^I\psi)]^{\f{1}{2}}\cdot\|\hat{\Gamma}^I(v\psi)\|_{L^2(\mathscr{H}^{\rm{in}}_\tau)}{\rm{d}}\tau,
\eeqn
where 
\beq\label{t(s)t(2)}
t(s)=\f{s^2+1}{2},\quad\quad t(2)=\f{5}{2}.
\eeq
For $|I|\le N$ and $\tau\in[2,s]$, by \eqref{s2: 2.2-41} in Lemma \ref{sC: 2.2-4}, \eqref{s4: psivL^2} and \eqref{s4: psivL^infty}, we have
\beqn\label{s4: higpsi}
&&\|\hat{\Gamma}^I(v\psi)\|_{L^2(\mathscr{H}^{\rm{in}}_\tau)}\lesssim\sum_{|I_1|+|I_2|\le N}\|\,|\Gamma^{I_1}v|\cdot|\hat{\Gamma}^{I_2}\psi|\,\|_{L^2(\mathscr{H}^{\rm{in}}_\tau)}\nonumber\\
&\lesssim&\sum_{\substack{|I_1|\le N-3\\|I_2|\le N}}\|(t/\tau)\Gamma^{I_1}v\|_{L^\infty(\mathscr{H}^{\rm{in}}_\tau)}\cdot\|(\tau/t)\hat{\Gamma}^{I_2}\psi\|_{L^2(\mathscr{H}^{\rm{in}}_\tau)}+\sum_{\substack{|I_2|\le N-3\\|I_1|\le N}}\|\Gamma^{I_1}v\|_{L^2(\mathscr{H}^{\rm{in}}_\tau)}\cdot\|\hat{\Gamma}^{I_2}\psi\|_{L^\infty(\mathscr{H}^{\rm{in}}_\tau)}\nonumber\\
&\lesssim& (C_2\ez)^2\tau^{-1+\dz}.
\eeqn
Note that the estimate of the second term on the right hand side of \eqref{s4: S1gzpsi} was given in \eqref{s3: psiexh}. It follows that
\be
\mathcal{E}^{{\rm{in}},h}_D(s,\hat{\Gamma}^I\psi)\lesssim \ez^2+(C_1\ez)^3+(C_2\ez)^2 s^\dz\sup_{\tau\in[2,s]}[\mathcal{E}^{{\rm{in}},h}_D(\tau,\hat{\Gamma}^I\psi)]^{\f{1}{2}},\quad |I|\le N,
\ee
which implies
\beq\label{s4: repsiN}
\sum_{|I|\le N}[\mathcal{E}^{{\rm{in}},h}_D(s,\hat{\Gamma}^I\psi)]^{\f{1}{2}}\lesssim\ez+(C_2\ez)^{\f{3}{2}}s^{\dz}.
\eeq
By Proposition \ref{s2: inhs}, for $|I|\le N$, we have
\beqn\label{s4: Einv}
\mathcal{E}^{{\rm{in}},h}_1(s,\Gamma^Iv)\lesssim \mathcal{E}^{{\rm{in}},h}_1(2,\Gamma^Iv)&+&\int_{\mathscr{l}_{[t(2),\,t(s)]}}\big(|G\Gamma^Iv|^2+|\Gamma^Iv|^2\big){\rm{d}}\sz\nonumber\\
&+&\int_2^s[\mathcal{E}^{{\rm{in}},h}_1(\tau,\Gamma^Iv)]^{\f{1}{2}}\cdot\|\Gamma^I(\psi^*\gz^0\psi)\|_{L^2(\mathscr{H}^{\rm{in}}_\tau)}{\rm{d}}\tau,
\eeqn
where $t(s), t(2)$ are as in \eqref{t(s)t(2)}. For $|I|\le N$ and $\tau\in[2,s]$, by \eqref{s2: 2.2-42} and \eqref{s2: 2.2-43} in Lemma \ref{sC: 2.2-4},  \eqref{s4: psivL^2} and \eqref{s4: psivL^infty}, we have
\beqn\label{s4: higv}
&&\|\Gamma^I(\psi^*\gz^0\psi)\|_{L^2(\mathscr{H}^{\rm{in}}_\tau)}\lesssim\sum_{|I_1|+|I_2|\le N}\|\,|(\hat{\Gamma}^{I_1}\psi)_{-}|\cdot|\hat{\Gamma}^{I_2}\psi|+(\tau^2/t^2)|\hat{\Gamma}^{I_1}\psi|\cdot|\hat{\Gamma}^{I_2}\psi|\,\|_{L^2(\mathscr{H}^{\rm{in}}_\tau)}\nonumber\\
&\lesssim&\sum_{\substack{|I_1|\le N-3\\|I_2|\le N}}\|(t/\tau)(\hat{\Gamma}^{I_1}\psi)_{-}\|_{L^\infty(\mathscr{H}^{\rm{in}}_\tau)}\cdot\|(\tau/t)\hat{\Gamma}^{I_2}\psi\|_{L^2(\mathscr{H}^{\rm{in}}_\tau)}+\sum_{\substack{|I_2|\le N-3\\|I_1|\le N}}\|(\hat{\Gamma}^{I_1}\psi)_{-}\|_{L^2(\mathscr{H}^{\rm{in}}_\tau)}\cdot\|\hat{\Gamma}^{I_2}\psi\|_{L^\infty(\mathscr{H}^{\rm{in}}_\tau)}\nonumber\\
&+&\sum_{\substack{|I_1|\le N-3\\|I_2|\le N}}\|\hat{\Gamma}^{I_1}\psi\|_{L^\infty(\mathscr{H}^{\rm{in}}_\tau)}\cdot\|(\tau/t)\hat{\Gamma}^{I_2}\psi\|_{L^2(\mathscr{H}^{\rm{in}}_\tau)}\lesssim (C_2\ez)^2\tau^{-1+\dz}.
\eeqn
We recall \eqref{s2: l[2,t]} and note that on $\mathscr{l}_{[2,\,t(s)]}$ we have $s^{-2\dz}\lesssim t^{-\dz}$. Hence by \eqref{s4: Einv}, \eqref{s3: vexh} and \eqref{s4: higv}, we have
\be
\mathcal{E}^{{\rm{in}},h}_1(s,\Gamma^Iv)\lesssim[\ez^2+(C_1\ez)^3]s^{2\dz}+(C_2\ez)^2 s^\dz\sup_{\tau\in[2,s]}[\mathcal{E}^{{\rm{in}},h}_1(\tau,\Gamma^Iv)]^{\f{1}{2}},\quad |I|\le N,
\ee
which implies 
\beq\label{s4: revN}
\sum_{|I|\le N}[\mathcal{E}^{{\rm{in}},h}_1(s,\Gamma^Iv)]^{\f{1}{2}}\lesssim[\ez+(C_2\ez)^{\f{3}{2}}]s^{\dz}.
\eeq

${\bf Step\ 2.}$ Next, we refine the lower order energy estimates of $\psi$. Let $\tilde{\psi}:=\psi-i\gz^\mu\partial_\mu(v\psi)$. For $|I|\le N-1$, by $i)$ in Lemma \ref{s2: nlt} and Proposition \ref{s2: inhsD}, we have
\beqn\label{s4: S2tipsiN-}
\mathcal{E}^{{\rm{in}},h}_D(s,\hat{\Gamma}^I\tilde{\psi})&\lesssim& \mathcal{E}^{{\rm{in}},h}_D(2,\hat{\Gamma}^I\tilde{\psi})+\int_{\mathscr{l}_{[t(2),\,t(s)]}}|[\hat{\Gamma}^I\tilde{\psi}]_{-}|^2{\rm{d}}\sz\nonumber\\
&+&\int_2^s[\mathcal{E}^{{\rm{in}},h}_D(\tau,\hat{\Gamma}^I\tilde{\psi})]^{\f{1}{2}}\cdot\|\hat{\Gamma}^I\tilde{F}_\psi\|_{L^2(\mathscr{H}^{\rm{in}}_\tau)}{\rm{d}}\tau,
\eeqn
where $t(s), t(2)$ are as in \eqref{t(s)t(2)} and $\tilde{F}_\psi$ is as in \eqref{s3: tildeG}, i.e.,
\be\label{s4: tildeFps}
\tilde{F}_\psi=(\psi^*\gz^0\psi)\psi-iv\gz^\mu\partial_\mu(v\psi)-2\partial_\az v\partial^\az\psi.
\ee
For $|I|\le N-1$ and $\tau\in[2,s]$, by \eqref{s4: psivL^2}, \eqref{s4: psivL^2N-3}, \eqref{s4: psivL^infty} and \eqref{s2: 2.2-41}, \eqref{s2: 2.2-42} and \eqref{s2: 2.2-43} in Lemma \ref{sC: 2.2-4}, we have
\beqn\label{s4: gzpsg0ps}
&&\|\hat{\Gamma}^I\big((\psi^*\gz^0\psi)\psi\big)\|_{L^2(\mathscr{H}^{\rm{in}}_\tau)}\lesssim\sum_{|I_1|+|I_2|+|I_3|\le N-1}\|\,|(\hat{\Gamma}^{I_1}\psi)^*\gz^0\hat{\Gamma}^{I_2}\psi|\cdot|\hat{\Gamma}^{I_3}\psi|\,\|_{L^2(\mathscr{H}^{\rm{in}}_\tau)}\nonumber\\
&\lesssim&\sum_{|I_1|+|I_2|+|I_3|\le N-1}\|\,|(\hat{\Gamma}^{I_1}\psi)_{-}|\cdot|\hat{\Gamma}^{I_2}\psi|\cdot|\hat{\Gamma}^{I_3}\psi|+(\tau^2/t^2)|\hat{\Gamma}^{I_1}\psi|\cdot|\hat{\Gamma}^{I_2}\psi|\cdot|\hat{\Gamma}^{I_3}\psi|\,\|_{L^2(\mathscr{H}^{\rm{in}}_\tau)}\nonumber\\
&\lesssim&\sum_{\substack{|I_1|,|I_2|\le N-3\\|I_3|\le N-1}}\|(t/\tau)(\hat{\Gamma}^{I_1}\psi)_{-}\|_{L^\infty(\mathscr{H}^{\rm{in}}_\tau)}\cdot\|\hat{\Gamma}^{I_2}\psi\|_{L^\infty(\mathscr{H}^{\rm{in}}_\tau)}\cdot\|(\tau/t)\hat{\Gamma}^{I_3}\psi\|_{L^2(\mathscr{H}^{\rm{in}}_\tau)}\nonumber\\
&+&\sum_{\substack{|I_2|,|I_3|\le N-3\\|I_1|\le N-1}}\|(\hat{\Gamma}^{I_1}\psi)_{-}\|_{L^2(\mathscr{H}^{\rm{in}}_\tau)}\cdot\|\hat{\Gamma}^{I_2}\psi\|_{L^\infty(\mathscr{H}^{\rm{in}}_\tau)}\cdot\|\hat{\Gamma}^{I_3}\psi\|_{L^\infty(\mathscr{H}^{\rm{in}}_\tau)}\nonumber\\
&+&\sum_{\substack{|I_1|,|I_2|\le N-3\\|I_3|\le N-1}}\|\hat{\Gamma}^{I_1}\psi\|_{L^\infty(\mathscr{H}^{\rm{in}}_\tau)}\cdot\|\hat{\Gamma}^{I_2}\psi\|_{L^\infty(\mathscr{H}^{\rm{in}}_\tau)}\cdot\|(\tau/t)\hat{\Gamma}^{I_3}\psi\|_{L^2(\mathscr{H}^{\rm{in}}_\tau)}\lesssim (C_2\ez)^3\tau^{-2}
\eeqn
and
\beqn\label{s4: gvgvgpsi}
&&\|\hat{\Gamma}^I\big(v\gz^\mu\partial_\mu(v\psi)\big)\|_{L^2(\mathscr{H}^{\rm{in}}_\tau)}\lesssim\sum_{|I_1|+|I_2|+|I_3|\le N}\|\,|\Gamma^{I_1}v|\cdot|\Gamma^{I_2}v|\cdot|\Gamma^{I_3}\psi|\,\|_{L^2(\mathscr{H}^{\rm{in}}_\tau)}\nonumber\\
&\lesssim&\sum_{\substack{|I_1|,|I_2|\le N-3\\|I_3|\le N}}\|(t/\tau)\Gamma^{I_1}v\|_{L^\infty(\mathscr{H}^{\rm{in}}_\tau)}\cdot\|\Gamma^{I_2}v\|_{L^\infty(\mathscr{H}^{\rm{in}}_\tau)}\cdot\|(\tau/t)\Gamma^{I_3}\psi\|_{L^2(\mathscr{H}^{\rm{in}}_\tau)}\nonumber\\
&+&\sum_{\substack{|I_2|,|I_3|\le N-3\\|I_1|\le N}}\|\Gamma^{I_1}v\|_{L^2(\mathscr{H}^{\rm{in}}_\tau)}\cdot\|\Gamma^{I_2}v\|_{L^\infty(\mathscr{H}^{\rm{in}}_\tau)}\cdot\|\Gamma^{I_3}\psi\|_{L^\infty(\mathscr{H}^{\rm{in}}_\tau)}\lesssim (C_2\ez)^3\tau^{-2+\dz}.
\eeqn
For $|I|\le N-1$ and $\tau\in[2,s]$, by Lemma \ref{s2: Q_0} and \eqref{s2: parpsi} in Lemma \ref{s2: Edecay}, we have
\beqn\label{s4: gzpavps}
\!\!\!&&\|\Gamma^I\big(\partial_\az v\partial^\az\psi\big)\|_{L^2(\mathscr{H}^{\rm{in}}_\tau)}\lesssim\sum_{|I_1|+|I_2|\le N-1}\|(\partial_\az \Gamma^{I_1}v)(\partial^\az\Gamma^{I_2}\psi)\|_{L^2(\mathscr{H}^{\rm{in}}_\tau)}\nonumber\\
\!\!\!&\lesssim&\sum_{\substack{|I_1|+|I_2|\le N-1\\|J|=|J'|=1}}\|t^{-1}|\Gamma^J\Gamma^{I_1}v|\cdot|\Gamma^{J'}\Gamma^{I_2}\psi|+((t-r)/t)|\partial\Gamma^{I_1}v|\cdot|\partial\Gamma^{I_2}\psi|\,\|_{L^2(\mathscr{H}^{\rm{in}}_\tau)}\nonumber\\
\!\!\!&\lesssim&\sum_{\substack{|I_1|+|I_2|\le N-1\\|J|=|J'|=1}}\|t^{-1}|\Gamma^J\Gamma^{I_1}v|\cdot|\Gamma^{J'}\hat{\Gamma}^{I_2}\psi|+|\Gamma^J\Gamma^{I_1}v|\cdot|\hat{\Gamma}^{I_2}(v\psi)|\,\|_{L^2(\mathscr{H}^{\rm{in}}_\tau)}\nonumber\\
\!\!\!&\lesssim&\!\!\!\sum_{\substack{|I_1|\le N-4\\|J|=1,|I_2|\le N}}\|\tau^{-1}\Gamma^J\Gamma^{I_1}v\|_{L^\infty(\mathscr{H}^{\rm{in}}_\tau)}\cdot\l\|\f{\tau}{t}\hat{\Gamma}^{I_2}\psi\r\|_{L^2(\mathscr{H}^{\rm{in}}_\tau)}+\sum_{\substack{|I_2|\le N-4\\|J'|=1,|I_1|\le N}}\|\Gamma^{I_1}v\|_{L^2(\mathscr{H}^{\rm{in}}_\tau)}\cdot\|t^{-1}\Gamma^{J'}\hat{\Gamma}^{I_2}\psi\|_{L^\infty(\mathscr{H}^{\rm{in}}_\tau)}\nonumber\\
\!\!\!&+&\sum_{|I_1|+|I_2|+|I_3|\le N}\|\,|\Gamma^{I_1}v|\cdot|\Gamma^{I_2}v|\cdot|\hat{\Gamma}^{I_3}\psi|\,\|_{L^2(\mathscr{H}^{\rm{in}}_\tau)}\lesssim(C_2\ez)^2\tau^{-2+\dz},
\eeqn
where we use \eqref{s4: psivL^2}, \eqref{s4: psivL^infty} and \eqref{s4: gvgvgpsi}. Combining \eqref{s4: gzpsg0ps}-\eqref{s4: gzpavps}, we derive
\beq\label{s4: tildeGL^2}
\|\hat{\Gamma}^I\tilde{F}_\psi\|_{L^2(\mathscr{H}^{\rm{in}}_\tau)}\lesssim(C_2\ez)^2\tau^{-2+\dz},\quad |I|\le N-1, \tau\in[2,s].
\eeq
We point out that the estimate of the second term on the right hand side of \eqref{s4: S2tipsiN-} was given by \eqref{s3: S1tipsi}. It follows that
\ben
\mathcal{E}^{{\rm{in}},h}_D(s,\hat{\Gamma}^I\tilde{\psi})\lesssim\ez^2+(C_1\ez)^4+(C_2\ez)^2\sup_{\tau\in[2,s]}[\mathcal{E}^{{\rm{in}},h}_D(\tau,\hat{\Gamma}^I\tilde{\psi})]^{\f{1}{2}},\quad |I|\le N-1,
\een
which implies
\beq\label{s4: S2imtipsi}
\sum_{|I|\le N-1}[\mathcal{E}^{{\rm{in}},h}_D(s,\hat{\Gamma}^I\tilde{\psi})]^{\f{1}{2}}\lesssim\ez+(C_2\ez)^{2}.
\eeq
We recall the definition \eqref{dinhsD} and that $\tilde{\psi}=\psi-i\gz^\mu\partial_\mu(v\psi)$. We note that
\ben
&&\sum_{|I|\le N-1}\|\hat{\Gamma}^I\big(i\gz^\mu\partial_\mu(v\psi)\big)\|_{L^2(\mathscr{H}^{\rm{in}}_s)}\lesssim\sum_{|I_1|+|I_2|\le N}\|\,|\Gamma^{I_1}v|\cdot|\hat{\Gamma}^{I_2}\psi|\,\|_{L^2(\mathscr{H}^{\rm{in}}_s)}\\
&\lesssim&\sum_{\substack{|I_1|\le N-3\\|I_2|\le N}}\l\|s^{\dz}\f{t}{s}\Gamma^{I_1}v\r\|_{L^\infty(\mathscr{H}^{\rm{in}}_s)}\cdot\l\|s^{-\dz}\f{s}{t}\hat{\Gamma}^{I_2}\psi\r\|_{L^2(\mathscr{H}^{\rm{in}}_s)}+\sum_{\substack{|I_2|\le N-3\\|I_1|\le N}}\|s^{-\dz}\Gamma^{I_1}v\|_{L^2(\mathscr{H}^{\rm{in}}_s)}\cdot\|s^{\dz}\hat{\Gamma}^{I_2}\psi\|_{L^\infty(\mathscr{H}^{\rm{in}}_s)}\\
&\lesssim& (C_2\ez)^2
\een
by \eqref{s4: psivL^2} and \eqref{s4: psivL^infty}. This together with \eqref{s4: S2imtipsi} implies
\beq\label{s4: S2psN-}
\sum_{|I|\le N-1}[\mathcal{E}^{{\rm{in}},h}_D(s,\hat{\Gamma}^I\psi)]^{\f{1}{2}}\lesssim\ez+(C_2\ez)^{2}.
\eeq

${\bf Step\ 3.}$ We turn to the lower order energy estimates of $v$. Let $\tilde{v}:=v-\psi^*\gz^0\psi$. For $|I|\le N-1$, by $ii)$ in Lemma \ref{s2: nlt} and Proposition \ref{s2: inhs}, we have
\beqn\label{s4: S3tivN-}
\mathcal{E}^{{\rm{in}},h}_1(s,\Gamma^I\tilde{v})&\lesssim& \mathcal{E}^{{\rm{in}},h}_1(2,\Gamma^{I}\tilde{v})+\int_{\mathscr{l}_{[t(2),\,t(s)]}}\big(|G\Gamma^{I}\tilde{v}|^2+|\Gamma^{I}\tilde{v}|^2\big){\rm{d}}\sz\nonumber\\
&+&\int_2^s[\mathcal{E}^{{\rm{in}},h}_1(\tau,\Gamma^{I}\tilde{v})]^{\f{1}{2}}\cdot\|\Gamma^{I}\tilde{F}_v\|_{L^2(\mathscr{H}^{\rm{in}}_\tau)}{\rm{d}}\tau,
\eeqn
where $t(s), t(2)$ are as in \eqref{t(s)t(2)} and $\tilde{F}_v$ is as in \eqref{s3: tildeF}, i.e.,
\be
\tilde{F}_v=-i\partial_\mu(v\psi^*)\gz^0\gz^\mu\psi+i\psi^*\gz^0\gz^\mu\partial_\mu(v\psi)+2\partial_\az\psi^*\gz^0\partial^\az\psi.
\ee
For $|I|\le N-1$ and $\tau\in[2,s]$, we have
\beqn\label{s4: gpvpsps}
&&\|\Gamma^{I}\big(\partial_\mu(v\psi^*)\gz^0\gz^\mu\psi\big)\|_{L^2(\mathscr{H}^{\rm{in}}_\tau)}\lesssim\sum_{|I_1|+|I_2|+|I_3|\le N}\|\,|\Gamma^{I_1}v|\cdot|\Gamma^{I_2}\psi|\cdot|\Gamma^{I_3}\psi|\,\|_{L^2(\mathscr{H}^{\rm{in}}_\tau)}\nonumber\\
&\lesssim&\sum_{\substack{|I_1|,|I_2|\le N-3\\|I_3|\le N}}\|(t/\tau)\Gamma^{I_1}v\|_{L^\infty(\mathscr{H}^{\rm{in}}_\tau)}\cdot\|\hat{\Gamma}^{I_2}\psi\|_{L^\infty(\mathscr{H}^{\rm{in}}_\tau)}\cdot\|(\tau/t)\hat{\Gamma}^{I_3}\psi\|_{L^2(\mathscr{H}^{\rm{in}}_\tau)}\nonumber\\
&+&\sum_{\substack{|I_2|,|I_3|\le N-3\\|I_1|\le N}}\|\Gamma^{I_1}v\|_{L^2(\mathscr{H}^{\rm{in}}_\tau)}\cdot\|\hat{\Gamma}^{I_2}\psi\|_{L^\infty(\mathscr{H}^{\rm{in}}_\tau)}\cdot\|\hat{\Gamma}^{I_3}\psi\|_{L^\infty(\mathscr{H}^{\rm{in}}_\tau)}\lesssim (C_2\ez)^3\tau^{-2+\dz}.
\eeqn
By Lemma \ref{s2: Q_0} and \eqref{s2: parpsi} in Lemma \ref{s2: Edecay}, for $|I|\le N-1$ and $\tau\in[2,s]$, we have
\ben\label{s4: Gppsig0psi}
&&\|\Gamma^{I}\big(\partial_\az\psi^*\gz^0\partial^\az\psi\big)\|_{L^2(\mathscr{H}^{\rm{in}}_\tau)}\lesssim\sum_{|I_1|+|I_2|\le N-1}\|(\partial_\az\Gamma^{I_1}\psi)^*\gz^0\partial^\az\Gamma^{I_2}\psi\|_{L^2(\mathscr{H}^{\rm{in}}_\tau)}\nonumber\\
&\lesssim&\sum_{\substack{|I_1|+|I_2|\le N-1\\1\le a\le 2}}\|t^{-1}|L_a\Gamma^{I_1}\psi|\cdot|\partial\Gamma^{I_2}\psi|+((t-r)/t)|\partial\Gamma^{I_1}\psi|\cdot|\partial\Gamma^{I_2}\psi|\,\|_{L^2(\mathscr{H}^{\rm{in}}_\tau)}\nonumber\\
&\lesssim&\sum_{\substack{|I_1|+|I_2|\le N-1\\|J|=|J'|=1}}\|t^{-1}|\Gamma^J\hat{\Gamma}^{I_1}\psi|\cdot|\Gamma^{J'}\hat{\Gamma}^{I_2}\psi|+|\Gamma^J\hat{\Gamma}^{I_1}\psi|\cdot|\hat{\Gamma}^{I_2}(v\psi)|\,\|_{L^2(\mathscr{H}^{\rm{in}}_\tau)}\nonumber\\
&\lesssim&\sum_{\substack{|I_1|\le N-4\\|J|=1,|I_2|\le N}}\|\tau^{-1}\Gamma^J\hat{\Gamma}^{I_1}\psi\|_{L^\infty(\mathscr{H}^{\rm{in}}_\tau)}\cdot\l\|\f{\tau}{t}\hat{\Gamma}^{I_2}\psi\r\|_{L^2(\mathscr{H}^{\rm{in}}_\tau)}+\sum_{|I_1|+|I_2|+|I_3|\le N}\|\,|\Gamma^{I_1}v|\cdot|\hat{\Gamma}^{I_2}\psi|\cdot|\hat{\Gamma}^{I_3}\psi|\,\|_{L^2(\mathscr{H}^{\rm{in}}_\tau)}\nonumber\\
&\lesssim&(C_2\ez)^2\tau^{-2+\dz},
\een
where we use \eqref{s4: psivL^2}, \eqref{s4: psivL^infty} and \eqref{s4: gpvpsps}. Combining the last two estimates, we obtain
\beq\label{s4: tildeFL^2}
\|\Gamma^{I}\tilde{F}_v\|_{L^2(\mathscr{H}^{\rm{in}}_\tau)}\lesssim (C_2\ez)^2\tau^{-2+\dz},\quad |I|\le N-1, \tau\in[2,s].
\eeq
We observe that the estimate of the second term on the right hand side of \eqref{s4: S3tivN-} was given by \eqref{s3: tivexh}. It follows that
\be
\mathcal{E}^{{\rm{in}},h}_1(s,\Gamma^I\tilde{v})
\lesssim\ez^2+(C_1\ez)^3+(C_2\ez)^2\sup_{\tau\in[2,s]}[\mathcal{E}^{{\rm{in}},h}_1(\tau,\Gamma^{I}\tilde{v})]^{\f{1}{2}},\quad |I|\le N-1,
\ee
which implies
\beq\label{s4: S3imtiv}
\sum_{|I|\le N-1}[\mathcal{E}^{{\rm{in}},h}_1(s,\Gamma^I\tilde{v})]^{\f{1}{2}}\lesssim\ez+(C_2\ez)^{\f{3}{2}}.
\eeq
We recall the definition \eqref{dhs0} and that $\tilde{v}=v-\psi^*\gz^0\psi$. By \eqref{s2: 2.2-42} and \eqref{s2: 2.2-43} in Lemma \ref{sC: 2.2-4}, we have
\ben
&&\sum_{|I|\le N}\|\Gamma^I(\psi^*\gz^0\psi)\|_{L^2(\mathscr{H}^{\rm{in}}_s)}\lesssim\sum_{|I_1|+|I_2|\le N}\|\,|(\hat{\Gamma}^{I_1}\psi)_{-}|\cdot|\hat{\Gamma}^{I_2}\psi|+(s^2/t^2)|\hat{\Gamma}^{I_1}\psi|\cdot|\hat{\Gamma}^{I_2}\psi|\,\|_{L^2(\mathscr{H}^{\rm{in}}_s)}\\
&\lesssim&\sum_{\substack{|I_1|\le N-3\\|I_2|\le N}}\l\|s^\dz\f{t}{s}(\hat{\Gamma}^{I_1}\psi)_{-}\r\|_{L^\infty(\mathscr{H}^{\rm{in}}_s)}\cdot\l\|s^{-\dz}\f{s}{t}\hat{\Gamma}^{I_2}\psi\r\|_{L^2(\mathscr{H}^{\rm{in}}_s)}+\sum_{\substack{|I_2|\le N-3\\|I_1|\le N}}\|s^{-\dz}(\hat{\Gamma}^{I_1}\psi)_{-}\|_{L^2(\mathscr{H}^{\rm{in}}_s)}\cdot\|s^{\dz}\hat{\Gamma}^{I_2}\psi\|_{L^\infty(\mathscr{H}^{\rm{in}}_s)}\\
&+&\sum_{\substack{|I_1|\le N-3\\|I_2|\le N}}\|s^\dz\hat{\Gamma}^{I_1}\psi\|_{L^\infty(\mathscr{H}^{\rm{in}}_s)}\cdot\l\|s^{-\dz}\f{s}{t}\hat{\Gamma}^{I_2}\psi\r\|_{L^2(\mathscr{H}^{\rm{in}}_s)}\lesssim(C_2\ez)^2,
\een
where we use \eqref{s4: psivL^2} and \eqref{s4: psivL^infty}. Hence \eqref{s4: S3imtiv} gives
\beq\label{s4: S3vN-}
\sum_{|I|\le N-1}[\mathcal{E}^{{\rm{in}},h}_1(s,\Gamma^Iv)]^{\f{1}{2}}\lesssim\ez+(C_2\ez)^{\f{3}{2}}.
\eeq
Combining \eqref{s4: repsiN}, \eqref{s4: revN}, \eqref{s4: S2psN-} and \eqref{s4: S3vN-}, we have strictly improved the bootstrap estimates \eqref{s4: bspsiv} and \eqref{s4: bspsivN-3}. Hence the proof of Proposition \ref{s4: maxtildeT} is completed. We conclude that the solution $(\psi,v)$ exists globally in the interior region $\mathscr{D}^{\rm{in}}$ and satisfies \eqref{s4: bspsiv}-\eqref{s4: bspsivN-3} for all $s\in[2,\infty)$.

\section{$n=3$: Uniform in $M$ existence}\label{s5}

In this section, we prove the existence result in Theorem \ref{thm2}, i.e., we prove uniform (in the mass parameter $M\in[0,1]$) global existence with unified time decay estimates for the solution to \eqref{s5: 3DDKG}-\eqref{s5: 3Dinitial}.

\subsection{Bootstrap setting and global existence in the exterior region}\label{s5.1}

Fix $N\ge 5$ and $\lz>0$. Let $C_1\gg 1$ and $0<\ez\ll C_1^{-1}$ be chosen later. We introduce the following bootstrap setting for the solution $(\psi,v)$ to \eqref{s5: 3DDKG}-\eqref{s5: 3Dinitial} in a time interval $[2,T]$:
\beq\label{s5: bspsi}
\sum_{|I|\le N}\big\{\|\hat{\Gamma}^I\psi\|_{Y^{{\rm{ex}},\lz}_{D,t}}+\|\Gamma^Iv\|_{Y^{{\rm{ex}},\lz}_{1,0,t}}\big\}\le C_1\ez,
\eeq
where (see \eqref{s2: exnormD} and \eqref{s2: exnorm})
\beq\label{s5: bspsi'}
\begin{split}
&\|\hat{\Gamma}^I\psi\|_{Y^{{\rm{ex}},\lz}_{D,t}}=[E^{{\rm{ex}},\lz}_D(t,\hat{\Gamma}^I\psi)]^{\f{1}{2}}+\l(\int_2^t\int_{\Sigma^{\rm{ex}}_\tau}(2+r-\tau)^\lz|[\hat{\Gamma}^I\psi]_{-}|^2{\rm{d}}x{\rm{d}}\tau\r)^{\f{1}{2}},\\
&\|\Gamma^Iv\|_{Y^{{\rm{ex}},\lz}_{1,0,t}}=[E^{{\rm{ex}},\lz}_{1}(t,\Gamma^Iv)]^{\f{1}{2}}+\l(\int_2^t\int_{\Sigma^{\rm{ex}}_\tau}(2+r-\tau)^\lz\big(|G\Gamma^Iv|^2+|\Gamma^Iv|^2\big){\rm{d}}x{\rm{d}}\tau\r)^{\f{1}{2}}
\end{split}
\eeq
and we recall \eqref{dexfD} and \eqref{dexfsaz}. In the sequel, the implied constants in $\lesssim$ do not depend on the mass parameter $M\in[0,1]$ and the constants $C_1$ and $\ez$ appearing in the bootstrap assumption \eqref{s5: bspsi}. Let 
\beqn\label{s5: l(tau)}
l(\tau):=(C_1\ez)^{-2}\sum_{|I|\le N}\int_{\Sigma^{\rm{ex}}_\tau}(2+r-\tau)^\lz\big(|[\hat{\Gamma}^I\psi]_{-}|^2+|\Gamma^Iv|^2\big){\rm{d}}x.
\eeqn
Then $l(\tau)\in L^1([2,t])$ and $\|l\|_{L^1([2,t])}\lesssim 1$. Under the assumption \eqref{s5: bspsi}, we have the following $L^2$-type estimates for the solution $(\psi,v)$ on the time interval $[2,T]$:
\beq\label{s5: L^2psiv}
\l\{\begin{split}
&\sum_{|I|\le N}\|(2+r-t)^{\f{1+\lz}{2}}\big(|\hat{\Gamma}^I\psi|+|\partial\Gamma^Iv|+|\Gamma^Iv|\big)\|_{L^2_x(\Sigma^{\rm{ex}}_t)}\lesssim C_1\ez,\\
&\sum_{|I|\le N}\|(2+r-\tau)^{\f{\lz}{2}}\big(|[\hat{\Gamma}^I\psi]_{-}|+|\Gamma^Iv|\big)\|_{L^2_x(\Sigma^{\rm{ex}}_\tau)}\lesssim C_1\ez\sqrt{l(\tau)},\quad\tau\in[2,t].
\end{split}\r.
\eeq
Similar to \eqref{s3: L^2psi-}, by \eqref{s5: L^2psiv}, we have
\beqn\label{s5: L^2psi-}
\sum_{|I|\le N-3,|J|\le 2}\big\{\|(2+r-\tau)^{\f{\lz}{2}}\hat{\Omega}^J[\hat{\Gamma}^I\psi]_{-}\|_{L^2_x(\Sigma^{\rm{ex}}_\tau)}&+&\|(2+r-\tau)^{\f{\lz}{2}}\partial_r\hat{\Omega}^J[\hat{\Gamma}^I\psi]_{-}\|_{L^2_x(\Sigma^{\rm{ex}}_\tau)}\big\}\nonumber\\
&\lesssim& C_1\ez\sqrt{l(\tau)},\quad\tau\in[2,t].
\eeqn
Using \eqref{s5: L^2psiv}, \eqref{s5: L^2psi-} and Lemma \ref{s2: Sobex}, we obtain the following pointwise estimates on the time interval $[2,T]$:
\beq\label{s5: pointpsiv}
\l\{\begin{split}
&\sum_{|I|\le N-3}\sup_{\Sigma^{\rm{ex}}_t}\ (2+r-t)^{\f{1+\lz}{2}} r\big(|\hat{\Gamma}^I\psi|+|\partial\Gamma^Iv|+|\Gamma^Iv|\big)\lesssim C_1\ez,\\
&\sum_{|I|\le N-3}\sup_{\Sigma^{\rm{ex}}_\tau}\ (2+r-\tau)^{\f{\lz}{2}} r\big(|[\hat{\Gamma}^I\psi]_{-}|+|\Gamma^Iv|\big)\lesssim C_1\ez\sqrt{l(\tau)},\quad\tau\in[2,t].
\end{split}\r.
\eeq
Arguing as in \eqref{s3: pointpsiv0*}, we have 
\beq\label{s5: pointpsiv0*}
\sum_{|I|\le N-3}\sup_{\Sigma^{\rm{ex}}_t}(2+r-t) r^{\f{1}{2}+\lz_0}\big(|\hat{\Gamma}^I\psi|+|\partial\Gamma^Iv|+|\Gamma^Iv|\big)\lesssim C_1\ez,
\eeq
where $\lz_0:=\f{1}{2}\min\{\lz,1\}$. Hence similar to \eqref{s3: exvextra*} (see \eqref{s3: pointve*} and \eqref{s3: pointvf}), by \eqref{s5: pointpsiv}, \eqref{s5: pointpsiv0*} and Lemma \ref{s2: Edecay}, we have
\beq\label{s5: pointve*}
|\Gamma^Iv|\lesssim C_1\ez r^{-\f{3}{2}}\quad \mathrm{on}\  \Sigma^{\rm{ex}}_t,\ \mathrm{for}\  |I|\le N-4.
\eeq
By Propositions \ref{s2: exfsD} and \ref{s2: exfs}, for $|I|\le N$, we have
\be
\begin{split}
\|\hat{\Gamma}^I\psi\|_{Y^{{\rm{ex}},\lz}_{D,t}}&\lesssim\|\hat{\Gamma}^I\psi\|_{Y^{{\rm{ex}},\lz}_{D,2}}+\int_2^t\|(2+r-\tau)^{\f{1+\lz}{2}}\hat{\Gamma}^I(v\psi)\|_{L^2_x(\Sigma^{\rm{ex}}_\tau)}{\rm{d}}\tau,\\
\|\Gamma^Iv\|_{Y^{{\rm{ex}},\lz}_{1,0,t}}&\lesssim\|\Gamma^Iv\|_{Y^{{\rm{ex}},\lz}_{1,0,2}}+\int_2^t\|(2+r-\tau)^{\f{1+\lz}{2}}\Gamma^I(\psi^*\gz^0\psi)\|_{L^2_x(\Sigma^{\rm{ex}}_\tau)}{\rm{d}}\tau.
\end{split}
\ee
For $|I|\le N$ and $\tau\in[2,t]$, by \eqref{s2: 2.2-41}, \eqref{s2: 2.2-42} and \eqref{s2: 2.2-43'} in Lemma \ref{sC: 2.2-4}, we have
\beqn\label{s5: psiNL^2}
\!\!\!\!\!\!\!\!\!&&\|(2+r-\tau)^{\f{1+\lz}{2}}\big(|\hat{\Gamma}^I(v\psi)|+|\Gamma^I(\psi^*\gz^0\psi)|\big)\|_{L^2_x(\Sigma^{\rm{ex}}_\tau)}\nonumber\\
\!\!\!\!\!\!\!\!\!&\lesssim&\sum_{|I_1|+|I_2|\le N}\|(2+r-\tau)^{\f{1+\lz}{2}}\big(|\Gamma^{I_1}v|+|[\hat{\Gamma}^{I_1}\psi]_{-}|\big)\cdot|\hat{\Gamma}^{I_2}\psi|\,\|_{L^2_x(\Sigma^{\rm{ex}}_\tau)}\nonumber\\
\!\!\!\!\!\!\!\!\!&\lesssim&\sum_{\substack{|I_1|\le N-3\\|I_2|\le N}}\|\,|\Gamma^{I_1}v|+|[\hat{\Gamma}^{I_1}\psi]_{-}|\,\|_{L^\infty_x(\Sigma^{\rm{ex}}_\tau)}\cdot\|(2+r-\tau)^{\f{1+\lz}{2}}\hat{\Gamma}^{I_2}\psi\|_{L^2_x(\Sigma^{\rm{ex}}_\tau)}\nonumber\\
\!\!\!\!\!\!\!\!\!&+&\sum_{\substack{|I_2|\le N-3\\|I_1|\le N}}\|\,|\Gamma^{I_1}v|+|[\hat{\Gamma}^{I_1}\psi]_{-}|\,\|_{L^2_x(\Sigma^{\rm{ex}}_\tau)}\cdot\|(2+r-\tau)^{\f{1+\lz}{2}}\hat{\Gamma}^{I_2}\psi\|_{L^\infty_x(\Sigma^{\rm{ex}}_\tau)}\lesssim(C_1\ez)^2\sqrt{l(\tau)}\tau^{-1},
\eeqn
where we use \eqref{s5: L^2psiv} and \eqref{s5: pointpsiv}. It follows that
\beq\label{s5: S1psi}
\sum_{|I|\le N}\big\{\|\hat{\Gamma}^I\psi\|_{Y^{{\rm{ex}},\lz}_{D,t}}+\|\Gamma^Iv\|_{Y^{{\rm{ex}},\lz}_{1,0,t}}\big\}\lesssim\ez+(C_1\ez)^2.
\eeq
Hence we have strictly improved the bootstrap estimate \eqref{s5: bspsi} for $C_1$ sufficiently large and $0<\ez\ll C_1^{-1}$ sufficiently small. We conclude that the solution $(\psi,v)$ exists globally in the exterior region and satisfies \eqref{s5: bspsi} for all $t\in[2,\infty)$. Moreover, by Proposition \ref{s2: exhsD} and Proposition \ref{s2: exhs} (with $\dz=0$), we obtain the energy estimates of the solution $(\psi,v)$ on truncated exterior hyperboloids $\mathscr{H}^{\rm{ex}}_s$ for any $s\in[2,\infty)$:
\beq\label{s5: psiexh}
\l\{\begin{split}
&\mathcal{E}^{{\rm{ex}},h}_D(s,\hat{\Gamma}^I\psi)+\int_{\mathscr{l}_{[2,\,t(s)]}}|[\hat{\Gamma}^I\psi]_{-}|^2{\rm{d}}\sz\lesssim\ez^2+(C_1\ez)^3,\quad |I|\le N,\\
&\mathcal{E}^{{\rm{ex}},h}_{1}(s,\Gamma^Iv)+\int_{\mathscr{l}_{[2,\,t(s)]}}\big(|G\Gamma^Iv|^2+|\Gamma^Iv|^2\big){\rm{d}}\sz\lesssim\ez^2+(C_1\ez)^3,\quad |I|\le N,
\end{split}\r.
\eeq
where $t(s)=\f{s^2+1}{2}$ and we use \eqref{s5: psiNL^2} and \eqref{s5: S1psi}. Below we give a mass-dependent estimate for $\psi$ in the exterior region. Acting the Dirac operator on both sides of the first equation in \eqref{s5: 3DDKG}, we obtain
\be
-\Box\psi+M^2\psi=-i\gz^\mu\partial_\mu(v\psi)+Mv\psi.
\ee
By Lemma \ref{s2: Edecay}, \eqref{s5: pointpsiv} and \eqref{s5: pointpsiv0*}, for any $t\in[2,\infty)$, on $\Sigma^{\rm{ex}}_t\cap\{r\le 3t\}$ we have
\beq\label{s5: pointpsie}
\sum_{|I|\le N-5}M^2|\Gamma^I\psi|\lesssim\f{2+r-t}{t}\sum_{|I|\le N-4}|\partial\Gamma^I\psi|+\sum_{|I|\le N-5}|\Gamma^I\big(-i\gz^\mu\partial_\mu(v\psi)+Mv\psi\big)|\lesssim C_1\ez r^{-\f{3}{2}}
\eeq
for any $M\in[0,1]$. By \eqref{s5: pointpsiv}, this also holds on $\Sigma^{\rm{ex}}_t\cap\{r\ge 3t\}$. Hence, we obtain the following estimate in the region $\mathscr{D}^{\rm{ex}}$:
\beq\label{s5: **pointpsie}
M^2|\Gamma^I\psi|\lesssim C_1\ez r^{-\f{3}{2}},\quad |I|\le N-5.
\eeq

\subsection{Bootstrap setting and global existence in the interior region}

Let $N\ge 5$, $C_1\gg 1$ be as in Section \ref{s5.1}. Let $C_2\gg C_1$ and $0<\ez\ll C_2^{-1}$ be chosen later. We assume the exterior estimate \eqref{s5: bspsi} globally in time as well as the following interior bounds for the solution $(\psi,v)$ to \eqref{s5: 3DDKG}-\eqref{s5: 3Dinitial}, for $s\in[2,\tilde{T}]$:
\beq\label{s52: bspsiv}
\sum_{|I|\le N} \big\{[\mathcal{E}^{{\rm{in}},h}_D(s,\hat{\Gamma}^I\psi)]^{\f{1}{2}}+[\mathcal{E}^{{\rm{in}},h}_1(s,\Gamma^Iv)]^{\f{1}{2}}\big\}\le C_2\ez,
\eeq
where we recall \eqref{dinhsD} and \eqref{dhs0}. In the sequel, the implied constants in $\lesssim$ do not depend on the mass parameter $M\in[0,1]$ and the constants $C_2$ and $\ez$ appearing in the bootstrap assumption \eqref{s52: bspsiv}. 

By \eqref{s5: psiexh} and \eqref{s52: bspsiv}, the following $L^2$-type estimates hold for $s\in[2,\tilde{T}]$:
\beq\label{s52: psivL^2}
\sum_{|I|\le N}\|\,|(\hat{\Gamma}^I\psi)_{-}|+|(s/t)\hat{\Gamma}^I\psi|+|(s/t)\partial\Gamma^Iv|+|\Gamma^Iv|\,\|_{L^2(\mathscr{H}_s)}\lesssim C_2\ez.
\eeq
By Lemmas \ref{s4: Sobin} and \ref{sC: DL2.3'}, we obtain the following pointwise estimates for $s\in[2,\tilde{T}]$:
\beq\label{s52: psivL^infty}
\sum_{|I|\le N-2}\|t^{\f{3}{2}}\big\{|(\hat{\Gamma}^I\psi)_{-}|+|(s/t)\hat{\Gamma}^I\psi|+|(s/t)\partial\Gamma^Iv|+|\Gamma^Iv|\big\}\|_{L^\infty(\mathscr{H}^{\rm{in}}_s)}\lesssim C_2\ez.
\eeq
By Propositions \ref{s2: inhsD} and \ref{s2: inhs}, for $|I|\le N$, we have
\be\label{s5: inhpsiv}
\begin{split}
&\mathcal{E}^{{\rm{in}},h}_D(s,\hat{\Gamma}^I\psi)\lesssim \mathcal{E}^{{\rm{in}},h}_D(2,\hat{\Gamma}^I\psi)+\mathcal{R}^I(s)+\int_2^s[\mathcal{E}^{{\rm{in}},h}_D(\tau,\hat{\Gamma}^I\psi)]^{\f{1}{2}}\cdot\|\hat{\Gamma}^I(v\psi)\|_{L^2(\mathscr{H}^{\rm{in}}_\tau)}{\rm{d}}\tau,\\
&\mathcal{E}^{{\rm{in}},h}_1(s,\Gamma^Iv)\lesssim \mathcal{E}^{{\rm{in}},h}_1(2,\Gamma^Iv)+\mathcal{\tilde{R}}^I(s)+\int_2^s[\mathcal{E}^{{\rm{in}},h}_1(\tau,\Gamma^Iv)]^{\f{1}{2}}\cdot\|\Gamma^I(\psi^*\gz^0\psi)\|_{L^2(\mathscr{H}^{\rm{in}}_\tau)}{\rm{d}}\tau,
\end{split}
\ee
where
\beq\label{s5: inhpsiv'}
\mathcal{R}^I(s):=\int_{\mathscr{l}_{[t(2),\,t(s)]}}|[\hat{\Gamma}^I\psi]_{-}|^2{\rm{d}}\sz,\quad\quad\mathcal{\tilde{R}}^I(s):=\int_{\mathscr{l}_{[t(2),\,t(s)]}}\big(|G\Gamma^Iv|^2+|\Gamma^Iv|^2\big){\rm{d}}\sz
\eeq
and $t(s), t(2)$ are as in \eqref{t(s)t(2)}. For $|I|\le N$ and $\tau\in[2,s]$, by \eqref{s2: 2.2-41}, \eqref{s2: 2.2-42} and \eqref{s2: 2.2-43} in Lemma \ref{sC: 2.2-4}, \eqref{s52: psivL^2} and \eqref{s52: psivL^infty}, we have
\beqn\label{s5: inpsiNL^2}
&&\|\hat{\Gamma}^I(v\psi)\|_{L^2(\mathscr{H}^{\rm{in}}_\tau)}\lesssim\sum_{|I_1|+|I_2|\le N}\|\,|\Gamma^{I_1}v|\cdot|\hat{\Gamma}^{I_2}\psi|\,\|_{L^2(\mathscr{H}^{\rm{in}}_\tau)}\nonumber\\
&\lesssim&\sum_{\substack{|I_1|\le N-2\\|I_2|\le N}}\|(t/\tau)\Gamma^{I_1}v\|_{L^\infty(\mathscr{H}^{\rm{in}}_\tau)}\cdot\|(\tau/t)\hat{\Gamma}^{I_2}\psi\|_{L^2(\mathscr{H}^{\rm{in}}_\tau)}+\sum_{\substack{|I_2|\le N-2\\|I_1|\le N}}\|\Gamma^{I_1}v\|_{L^2(\mathscr{H}^{\rm{in}}_\tau)}\cdot\|\hat{\Gamma}^{I_2}\psi\|_{L^\infty(\mathscr{H}^{\rm{in}}_\tau)}\nonumber\\
&\lesssim&(C_2\ez)^2\tau^{-\f{3}{2}}
\eeqn
and
\beqn\label{s5: invNL^2}
&&\|\Gamma^I(\psi^*\gz^0\psi)\|_{L^2(\mathscr{H}^{\rm{in}}_\tau)}\lesssim\sum_{|I_1|+|I_2|\le N}\|\,|(\hat{\Gamma}^{I_1}\psi)_{-}|\cdot|\hat{\Gamma}^{I_2}\psi|+(\tau^2/t^2)|\hat{\Gamma}^{I_1}\psi|\cdot|\hat{\Gamma}^{I_2}\psi|\,\|_{L^2(\mathscr{H}^{\rm{in}}_\tau)}\nonumber\\
&\lesssim&\sum_{\substack{|I_1|\le N-2\\|I_2|\le N}}\|(t/\tau)(\hat{\Gamma}^{I_1}\psi)_{-}\|_{L^\infty(\mathscr{H}^{\rm{in}}_\tau)}\cdot\|(\tau/t)\hat{\Gamma}^{I_2}\psi\|_{L^2(\mathscr{H}^{\rm{in}}_\tau)}+\sum_{\substack{|I_2|\le N-2\\|I_1|\le N}}\|(\hat{\Gamma}^{I_1}\psi)_{-}\|_{L^2(\mathscr{H}^{\rm{in}}_\tau)}\cdot\|\hat{\Gamma}^{I_2}\psi\|_{L^\infty(\mathscr{H}^{\rm{in}}_\tau)}\nonumber\\
&+&\sum_{\substack{|I_1|\le N-2\\|I_2|\le N}}\|\hat{\Gamma}^{I_1}\psi\|_{L^\infty(\mathscr{H}^{\rm{in}}_\tau)}\cdot\|(\tau/t)\hat{\Gamma}^{I_2}\psi\|_{L^2(\mathscr{H}^{\rm{in}}_\tau)}\lesssim (C_2\ez)^2\tau^{-\f{3}{2}}.
\eeqn
We note that the estimates of $\mathcal{R}^I(s)$ and $\mathcal{\tilde{R}}^I(s)$ in \eqref{s5: inhpsiv'} were given by \eqref{s5: psiexh}. It follows that
\be\label{s52: repsiN}
\sum_{|I|\le N}\big\{[\mathcal{E}^{{\rm{in}},h}_D(s,\hat{\Gamma}^I\psi)]^{\f{1}{2}}+[\mathcal{E}^{{\rm{in}},h}_1(s,\Gamma^Iv)]^{\f{1}{2}}\big\}\lesssim\ez+(C_2\ez)^{\f{3}{2}}.
\ee
Hence we have strictly improved the bootstrap estimate \eqref{s52: bspsiv} for $C_2$ sufficiently large and $0<\ez\ll C_2^{-1}$ sufficiently small. We conclude that the solution $(\psi,v)$ exists globally in the interior region and satisfies \eqref{s52: bspsiv} for all $s\in[2,\infty)$. We next derive a mass-dependent  estimate for $\psi$ in the interior region. By \eqref{s5: 3DDKG}, for any $|I|\le N-3$, we have
\be
i\gz^\mu\partial_\mu\hat{\Gamma}^I\psi+M\hat{\Gamma}^I\psi=\hat{\Gamma}^I(v\psi).
\ee
We rewrite the above equation as
\ben
i\l(\gz^0-\f{x_a}{t}\gz^a\r)\partial_t\hat{\Gamma}^I\psi+it^{-1}\gz^aL_a\hat{\Gamma}^I\psi+M\hat{\Gamma}^I\psi=\hat{\Gamma}^I(v\psi).
\een
By \eqref{s52: psivL^infty}, we obtain the following estimate in the region $\mathscr{D}^{\rm{in}}$:
\beq\label{s5: inpsiextra}
M|\hat{\Gamma}^I\psi|\lesssim|(\partial_t\hat{\Gamma}^I\psi)_{-}|+t^{-1}|L_a\hat{\Gamma}^I\psi|+|\hat{\Gamma}^I(v\psi)|\lesssim C_2\ez t^{-\f{3}{2}}, \quad |I|\le N-3.
\eeq

\section{Scattering results}\label{s6}

\subsection{Scattering for the $3D$ case}

We first prove the scattering result for the solution to \eqref{s5: 3DDKG}-\eqref{s5: 3Dinitial} in Theorem \ref{thm2}

By \eqref{s5: psiNL^2}, \eqref{s5: inpsiNL^2} and \eqref{s5: invNL^2}, for any $\dz>0$ sufficiently small and any $4\le T_1<\infty$, we have
\ben
\l(\int_{T_1}^{+\infty}\|F\|_{L^2_x(\Sigma^{\rm{ex}}_\tau)}^2 \cdot\tau^{1+\dz}{\rm{d}}\tau+\int_{T_1^{\f{1}{2}}}^{+\infty}\|F\|_{L^2(\mathscr{H}^{\rm{in}}_s)}^2\cdot s^{1+2\dz}{\rm{d}}s\r)^{\f{1}{2}}\lesssim T_1^{-\f{1}{4}+\f{\dz}{2}},
\een
where $F:=\sum_{k=0}^{N}\big(|\nabla^{k}F_\psi|+|\nabla^kF_v|\big)$ and $\nabla=(\partial_1,\partial_2,\partial_3)$. By Lemma \ref{Sca}, the solution $(\psi,v)$ scatters to a free solution in $H^{N}(\mathbb{R}^3)\times H^{N+1}(\mathbb{R}^3)\times H^{N}(\mathbb{R}^3)$.

\subsection{Scattering for the $2D$ case}

We next prove the scattering result for the solution to \eqref{s1: DKG}-\eqref{s1: initial} in Theorem \ref{thm1}. Recall that by Lemma \ref{s2: nlt}, the functions $\tilde{\psi}=\psi-i\gz^\mu\partial_\mu(v\psi)$ and $\tilde{v}=v-\psi^*\gz^0\psi$ solve the equations
\be
i\gz^\mu\partial_\mu\tilde{\psi}=\tilde{F}_\psi,\quad\quad-\Box\tilde{v}+\tilde{v}=\tilde{F}_v,
\ee
where $\tilde{F}_\psi$ and $\tilde{F}_v$ are as in \eqref{s3: tildeG} and \eqref{s3: tildeF} respectively. By \eqref{s3: stildeG}, \eqref{s3: stildevF}, \eqref{s4: tildeGL^2} and \eqref{s4: tildeFL^2}, we have
\be
\|\tilde{F}\|_{L^2_x(\Sigma^{\rm{ex}}_\tau)}\lesssim \sqrt{l(\tau)}\tau^{-1+\f{\dz}{2}}+\tau^{-1-\lz_0+\f{\dz}{2}}\quad \quad\mathrm{and}\quad\quad\|\tilde{F}\|_{L^2(\mathscr{H}^{\rm{in}}_\tau)}\lesssim\tau^{-2+\dz}
\ee
for any $\tau\in[2,\infty)$, where $\tilde{F}:=\sum_{k=0}^{N-1}\big(|\nabla^{k}\tilde{F}_\psi|+|\nabla^k\tilde{F}_v|\big)$, $\nabla=(\partial_1,\partial_2)$, $\lz_0=\f{1}{2}\min\{\lz,1\}$ is as in \eqref{s3: deflz_0} and $0<\dz\ll\lz_0$. Hence for any $4\le T_1<\infty$, we have
\be
\l(\int_{T_1}^{+\infty}\|\tilde{F}\|_{L^2_x(\Sigma^{\rm{ex}}_\tau)}^2 \cdot\tau^{1+\dz}{\rm{d}}\tau+\int_{T_1^{\f{1}{2}}}^{+\infty}\|\tilde{F}\|_{L^2(\mathscr{H}^{\rm{in}}_s)}^2\cdot s^{1+2\dz}{\rm{d}}s\r)^{\f{1}{2}}\lesssim T_1^{-\lz_0+\dz}.
\ee
By Lemma \ref{Sca}, $(\tilde{\psi},\tilde{v})$ scatters to a free solution in $\mathcal{X}_{N-1}(\mathbb{R}^2):=H^{N-1}(\mathbb{R}^2)\times H^{N}(\mathbb{R}^2)\times H^{N-1}(\mathbb{R}^2)$. We claim that
\beq\label{s6: psivsca}
\lim_{t\to+\infty}\|(\psi,v,\partial_tv)-(\tilde{\psi},\tilde{v},\partial_t\tilde{v})\|_{\mathcal{X}_{N-1}(\mathbb{R}^2)}=0.
\eeq
This implies the scattering result in Theorem \ref{thm1}. Next we prove \eqref{s6: psivsca}. We recall \eqref{s3: vNL^2}, \eqref{s4: higpsi} and \eqref{s4: higv}, and note that a similar estimate as in \eqref{s3: vNL^2} holds for $\|(2+r-\tau)^{\f{1+\lz}{2}}\hat{\Gamma}^I(v\psi)\|_{{L^2_x(\Sigma^{\rm{ex}}_\tau)}}$, where $|I|\le N$, $\tau\in[2,\infty)$. Hence, we obtain
\beq\label{s6: FL2}
\|F\|_{L^2_x(\Sigma^{\rm{ex}}_\tau)}\lesssim \sqrt{l(\tau)}\tau^{-\f{1}{2}+\f{\dz}{2}}\quad\quad\mathrm{and}\quad\quad \|F\|_{L^2(\mathscr{H}^{\rm{in}}_\tau)}\lesssim\tau^{-1+\dz}
\eeq
for any $\tau\in[2,\infty)$, where we denote $F:=\sum_{k=0}^N \big(|\nabla^{k}F_\psi|+|\nabla^kF_v|\big)$. By \eqref{s6: FL2} and the proof of \eqref{s5: II}, for $4\le t<\infty$, we have
\be
\begin{split}
&\quad\ \int_4^t\|F(\tau)\|_{L^2_x(\mathbb{R}^2)}\cdot\tau^{\f{1+\dz}{2}}\cdot\tau^{-\f{1+\dz}{2}}{\rm{d}}\tau\lesssim\l(\int_4^t\|F(\tau)\|_{L^2_x(\mathbb{R}^2)}^2\cdot\tau^{1+\dz}{\rm{d}}\tau\r)^{\f{1}{2}}\\
&\lesssim\l(\int_4^t\l(\int_{r\ge\tau-1}+\int_{r<\tau-1}\r)|F|^2(\tau,x)\cdot\tau^{1+\dz}{\rm{d}}x{\rm{d}}\tau\r)^{\f{1}{2}}\\
&\lesssim\l(\int_{4}^{t}\|F\|_{L^2_x(\Sigma^{\rm{ex}}_\tau)}^2 \cdot\tau^{1+\dz}{\rm{d}}\tau\r)^{\f{1}{2}}+\l(\int_{2}^{t}\|F\|_{L^2(\mathscr{H}^{\rm{in}}_s)}^2\cdot s^{1+2\dz}{\rm{d}}s\r)^{\f{1}{2}}\lesssim t^{2\dz}.
\end{split}
\ee
Hence, by Remarks \ref{RemK} and \ref{RemD}, for any $t\ge 4$, we have
\beq\label{s6: L^2psiv}
\sum_{k=0}^N\|\,|\nabla^k\psi(t)|+|\partial\nabla^kv(t)|+|\nabla^kv(t)|\,\|_{L^2_x(\mathbb{R}^2)}\lesssim\ez+\int_2^t\|F(\tau)\|_{L^2_x(\mathbb{R}^2)}{\rm{d}}\tau\lesssim t^{2\dz}.
\eeq
Using \eqref{s6: L^2psiv} and the pointwise estimates \eqref{s3: pointpsiv} and \eqref{s4: psivL^infty} in Sections \ref{s3} - \ref{s4}, we have
\ben
&&\sum_{k=0}^{N-1}\|\nabla^k\big(i\gz^\mu\partial_\mu(v\psi)\big)\|_{L^2_x(\mathbb{R}^2)}\lesssim\sum_{k_1+k_2\le N-1}\|\,|\partial\nabla^{k_1}v|\cdot|\nabla^{k_2}\psi|+|\nabla^{k_1}v|\cdot|\nabla^{k_2}(v\psi)|\,\|_{L^2_x(\mathbb{R}^2)}\\
&\lesssim&\sum_{\substack{k_1\le N-3\\k_2\le N-1}}\|t^{2\dz}\partial\nabla^{k_1}v\|_{L^\infty_x(\mathbb{R}^2)}\cdot\|t^{-2\dz}\nabla^{k_2}\psi\|_{L^2_x(\mathbb{R}^2)}+\sum_{\substack{k_2\le N-3\\k_1\le N-1}}\|t^{-2\dz}\partial\nabla^{k_1}v\|_{L^2_x(\mathbb{R}^2)}\cdot\|t^{2\dz}\nabla^{k_2}\psi\|_{L^\infty_x(\mathbb{R}^2)}\\
&+&\sum_{\substack{k_1\le N-3\\k_2,k_3\le N-1}}\|t^{4\dz}\nabla^{k_1}v\|_{L^\infty_x(\mathbb{R}^2)}\cdot\|t^{-2\dz}\nabla^{k_2}v\|_{L^2_x(\mathbb{R}^2)}\cdot\|t^{-2\dz}\nabla^{k_3}\psi\|_{L^2_x(\mathbb{R}^2)}\\
&+&\sum_{\substack{k_3\le N-3\\k_1,k_2\le N-1}}\|t^{-2\dz}\nabla^{k_1}v\|_{L^2_x(\mathbb{R}^2)}\cdot\|t^{-2\dz}\nabla^{k_2}v\|_{L^2_x(\mathbb{R}^2)}\cdot\|t^{4\dz}\nabla^{k_3}\psi\|_{L^\infty_x(\mathbb{R}^2)}\lesssim t^{-\f{1}{2}+4\dz}
\een
and similarly for $\sum_{k=0}^N\|\nabla^k(\psi^*\gz^0\psi)\|_{L^2_x(\mathbb{R}^2)}+\sum_{k=0}^{N-1}\|\nabla^k\partial_t(\psi^*\gz^0\psi)\|_{L^2_x(\mathbb{R}^2)}$. Hence \eqref{s6: psivsca} holds, which implies that the solution $(\psi,v)$ scatters to a free solution in $\mathcal{X}_{N-1}(\mathbb{R}^2)$. $\hfill\square$

Qian Zhang

School of Mathematical Sciences, Beijing Normal University

Laboratory of Mathematics and Complex Systems, Ministry of Education Beijing 100875, China

Email: qianzhang@bnu.edu.cn

\end{document}